\documentclass[a4paper,11pt]{article} 
 
\usepackage{mathrsfs}
\usepackage{amsmath} 
\usepackage[applemac]{inputenc}
\usepackage{amsfonts}
\usepackage{amssymb}
\usepackage{amsthm}
\usepackage{stmaryrd}
\usepackage{color}
\def\cC{{\mathcal{C}}}
\def\w{{\mathbf{w}}}

\newcommand{\F}{{\bf F}}  
\renewcommand{\P}{\mathbb{P}}

\newcommand{\N}{\mathbb{N}}  
  
\newcommand{\T}{\mathbb{T}}

\newcommand{\R}{\mathbb{R}}

\def\bphi{{\boldsymbol{\phi}}}
\def\bpsi{{\boldsymbol{\psi}}}

\newtheorem{REM}{Remark}

\def\1{\mbox{1\hspace{-0.25em}l}}

\def\mI{{{\mathbb I}}}
\def\cA{{{\mathcal A}}}
\def\cL{{{\mathcal L}}}

\def \t {\tau}
\def \dif {{\rm d}}
\def\tr{\mathrm {tr}}

\def\btheta{{\boldsymbol{\theta}}}
\def\bTheta{{\boldsymbol{\Sigma}}}
\def\bvtheta{{\boldsymbol{\vartheta}}}
\def\bxi{{\boldsymbol{\xi}}}

\def \cA{\mathcal A}

\def\X{{\bf{X}}}
\def\bx{{\mathbf{x}}}
\def\x{{\boldsymbol{x}}}
\def\p{{\mathbf{p}}}

\def\y{{\boldsymbol{y}}}
\def\z{{\boldsymbol{z}}}
\def\l{\lambda}
\def\<{\langle}
\def\>{\rangle}

\numberwithin{equation}{section}

\newcommand{\be}{\begin{eqnarray}}
\newcommand{\ee}{\end{eqnarray}}
\newcommand{\ce}{\begin{eqnarray*}}
\newcommand{\de}{\end{eqnarray*}}
\newtheorem{theorem}{Theorem}[section]
\newtheorem{lemma}[theorem]{Lemma}
\newtheorem{remark}[theorem]{Remark}
\newtheorem{definition}[theorem]{Definition}
\newtheorem{proposition}[theorem]{Proposition}
\newtheorem{Examples}[theorem]{Example}
\newtheorem{corollary}[theorem]{Corollary}

\def\eps{\varepsilon}
\def\t{\tau}
\def\e{\mathrm{e}}

\def\p{\partial}

\def\l{\lambda}

\def\[{{\Big[}}
\def\]{{\Big]}}
\def\<{{\langle}}
\def\>{{\rangle}}
\def\({{\big(}}
\def\){{\big)}}

\def\by{{\mathbf{y}}}
\def\bx{{\mathbf{x}}}
\def\tr{\mathrm {tr}}

\def\dif{{\mathord{{\rm d}}}}

\def\bb2{{\boldsymbol{2}}}
\def\no{\nonumber}
\def\={&\!\!=\!\!&}

\def\cA{{\mathcal A}}

\def\cC{{\mathcal C}}

\def\cE{{\mathcal E}}
\def\cF{{\mathcal F}}

\def\cH{{\mathcal H}}

\def\cJ{{\mathcal J}}

\def\cL{{\mathcal L}}

\def\cT{{\mathcal T}}

\def\mE{{\mathbb E}}

\def\mI{{\mathbb I}}

\def\mN{{\mathbb N}}

\def\mR{{\mathbb R}}

\def\1{{\mathbf{1}}}

\def\sA{{\mathscr A}}
\def\sB{{\mathscr B}}
\def\sC{{\mathscr C}}
\def\sD{{\mathscr D}}

\def\sI{{\mathscr I}}
\def\sJ{{\mathscr J}}
\def\sK{{\mathscr K}}

\def\sN{{\mathscr N}}

\def\geq{\geqslant}
\def\leq{\leqslant}
\def\ge{\geqslant}
\def\le{\leqslant}

\def\div{\mathord{{\rm div}}}

\def\eps{\varepsilon}
\def\t{\tau}
\def\e{\mathrm{e}}

\def\p{\partial}

\def\l{\lambda}

\def\[{{\Big[}}
\def\]{{\Big]}}
\def\<{{\langle}}
\def\>{{\rangle}}

\def\by{{\mathbf{y}}}
\def\bx{{\mathbf{x}}}
\def\tr{\mathrm {tr}}

\def\dif{{\mathord{{\rm d}}}}

\def\no{\nonumber}
\def\={&\!\!=\!\!&}
\def\bt{\begin{theorem}}
\def\et{\end{theorem}}
\def\bl{\begin{lemma}}
\def\el{\end{lemma}}
\def\br{\begin{remark}}
\def\er{\end{remark}}
\def\bx{\begin{Examples}}
\def\ex{\end{Examples}}
\def\bd{\begin{definition}}
\def\ed{\end{definition}}
\def\bp{\begin{proposition}}
\def\ep{\end{proposition}}
\def\bc{\begin{corollary}}
\def\ec{\end{corollary}}

\def\geq{\geqslant}
\def\leq{\leqslant}
\def\ge{\geqslant}
\def\le{\leqslant}

\def\div{\mathord{{\rm div}}}

\def\bA{{\mathbf A}}
 \def\R{\mathbb R}
 \def\R{\mathbb R}    
\def\N{\mathbb N}  
   
\def\<{\langle} \def\>{\rangle}
\def\wt{\widetilde}

\def\0{{\mathbf{0}}}

\def\gF{{\bf F}}

\def\K{{\mathbf{K}}}
\def\gR{{\mathbf{R}}}

\newcommand{\green}[1]{\textcolor{black}{#1}}

\allowdisplaybreaks

\title{\textbf{Heat kernel and gradient estimates for kinetic SDEs with low regularity coefficients}}
\author{\textbf{P.E. Chaudru de Raynal}\footnote{Laboratoire de Math\'ematiques Jean Leray, 2, rue de la Houssini\`ere BP 92208
F-44322 Nantes Cedex 3. 
Paul-Eric.Deraynal@univ-nantes.fr}, \textbf{S. Menozzi}\footnote{Laboratoire de Mod\'elisation Math\'ematique d'Evry (LaMME), Universit\'e d'Evry Val d'Essonne, 23 Boulevard de France 91037 Evry, France and Laboratory of Stochastic Analysis, HSE,
Pokrovsky Blvd, 11, Moscow, Russian Federation. stephane.menozzi@univ-evry.fr}, \textbf{A. Pesce}\footnote{Dipartimento di Matematica, Piazza di Porta San Donato, 5 Bologna (Italy), antonello.pesce2@unibo.it}, \textbf{X. Zhang}\footnote{School of Mathematics and Statistics, Wuhan University, Wuhan, Hubei 430072, P.R. China, Email: XichengZhang@gmail.com}}

\begin{document}
\maketitle
\begin{abstract}
We establish heat kernel and gradient estimates for the density of kinetic degenerate Kolmogorov stochastic differential equations. Our results are established under somehow minimal assumptions that guarantee the SDE is weakly well posed.
\end{abstract}



\section{Introduction}
\subsection{Statement of the problem}

We are interested in providing Aronson-like bounds and pointwise estimates for the full gradient of the transition probability density 
of the following {\it kinetic} system of SDEs:
\begin{equation}
\label{SDE}
\left\{
\begin{aligned}
\dif X_t^1 &= F_1(t,X_t^1,X_t^2) \dif t + \sigma(t,X_{t}^1,X_{t}^2) \dif W_t,\\
\dif X_t^2 &= F_2(t,X_t^1,X_t^2) \dif t, 
\end{aligned}
\right.
\end{equation}
where $(W_t)_{t \geq 0}$ stands for a $d$-dimensional Brownian motion on \textcolor{black}{some} stochastic basis  
$(\Omega,\cF,(\cF_t)_{t\ge 0},\P)$ and for $i\in \{ 1,2\}$, $t\ge 0 $ the component $X_t^i $ is $\R^d $-valued. This  equivalently amounts to establish  
the announced estimates for \green{the} fundamental solution of the parabolic PDE associated with \eqref{SDE}, which writes:
\begin{align}
\begin{cases}
\displaystyle \partial_s p\big(s,\x;t,\y\big)+\left\langle F_1(t,\x) , \nabla_{x_1} p\big(s,\x;t,\y\big)\right\rangle +\left\langle F_2(t,\x) , \nabla_{x_2} p\big(s,\x;t,\y\big)\right\rangle \\
\displaystyle\hspace*{.2cm}+\frac 12{\rm{Tr}}\Big(\sigma\sigma^*(s,\x)\nabla_{x_1}^2 p\big(s,\x;t,\y\big)\Big)=0, 0\le s<t,\ \x=(x_1,x_2),\y=(y_1,y_2)\in \R^{2d},\\
\displaystyle p\big(s,\cdot;t,\textcolor{black}{\y}\big) \underset{s\textcolor{black}{\uparrow} t}{\longrightarrow} \delta_{\textcolor{black}{\y}}(\cdot).
\end{cases}
\label{PDE_FUNDA}
\end{align}

Importantly, we aim at obtaining such estimates under somehow \textit{minimal} assumptions, relying on \textit{minimal} conditions required on the SDE to be well posed, in a weak sense. For our approach to work, we will assume a  kind of {\it weak} H\"ormander condition - the Jacobian $(\nabla_{x_1}F_2)$ has full rank and the diffusion coefficient $\sigma$ is bounded and separated from zero - however
 the coefficients can be rather rough in their entries, namely, measurable with respect to the time variable and belonging to suitable \textit{anisotropic H\"older spaces} in the spatial variables. \textcolor{black}{In particular, we emphasize that the drift term $F$ can be unbounded in all its variables and entries.} Through the analysis, some thresholds for the H\"older regularity of the drift with respect to the second (whence degenerate) variable $x_2$ will appear. Such thresholds are related to the degenerate nature of the system of interest and appear to be rather sharp as they are precisely the ones which provide a sufficient and (almost) necessary condition for the system to be well posed{\color{black}, see \cite{Raynal2017RegularizationEO}}.\\

From now on we shall use bold letters $\X$ and $\F$ to denote vectors $(X^1,X^2)$ and $(F_1,F_2)$ in $\mR^{2d}$. 
Let $B=(\mI_{d\times d},0_{d\times d})^*$ be a $2d\times d$-matrix, where $*$ stands for the transpose. 
Using these notations,  we can rewrite SDE \eqref{SDE} in the following compact form:
\begin{equation}\label{SDE0}
\dif \X_t = \F(t,\X_t) \dif t + B\sigma(t,\X_{t}) \dif W_t.
\end{equation}
\smallskip

\noindent\textbf{Related applications}:
 These kinds of \textit{kinetic} (or speed/position) systems appear in several application fields. For instance \eqref{SDE} describes the dynamics of some Hamiltonian systems. For a Hamilton function $H(\x)=V(x_2)+|x_1|^2/2 $, where $V$ is a potential and $|x_1|^2/2 $ corresponds to the kinetic energy of a particle with unit mass, the corresponding drift $\gF_H $ would write $\gF_H(\x)=(-\nabla_{x_2} V(x_2), x_1)^*$. Adding a damping term $\mathbf D(\x)$, i.e. for $\gF(\x)=(\gF_H-{\mathbf D})(\x)$, leads to investigate the long time behavior of the system, we can e.g. refer to the works \cite{gada:panl:14}, \cite{hera:nier:04} for related discussions, to the monograph \cite{soiz:94} for applications in mechanics  or to \cite{tala:02}, \cite{matt:stua:high:02} for numerical approximations of the invariant measures. 

In mathematical finance, \textcolor{black}{equation \eqref{SDE}} can be related to the model used to price path-dependent contracts, such as Asian options (see, \cite{baru:poli:vesp:01} or \cite{cibe:poli:ross:19} for recent developments).\\
 
We choose here to focus on the very object behind, the density, over a finite time interval. To expose some of the particular features of the model, let us start our discussion \textcolor{black}{with} the (striking) Gaussian setting.\\

\noindent\textbf{Gaussian case and the H\"ormander condition}:
For illustrative purposes let us examine the case $F_1\equiv 0$, $\sigma\equiv 1$ and $F_2(t,X_t^1,X_t^2)=X^1_t$, which corresponds to the Langevin dynamics in its simplest form: 
\begin{equation}\label{Langevin_SDE}
\dif X_t^1=\dif W_t, \ \  \dif X_t^2=X_t^1 \dif t, \qquad t\geq 0,
\end{equation}
which equivalently rewrites in the short form \eqref{SDE0}
\begin{equation}\label{Langevin_SDE_SHORT}
\dif \X_t={\mathbf A}\X_t \dif t+ B  \dif W_t, {\mathbf A}=\begin{pmatrix}0_{d\times d} &0_{d\times d}\\
\mathbb I_{d\times d}& 0_{d\times d}\end{pmatrix}, \qquad t\geq 0.
\end{equation}

In his seminal work \cite{kolm:33}, Kolmogorov derived the fundamental solution for the PDE \eqref{PDE_FUNDA} associated with the above process.  For an initial value $\x=(x_1,x_2)\in\mR^{2d}$, we have for $t\geq 0$,
$$
\X_t=(X^1_t,X_t^2)=\left(x_1+ W_t, x_2+x_1t+\int_0^tW_s \dif s\right)=\exp(\mathbf A t)\x+\int_0^t \exp\big(\mathbf A(t-s)\big)B\dif W_s,
$$ 
which is a Gaussian process with mean 
and covariance matrix respectively given by
\begin{equation}
\btheta_t(\x)=\exp({\mathbf A} t){\textcolor{black}{\x}} =(x_1, x_2+x_1t), \qquad 
\K_t=\begin{pmatrix}
t \mI_{d\times d}& \frac{t^2}{2}\mI_{d\times d}\\ \frac{t^2}{2}\mI_{d\times d} & \frac{t^3}{3}\mI_{d\times d}
\end{pmatrix}.
\end{equation}
The matrix $\K_t$ is positive definite for every $t> 0$ and therefore the process admits a density for every $t> 0$, explicitly given by
\begin{equation}
\y \mapsto \Big(\frac{\sqrt3}{\pi t^2}\Big)^d\exp\left(-\tfrac{1}{2}|\K_t^{-\frac 12}(\y-\btheta_t(\x))|^2\right)=p(0,\x;t,\y)=:p(\x;t,\y).
\end{equation}
In particular, there exist constants $0<c_- < c_+$ such that
\begin{eqnarray}
\Big(\frac{\sqrt3}{\lambda\pi t^2}\Big)^d\exp{\left(-c_+|\T_t^{\textcolor{black}{-1}}(\y-\btheta_t(\x))|^2\right)}\leq p(\x;t,\y) \leq\Big(\frac{\sqrt3}{\lambda\pi t^2}\Big)^d\exp{\left(-c_-|\T_t^{\textcolor{black}{-1}}(\y-\btheta_t(\x))|^2\right)},
\end{eqnarray}
where, for $t>0$,
\begin{align}\label{DEF_T}
\green{
\T_{t}=\begin{pmatrix} t^{\frac 12}\mI_{d\times d}& 0_{d\times d}\\
0_{d\times d}& t^{\frac 32}\mI_{d\times d}
\end{pmatrix},
}
\end{align}
is the scale matrix which precisely reflects the multi-scale behavior of the components.\\

Let us first remark that, in H\"ormander form, the generator of the process in \eqref{Langevin_SDE} writes
\begin{align}
L=A_1^2+A_0,\ A_1= \nabla_{x_1}, \ A_0(\x) := x_1\nabla_{x_2},
\end{align}
so that, denoting by $[\cdot,\cdot] $ the Lie bracket, $[A_1,A_0]=  \nabla_{x_2}$ and ${\rm Span}\{A_1,[A_1,A_0]\}=\R^{2d} $. Importantly, we see that the drift is really needed to span the whole space. This is why we speak
about \textit{weak} H\"ormander condition. As we have seen, this kind of assumption leads to a multi-scale behavior of the heat-kernel as opposed to the strong H\"ormander condition, i.e. when the diffusive vector fields and their Lie brackets span the space. In that case, two-sided heat kernel bounds, which exhibit a usual parabolic scaling, in $\sqrt{t} $ w.r.t. the Carnot metric induced by the vector fields, are available in \cite{kusu:stro:87}. There is therefore a drastic difference between these two types of assumption\textcolor{black}{s}. 

Let us eventually mention that, since we are going to consider \textit{rough} coefficients, we will not be able to perform Lie bracketing to justify the existence of the density from the H\"ormander condition. The non degeneracy of $\nabla_{x_1}F_2 $ can somehow be seen as a \textit{\textcolor{black}{mild} weak} H\"ormander type condition\footnote{\textcolor{black}{In the document, we will further refer with a slight terminology abuse to this assumption as the weak type H\"ormander condition.}}. We emphasize that it is precisely this term which makes the Kolmogorov example work, because it precisely allows the noise on the first component to propagate to the second one.\\

Secondly, it can be observed from this example that the time-scale in the density is not diffusive: the fluctuations of the two components $X^1$, $X^2$ are of order of $t^{1/2}$ and $t^{3/2}$ respectively, which corresponds to the intuition that the typical time-scale of an integrated Brownian motion should be equal to the integral of the time-scale of the Brownian motion. 
This phenomenon also appears in the deviation term in the exponential. The growth rate is different for the two components. Also, the unbounded drift term induces deviations w.r.t. to the transport of the initial condition by the underlying deterministic differential system $\btheta_t(\x) $ and not the starting point itself $\x $, \textcolor{black}{normalized w.r.t. the previous intrinsic time scales}.

Similarly, it is seen that there exists $C\ge 1$ s.t. for $i\in \{1,2\} $,
\begin{eqnarray*}
|\nabla_{x_i} p(\x;t,\y)|&\le& |\big( (\K_t^{-\frac 12}\nabla \btheta(\x))^*\K_t^{-\frac 12}(\btheta_t(\x)-\y))\big)_i|  p(\x;t,\y)\\
&\le& \frac{C}{ t^{\frac{2i-1} 2}}\Big(\frac{\sqrt3}{\pi t^2}\Big)^d\exp\left({-C^{-1}|\T_t^{\textcolor{black}{-1}}(\y-\btheta_t(\x))|^2}\right).
\end{eqnarray*}
Namely, a differentiation with respect to the non-degenerate variable induces an additional time singularity of the corresponding typical order rate $t^{-1/2}$ whereas a differentiation with respect to the degenerate variable  induces a time singularity in $t^{-3/2}$ i.e. at its corresponding typical rate.

Our goal in the current work is to extend those heat kernel and gradient bounds, obtained directly for the Kolmogorov example, to the densities of SDEs solving \eqref{SDE} under somehow minimal smoothness assumptions on the coefficients.\\

{\color{black} \noindent\textbf{Available results and \emph{minimal} conditions from a regularization by noise perspective.}  Since the seminal work of Kolmogorov \cite{kolm:33}, such equations have been thoroughly investigated in the literature both from the analytic or probabilistic viewpoint\footnote{Let us also recall that the work \cite{kolm:33} was recalled by H\"ormander as the starting point of his general theory of hypoellipticity \cite{horm:67}.}.  Existence of fundamental solution for the underlying parabolic PDE \eqref{PDE_FUNDA} was first obtained through a \textit{parametrix} type perturbation technique, for smooth enough coefficients,  by Weber \cite{webe:51}. We can also refer to Sonin for further results in that direction \cite{soni:67}. On the other hand, density estimates were derived in \cite{kona:meno:molc:10} (global upper and lower \textit{diagonal} bound) and then extended to more general models of SDEs that can be seen as perturbed ODEs for which a noise acting on the first component will transmit to the whole chain of ODEs through a weak H\"ormander like condition. For such models we refer to \cite{dela:meno:10}, \cite{meno:10}. From those works one can derive  two-sided heat kernel bounds for the density of \eqref{SDE} when the drift is globally Lipschitz in space, i.e. when the drift part of the dynamics in \eqref{SDE} can be associated with a usual well-posed ODE, {\color{black} and when the diffusion coefficient is H\"older continuous in space. Eventually, in a smooth framework, Pigato derived in \cite{piga:22} heat kernel and gradient estimates for the system as well as short time asymptotics. Except for the diffusion coefficient, the aforementioned works do not take advantage of the propagation of the noise through the system, in the sense that the drift is always assumed to be (at least) a Lipschitz in space function whereas rougher drifts may be of interest see e.g. Section 1.3 in \cite{IM21}. As already claimed, our aim consists in deriving those bounds under rather \textit{minimal} conditions, meaning that we manage to benefit from the regularization by noise phenomenon{\color{black}, see the Saint Flour lectures notes \cite{flandoli_random_2011} for an overview of such kind of phenomenom.}}\\

 Regularization by noise for degenerate system was investigated in  e.g.  \cite{chau:17} from a strong point of view, meaning that the author exhibited therein, within the framework of H\"older spaces, some minimal thresholds on $F$ that guarantee strong well posedness of equation \eqref{SDE}, in spite of a non Lipschitz drift. Generalization to H\"older-Dini coefficients were investigated in \cite{wang:zhan:16} with other techniques {\color{black} but somehow similar threholds}. For rougher drifts, corresponding namely to Bessel potentials, but for $F_2(t,x)=x_1$, the strong well posedness was also derived (for the same previous \textit{regularity threshold on the degenerate variable}) in \cite{fedr:flan:prio:vove:17} and \cite{zhan:16} where the critical case for the regularity index is considered. However, it is well known that existence of the density does not rely on strong well-posedness of the system, but more generally on its weak well-posedness.\\

In the current degenerate setting, the weak regularization by noise was investigated for the kinetic case in \cite{chau:16}. The author derived therein, still in the setting of H\"older spaces, smaller thresholds on the regularity index of the drift w.r.t. the degenerate variable that yield weak uniqueness. It is also importantly shown that these thresholds are (almost) sharp in the sense that, when the drift of the degenerate component has H\"older regularity below the threshold, there are counter-examples to weak uniqueness. This weak well-posedness and associated counter-examples have then been extended to a full chain of perturbed ODEs in \cite{Raynal2017RegularizationEO}. In this latter work, the Authors exhibit an (almost) sharp characterization, at each level of the ODE chain, of the regularity index needed in each variable to restore weak well-posedness. In comparison with our current setting, it is proved therein that the drift of each component feels differently the degenerate variable: whereas one only needs positive regularity indexes for the drift of the non degenerate component w.r.t. all variables for weak uniqueness to hold, a threshold of $1/3$ appears for the regularity index of the drift of the degenerate component w.r.t. to the degenerate variable. We will try to match this setting as much as possible for what we aim at doing here.
We refer to Theorem \ref{MAIN_THM} and \ref{Grad1} and associated remarks for details.}\\

{\color{black} \noindent\textbf{Objective and strategy.}} Many other type of estimates have been established for the SDE \eqref{SDE} or its formal generator. Let us mention among them: Harnack inequalities \cite{lanc:poli:94}, \cite{pasc:poli:06}; related heat kernel estimates for operators in divergence form with measurable coefficients \cite{lanc:pasc:17}, \cite{lanc:pasc:poli:20}; Schauder estimates, see \cite{IM21} for the current framework and  \cite{Raynal2018SharpSE} for an even more general case (one can also refer to \cite{hao:wu:zhan:20} for an extension to kinetic non-local operators - with an application to strong well posedness - and to \cite{Ma20} as well for a more general framework); $L^p$ estimates, see \cite{huan:meno:prio:16}, \cite{meno:17} or \cite{Chen:Zhang18}. Let us eventually mention the work \cite{pasc:pesc:22} which deals with the associated Stochastic PDE in the two-dimensional case or \cite{zhan:21} which investigates the well-posedness of a McKean-Vlasov version of \eqref{PDE_FUNDA} through the De Giorgi approach.\\

We will focus here on the density/heat kernel and will adapt the approach already considered in \cite{MPZ20}, \cite{meno:zhan:20}, for non degenerate SDEs with unbounded drift respectively Brownian and stable driven,  to the current degenerate case. The first step consists in obtaining two-sided estimates. This is done using \textit{forward} type parametrix or Duhamel type expansions, as e.g. considered in the classic non-degenerate case in \cite{frie:64} or \cite{mcke:sing:67}, with the additional difficulty that, because of the unbounded drift the parametrix series needs to be truncated. The tails of the series {\color{black}are} controlled through stochastic control arguments (see \cite{dela:meno:10}, \cite{zhan:09}). For the estimates on the derivatives
the idea consists  in mixing \textit{forward} and \textit{backward} Duhamel expansions and to consider suitable normalizations which exploit thoroughly the underlying two-sided estimates.

We restrict in this work to the kinetic case for simplicity. We believe that our main results would extend to the full perturbed chain of ODE as considered in \cite{Raynal2017RegularizationEO} under suitable assumptions on the coefficients. In this more general setting the idea would be to couple the {\color{black}current approach} with the computations performed in \cite{chau:hono:meno:18} to derive strong well posedness for the full chain.\\



The article is organized as follows. We state our main results (Theorems \ref{MAIN_THM} and \ref{Grad1}) in Section \ref{MR}. Section \ref{SEC_PREL} gathers some technical results about mollified flows associated with H\"older in space coefficients and also addresses the corresponding deterministic control problem which will be useful for the two-sided heat kernel estimate. We will establish in Section \ref{SEC_DENS_SMOOTH} our main results for smooth coefficients which satisfy the assumptions of Theorems \ref{MAIN_THM} and \ref{Grad1} and we will carefully prove that the constants in the estimates obtained do not depend on such smoothness. This is precisely why we then derive our main result\textcolor{black}{s} through compactness arguments detailed in Section \ref{SEC_COMP}.

\subsection{Statement of main results}
\label{MR}
Let $d,l\in\mN$. For $j\in\{0\}\cup\N$ and $\gamma\in \textcolor{black}{[0,1)}$, let $\sC^{j+\gamma}(\R^d;\R^{l})$ be the space of H\"older functions from $\mR^d$ to $\mR^l$ defined by
$$
\sC^{j+\gamma}(\R^d;\R^{l}):=\left\{f: \|f\|_{\sC^{j+\gamma}(\R^d;\R^{l})}:=\sum_{k=1}^j \|\nabla^k f\|_{L^\infty(\R^d;\R^l)}+\sup_{x\not=y,|x-y|\leq 1}\frac{|\nabla^j f(x)-\nabla^jf(y)|}{|x-y|^\gamma}<\infty\right\},
$$
where $\nabla^k$ stands for the $k$-order gradient. \textcolor{black}{Note that the functions in $\sC^{j+\gamma}(\R^d;\R^{l})$ can be unbounded and have sublinear growth}.

\textcolor{black}{Importantly the functions in  $\sC^0(\R^d;\R^{l})$ can  be possibly discontinuous and also satisfy (see \cite{wang:zhan:16}, Lemma 2.3)
$$
\sup_{x\neq y} \frac{|f(x)-f(y)|}{1+|x-y|}<+\infty.
$$
}
We denote by $\x=(x_1,x_2)$ a point in $\R^d\times\R^d$ and by $\nabla_{x_1}$, $\nabla_{x_2}$ the gradients with respect to the first and second set of variables, respectively.
Following the previous discussion, it is natural to endow $\R^{2d}$ with an anisotropic distance, corresponding to the intrinsic scale matrix \textcolor{black}{\eqref{DEF_T}}:
\begin{equation}\label{norm}
|\x|_{\mathbf{d}}:=|x_1|+|x_2|^{\frac 13}, \qquad \x=(x_1,x_2)\in\R^{2d}.
\end{equation}
Next, we recall the definition of a \emph{anisotropic H\"older spaces} associated with \eqref{norm} (see, for instance \cite{Lunardi:97}, \cite{Raynal2018SharpSE}).
We say that a vector valued function $f\in \sC_{\mathbf{d}}^{j+\gamma}(\R^{2d};\R^l)$ if
\begin{equation}
\|f\|_{\sC_{\mathbf{d}}^{j+\gamma}(\R^{2d};\R^l)}:=\sup_{x_2\in\R^{d}}\|f(\cdot,x_2)\|_{\sC^{j+\gamma}(\R^d;\R^l)}
+\sup_{x_1\in\R^{d}}\|f(x_1,\cdot)\|_{\sC^{(j+\gamma)/3}(\R^d;\R^l)}<\infty.
\end{equation}
In particular, for $f\in\sC_{\mathbf{d}}^{1+\gamma}(\R^{2d};\mR)$, by Taylor's expansion, we have
\begin{align}\label{Taylor}
|\cT_f(\x,\y):=f(\x)-f(\y)-\nabla_{x_1}f(\y) (\x-\y)_1|\leq C_\gamma\|f\|_{\sC^{1+\gamma}_{\bf d}}|\x-\y|^{1+\gamma}_{\bf d}.
\end{align}



We assume the following conditions to hold:
\begin{itemize}
\item[{\bf(H$^\gamma_\sigma$)}] 
  There exist $\gamma\in(0,1]$ and $\kappa_0\ge 1$ such that for all $\green{(t,\x)}\in \R_+\times \R^{2d}$ and $\xi\in \R^d$,
$$ \kappa_0^{-1}|\xi|^2\le \langle \sigma\sigma^*\green{(t,\x)} \xi,\xi\rangle \le \kappa_0 |\xi|^2
$$ 
and
$$
|\sigma\green{(t,\x)}-\sigma(t,\y)|\leq \kappa_0|\x-\y|_{\bf d}^\gamma,
$$
where  $|\cdot|$ denotes the Euclidean norm, $\langle \cdot, \cdot \rangle $ is the inner product and  ${}^* $ stands for the transpose.

\item[{\bf(H$^\gamma_{\F}$)}] For some $\gamma\in \textcolor{black}{(0,1]}$ and $\kappa_1,\kappa_2>0$, it holds that
$$
|F_i(t,{\bf 0})|\le \kappa_i,\ i=1,2,\ \|F_1(t,\cdot)\|_{\sC^0_{\bf d}}\leq \kappa_1,\ \|F_2(t,\cdot)\|_{\sC^{1+\gamma}_{\bf d}}\leq \kappa_2.
$$
Moreover, there exists a closed convex subset ${\mathcal E}\subset GL_{d}(\R)$
(the set of all invertible $d \times d$  matrices)
 such that $\nabla_{x_{1}} \green{F_2(t,\x)}\in {\mathcal E}$ for all $t\geq 0$ and $\x\in\R^{2d}$. 
\end{itemize}

We introduce the following notation
\begin{align}\label{GG2}
g_\lambda\green{(t,\x)}:=t^{-2d}\e^{- |\T^{-1}_t\x|^2/(2\lambda)}
\end{align}
as well as the following set of parameters for later use: for $T>0$,
$$
\Theta_T:=(T,\kappa_0,\kappa_1,\kappa_2,d,\gamma,\cE).
$$

\textcolor{black}{Eventually, to state our main results we need to introduce a \textit{mollified} flow associated with the drift $F$ in \eqref{SDE}, which under {\bf(H$^\gamma_{\F}$)} is \textit{rough}. Namely,
\begin{equation}\label{FLOW_MACRO}
{\dot{\wt\btheta}}_{t,s}(\x)=\wt\gF(t,\wt\btheta_{t,s}(\x)), \quad \wt\btheta_{s,s}(\x)=\x,\ 
\end{equation}
where 
$$
\wt\gF(t,\x)=\Big([F_1(t,\cdot)*\rho_1](\x), [F_2(t,\cdot)*\rho_{|t-s|^{3/2}}](\x)\Big),
$$
and $\rho_\eps(\x)=\eps^{-2d}\rho(\eps^{-1}\x)$ and $\rho$ is a smooth density function with compact support and $*$ stands for the convolution in space. The first regularization, performed at a macro level, is very natural to introduce a flow since the initial drift coefficient is not necessarily smooth. Regularizing the second component,  allows as well to have a  flow defined in the classical sense. The regularization parameter, corresponding to the intrinsic time scale of the component allows to have the \textit{equivalence} between this flow and any other regularized flow (see Remark \ref{RK_PEANO} and Lemma \ref{lemme:bilipflow} for details).}

Our first main result of this paper is stated as follows.
\bt\label{MAIN_THM}
Under {\bf(H$^\gamma_\sigma$)} and {\bf(H$^\gamma_{\F}$)} with $\gamma\in(0,1]$,
for any $T>0$ and $0\le s<t\le T$, there exists a unique weak solution ${\bf X}_{t,s}(\x)$ of \eqref{SDE} 
starting from $\x$ at time $s$ which admits a density $p(s,\x;t,\y)$ continuous in $\x,\y\in\R^{2d}$.
Moreover, $p(s,\x;t,\y)$ enjoys the following estimates: 
\begin{itemize}
\item[(i)] (Two sided estimates) There are $\lambda_0, C_0\ge 1$ depending on $\Theta_T$ such that for all $0\le s<t\le T$ and $\x,\y\in\R^{2d}$,
\begin{align}
C_0^{-1}g_{\lambda^{-1}_0}\(t-s, \wt\btheta_{t,s}(\x)-\y\)
 \le p(s,\x;t,\y)\le C_0g_{\lambda_0}\(t-s, \wt\btheta_{t,s}(\x)-\y\). \label{Density_bounds_THM}
\end{align}
\item[(ii)] (Gradient estimate in $x_1$) There exist constants $\lambda_1,\ C_1\ge 1$ depending on $\Theta_T$ such that
for any $0\le s<t\le T$ and $\x,\y\in\R^{2d}$,
\begin{align}
\left|\nabla_{x_1}p(s,\x;t,\y) \right|&\le {\textcolor{black}{C_1}}(t-s)^{-\frac{1}{2}}g_{\lambda_1}\(t-s, \wt\btheta_{t,s}(\x)-\y\). \label{Derivatives_x1_THM}
\end{align}
\item[(iii)] (H\"older estimate in $\x$) Let $\eta_0,\eta_1\in(0,1)$. For $j=0,1$, there exist constants $\lambda_j,\ C_j\ge 1$ depending on $\Theta_T$ and $\eta_j$ such that
for any $0\le s<t\le T$ and $\x,\x',\y\in\R^{2d}$,
\begin{align}\label{Holder}
\begin{split}
&\left|\nabla^j_{x_1}p(s,\x;t,\y)-\nabla^j_{x_1}p(s,\x';t,\y) \right|\le {C_j}|\x-\x'|_{\bf d}^{\eta_j}(t-s)^{-(\frac{j}{2}\textcolor{black}{+\eta_j})}\\
&\qquad\times
\Big(g_{\lambda_j}\(t-s, \wt\btheta_{t,s}(\x)-\y\)+g_{\lambda_j}\(t-s, \wt\btheta_{t,s}(\x')-\y\)\Big). 
\end{split}
\end{align}
\end{itemize}
\et
\br\rm
Under {\bf(H$^\gamma_{\F}$)}, $F_1$ may be unbounded and discontinuous. For instance, if $F_1(\x)=F_{11}(\x)+F_{12}(\x)$ and $F_2(x_1,x_2)=x_1$, where $F_{11}$ is bounded measurable and
$F_{12}$ is global Lipschitz, then $\F=(F_1,F_2)$ satisfies {\bf(H$^{1}_{\F}$)}. This example corresponds to the standard kinetic SDEs.

\textcolor{black}{Let us emphasize that this Theorem holds under the assumptions that have been shown in  \cite{Raynal2017RegularizationEO} to be minimal to guarantee weak uniqueness for the solution of \eqref{SDE}. In particular the two sided estimates \eqref{Density_bounds_THM} specify the Krylov type estimate of  \cite{Raynal2017RegularizationEO} which roughly said that in $L^q-L^p$ norms (for suitable indexes $p,q$) the density behaved \textit{as} the Kolmogorov one appearing in \eqref{Density_bounds_THM}}.
\er

{\color{black}\br[About the flow in the above estimates] \label{RK_PEANO}\rm 
We point out that the above estimates could also be stated replacing $\widetilde \btheta_{t,s}(\x) $ introduced in \eqref{FLOW_MACRO} by any Peano flow $\btheta_{t,s}^{(1)}(\x) $ solving:
\begin{equation*}
{\dot{\btheta}}_{t,s}^{(1)}(\x)=\gF^{(1)}(t,\btheta_{t,s}^{(1)}(\x)), \quad \btheta_{s,s}^{(1)}(\x)=\x,\ 
\end{equation*}
where 
$$
\gF^{(1)}(t,\x)=\Big([F_1(t,\cdot)*\rho_1](\x),F_2(t,\cdot)(\x)\Big).
$$
Pay attention that  $F_2$ is not regularized whereas the possibly discontinuous component $F_1$ still needs to be to define an underlying flow.\\

Indeed it can be shown, similarly to the estimates established in Lemma \ref{lemme:bilipflow} below, that there exists a constant $C:=C(\Theta_T)$ s.t.:
$$\big|\T_{t-s}^{-1}\big(\btheta_{t,s}^{(1)}(\x)-\wt \btheta_{t,s}(\x) \big)\big|\le C,$$
i.e. $\btheta_{t,s}^{(1)}(\x)$ and $\wt \btheta_{t,s}(\x)$ are equivalent with respect to the intrinsic scales.
\er
}
\vspace*{.3cm}
In the above theorem, the regularity assumptions on $\sigma$ and $\F$ are almost sharp. To obtain the second order derivative estimate in $x_1$ and the first order gradient estimate in $x_2$, we have to
make further regularity assumption\textcolor{black}{s} as stated below.
\bt\label{Grad1}
In the situation of Theorem \ref{MAIN_THM}, we also assume that for the same $\gamma\in(0,1]$,
\begin{align}\label{F11}
\|F_1(t,\cdot)\|_{\sC^\gamma_{\bf d}}\leq\kappa_1,\ \ t\geq 0.
\end{align}
\begin{itemize}
\item[(i)] (Second order derivative estimate in $x_1$) There exist constants $\lambda_1,\ C_1\ge 1$ depending on $\Theta_T$ such that
for any $0\le s<t\le T$ and $\x,\y\in\R^{2d}$,
\begin{align}
\left|\nabla^2_{x_1}p(s,\x;t,\y) \right|&\le {C_1}(t-s)^{-1}g_{\lambda_1}\(t-s, \wt\btheta_{t,s}(\x)-\y\). \label{Sec1}
\end{align}
\item[(ii)] (H\"older estimate of $\nabla^2_{x_1}p$ in $\x$) For any $\eta_2\in(0,\gamma)$, there exist constants $\lambda_2,\ C_2\ge 1$ depending on $\Theta_T$ and $\eta_2$ such that
for any $0\le s<t\le T$ and $\x,\x'\y\in\R^{2d}$,
\begin{align}\label{Holder1}
\begin{split}
&\left|\nabla^2_{x_1}p(s,\x;t,\y)-\nabla^2_{x_1}p(s,\x';t,\y) \right|\le {C_2}|\x-\x'|_{\bf d}^{\eta_2}(t-s)^{-1-\frac{\eta_2}{2}}\\
&\qquad\times\Big(g_{\lambda_2}\(t-s, \wt\btheta_{t,s}(\x)-\y\)+g_{\lambda_2}\(t-s, \wt\btheta_{t,s}(\x')-\y\)\Big). 
\end{split}
\end{align}
\end{itemize}
 If $\sigma$ also satisfies that 
\begin{align}\label{Re1}
|\sigma(t,\x)-\sigma(t,\y)|\leq \kappa_0(|(\x-\y)_1|^{\alpha}+(\x-\y)_2|^{\frac{1+\gamma}{3}}),
\end{align}
where $\gamma$ is the same as in Theorem \ref{MAIN_THM} and $\alpha\in((1-\gamma)\vee\gamma,1]$, and
\begin{align}\label{Re2}
|F_1(t,\x)-F_1(t,\y)|\leq \kappa_1(|(\x-\y)_1|^{\gamma}+(\x-\y)_2|^{\frac{1+\gamma}{3}}),
\end{align}
\textcolor{black}{then}

\begin{itemize}
\item[(iii)] (Gradient estimate in $x_2$)
 there exist $\lambda_3,\ C_3\ge 1$ depending on $\Theta_T$ and $\alpha$ such that 
for any $0\le s<t\le T$ and $\x,\y\in\R^{2d}$,
\begin{align}
\left|\nabla_{x_2}p(s,\x;t,\y) \right|&\le {C_3}(t-s)^{-\frac{3}{2}}g_{\lambda_3}\(t-s, \wt\btheta_{t,s}(\x)-\y\).\label{Derivatives_x2_THM}
\end{align}
\item[(iv)] (H\"older estimate of $\nabla_{x_2}p$ in $\x$) For any $\eta_3\in(0,(\alpha-\gamma)\wedge (\alpha+\gamma-1))$, there exist constants $\lambda_4,\ C_4\ge 1$ depending on $\Theta_T$, $\alpha$ and $\eta_3$ such that
for any $0\le s<t\le T$ and $\x,\x',\y\in\R^{2d}$,
\begin{align}\label{Holder2}
\begin{split}
&\left|\nabla_{x_2}p(s,\x;t,\y)-\nabla_{x_2}p(s,\x';t,\y) \right|\le {C_4}|\x-\x'|_{\bf d}^{\eta_3}(t-s)^{-\frac{3+\eta_3}{2}}\\
&\qquad\times\Big(g_{\lambda_4}\(t-s, \wt\btheta_{t,s}(\x)-\y\)+g_{\lambda_4}\(t-s, \wt\btheta_{t,s}(\x')-\y\)\Big). 
\end{split}
\end{align}
\end{itemize}

\et

\br\rm {\color{black}
For gradient estimates in the degenerate component $x_2$, we need extra regularities \textcolor{black}{\eqref{Re1} and \eqref{Re2}} since for kinetic operators, we only have $\frac23$-gain of regularity in $x_2$. Let us  briefly comment this additional regularity, which might not seem sharp at first sight.  Our feeling is that such assumptions are actually sharp with respect to the methodology employed here. Indeed, as our starting point to estimate the density relies on a first order parametrix expansion, see \eqref{D0}, we end up with an implicit representation of the density involving as well all the coefficients of the system. The crucial point is that, when estimating the gradient in the degenerate directions of the implicit representation of the density, we make all the coefficients feel the differentiation w.r.t. the degenerate variables. This roughly explains why we impose assumptions \eqref{Re1}, \eqref{Re2} which lead to similar (w.r.t. the degenerate variables) conditions as the one previously imposed on $F_2$ in {\bf(H$^\gamma_{\F}$)}. Assuming the same regularity in the degenerate directions for the whole drift $\mathbf F = (F_1,F_2)$ already appeared in \cite{chau:16} in connection with weak uniqueness for \eqref{SDE} as well as in \cite{chau:hono:meno:18} and \cite{hao:wu:zhan:20} to derive Schauder like estimates for strong uniqueness purposes. Note that in those frameworks, since the diffusion coefficient was assumed to be Lipschitz, assumption \eqref{Re1} did not explicitly appears. Let us eventually conclude by emphasizing that, under similar assumptions as the one of Theorem \ref{MAIN_THM}, the Authors in \cite{Raynal2018SharpSE}  only succeeded in deriving gradient estimates in the non degenerate directions, but only H\"older estimates in the degenerate ones. 
}
%
%
%
\er

For  two quantities $Q_1$ and $Q_2$, we will frequently use the notation $Q_1\lesssim Q_2 $ meaning that that there exists $C:=C(\Theta_T)$  such that $Q_1\le CQ_2 $.
\section{Preliminaries}
\label{SEC_PREL}

In this section we assume that $\F=(F_1,F_2)$ satisfies {\bf(H$^0_{\F}$)} and temporarily assume that
$$
\|\nabla_\x \F\|_\infty<\infty.
$$
In particular, for some $\kappa_1>0$,
\begin{align}\label{AX1}
|F_1(t,\x)-F_1(t,\y)|\leq \kappa_1(1+|\x-\y|),
\end{align}
and for some $\kappa_2>0$,
\begin{align}\label{AX2}
|F_2(t,\x)-F_2(t,\y)|\leq \kappa_2( |(\x-\y)_1|+|(\x-\y)_2|^{\frac13}+|(\x-\y)_2|).
\end{align}
For $s,t\geq 0$ and $\x\in\R^{2d}$, let $\btheta_{t,s}(\x)$ be the regularization flow defined by the differential system
\begin{equation}\label{FLOW}
{\dot\btheta}_{t,s}(\x)=\gF(t,\btheta_{t,s}(\x)), \quad \btheta_{s,s}(\x)=\x.
\end{equation}
Here $(\btheta_{t,s}(\x))_{t\geq s}$ stands for a forward flow, while $(\btheta_{t,s}(\x))_{t\leq s}$ stands for a backward flow. 
In particular, let $(\btheta_{t,s}(\x))^{-1}$ be the inverse of $x\mapsto \btheta_{t,s}(\x)$. Then
\begin{equation}\label{INV}
(\btheta_{t,s}(\x))^{-1}=\btheta_{s,t}(\x).
\end{equation}

\subsection{Equivalence of measurable flow for ODEs}

We recall the following \textcolor{black}{(sublinear)} Gronwall type lemma.
\bl\label{Le11}
Let $f(t):[0,\infty)\to[0,\infty)$ be a continuous function. Suppose that for some $\alpha\in(0,1)$ and $c_1,c_2\geq0$,
$$
f(t)\leq f(0)+c_1\int^t_0 f(s)^\alpha\dif s+c_2\int^t_0 f(s)\dif s,\ t>0.
$$
Then
$$
f(t)\leq\e^{c_2t} f(0)+(c_1\e^{c_2t}(1-\alpha)t)^{\frac{1}{1-\alpha}},\ t>0.
$$
\el
We have the following crucial lemma, which corresponds to Lemma 1.1 of \cite{MPZ20} and \cite{Raynal2017RegularizationEO}.
\bl\label{lemme:bilipflow}
Under {\bf(H$^0_{\F}$)}, for any $T>0$, there exist constants $\kappa_3,\geq 1$ only depending on $\kappa_1,\kappa_2,d, T$  such that for all 
$0\leq s\leq r<t\leq T $
and $\x,\y\in\mR^{2d}$,
\begin{equation}
\label{EQ_EQUIV_FLOW}
\kappa^{-1}_3\(|\T^{-1}_{t-s}(\x-\btheta_{r,t}(\y))|-1\)\leq |\T^{-1}_{t-s}(\btheta_{t,r}(\x)-\y)|
\leq \kappa_3\( |\T^{-1}_{t-s}(\x-\btheta_{r,t}(\y))|+ 1\).
\end{equation}
\el
\begin{proof}
To show \eqref{EQ_EQUIV_FLOW}, by \eqref{INV} and the symmetry, it suffices to show
$$
|\T^{-1}_{t-s}(\btheta_{t,r}(\x)-\btheta_{t,r}(\y))|
\leq \kappa_3\(|\T^{-1}_{t-s}(\x-\y)| +1\).
$$
Without loss of generality, we may assume $0=r<t\leq T$, and 
write for $\x,\y\in\R^{2d}$ and $t\geq 0$,
$$
\ell_i(t):=\big|\(\btheta_{t,0}(\x)-\btheta_{t,0}(\y)\)_i\big|,\ \ i=1,2.
$$ 
Using the above notation and by definition \textcolor{black}{\eqref{DEF_T}}, we only need to show that for all $0<t\leq t_0\leq T$,
\begin{align}\label{AA2}
t_0\ell_1(t)+\ell_2(t)\leq \kappa_3\(t_0\ell_1(0)+\ell_2(0)+t_0^{\frac32}\).
\end{align}
For $i=1$, by \eqref{AX1} we have
\begin{align*}
\ell_1(t)&\leq\ell_1(0)+\int^t_0|F_1(r,\btheta_{r,0}(\x))-F_1(r,\btheta_{r,0}(\y))|\dif r\\
&\leq \ell_1(0)+\kappa_1t+\kappa_1\int^t_0(\ell_1(r)+\ell_2(r))\dif r,
\end{align*}
which implies by Gronwall's inequality that
\begin{align}
\ell_1(t)\leq\e^{\kappa_1 t}(\ell_1(0)+\kappa_1t)+\kappa_1\e^{\kappa_1 t}\int^t_0\ell_2(r)\dif r.\label{AA1}
\end{align}
For $i=2$,  by \eqref{AX2} and \eqref{AA1}, we have
\begin{align}\label{AX6}
\begin{split}
\ell_2(t)&\leq\ell_2(0)+\int^t_0|F_2(r,\btheta_{r,0}(\x))-F_2(r,\btheta_{r,0}(\y))|\dif r\\
&\leq\ell_2(0)+\int^t_0\ell_1(r)\dif r+\kappa_2\int^t_0\(\ell_2(r)^{\frac13}+\ell_2(r)\)\dif r\\
&\leq\ell_2(0)+\e^{\kappa_1 t}(t\ell_1(0)+\kappa_1t^2)+\textcolor{black}{\kappa_2}\int^t_0\ell_2(r)^{\frac13}\dif r+\textcolor{black}{(\kappa_2+e^{\kappa_1 t })}\int_0^t\ell_2(r)\dif r.
\end{split}
\end{align}
In particular, by  Lemma \ref{Le11}, for $\textcolor{black}{\tilde \kappa_2=\kappa_2+e^{\kappa_1 t }} $, we have
\begin{align*}
\ell_2(t)&\leq \e^{\textcolor{black}{\tilde \kappa_2}t}\(\ell_2(0)+\e^{\kappa_1 t}(t\ell_1(0)+t^2)\)+\textcolor{black}{(\kappa_2\e^{\tilde \kappa_2t}\tfrac 23t)^{\frac{3}{2}}}\\
&\leq \e^{\textcolor{black}{\tilde \kappa_2}t}\ell_2(0)+\e^{(\textcolor{black}{\tilde \kappa_2}+\kappa_1)t}t\ell_1(0)+c_2t^{\frac{3}{2}},
\end{align*}
which together with \eqref{AA1} yields \eqref{AA2}.
\end{proof}


\br\rm
By \eqref{EQ_EQUIV_FLOW} and the flow property $\btheta_{t,s}(\x)=\btheta_{t,r}\circ \btheta_{r,s}(\x)$, we also have
\begin{align}\label{EZ1}
|\T^{-1}_{t-s}(\btheta_{r,s}(\x)-\btheta_{r,t}(\y))|-1\lesssim |\T^{-1}_{t-s}(\btheta_{t,s}(\x)-\y)|
\lesssim |\T^{-1}_{t-s}(\btheta_{r,s}(\x)-\btheta_{r,t}(\y))|+ 1.
\end{align}
Moreover, if $|\x-\x'|_{\bf d}\leq C_0(t-s)^{1/2}$ for some $C_0>0$, then
\begin{align}\label{VX2}
|\T^{-1}_{t-s}(\btheta_{t,s}(\x')-\y)|-1\lesssim |\T^{-1}_{t-s}(\btheta_{t,s}(\x)-\y)|
\lesssim |\T^{-1}_{t-s}(\btheta_{t,s}(\x')-\y)|+ 1.
\end{align}
Indeed, by \eqref{EQ_EQUIV_FLOW} we have
\begin{align*}
|\T^{-1}_{t-s}(\btheta_{t,s}(\x')-\y)|
&\lesssim |\T^{-1}_{t-s}(\x'-\btheta_{s,t}(\y))|+1\\
&\leq |\T^{-1}_{t-s}(\x-\btheta_{s,t}(\y))|+|\T^{-1}_{t-s}(\x'-\x)|+1\\
&\leq |\T^{-1}_{t-s}(\x-\btheta_{s,t}(\y))|+C_0+1\\
&\lesssim |\T^{-1}_{t-s}(\btheta_{t,s}(\x)-\y)|+1.
\end{align*}
\textcolor{black}{The other inequality in \eqref{EZ1} is derived by symmetry}. 
\er

\subsection{Gram matrix}

Let $\green{A_t,\sigma_t}: \mR_+\to\mR ^{d}\times\mR^d$ be two measurable maps. Suppose that
for some closed convex subset $\cE$ of $GL(\mR^d)$ and $\kappa_0>0$, and for all $t\geq 0$,
\begin{align}\label{CC1}
A_t\in\cE,\ \ \kappa_0^{-1}|\xi|\leq |\sigma_t\xi|\leq \kappa_0|\xi|.
\end{align}
Define
$$
\bA_t:=\begin{pmatrix}
0_{d\times d}&0_{d\times d}\\
A_t & 0_{d\times d} \end{pmatrix},\ \ \ 
\bTheta_t:=\begin{pmatrix}
\sigma_t\\
0_{d\times d}\end{pmatrix}.
$$
For $s,t\geq 0$, let $\gR_{t,s}$ be the resolvent of $\bA_t$, that is,
\begin{align}\label{Res11}
\p_t\gR_{t,s}=\bA_t\gR_{t,s},\ \ \gR_{s,s}=\mI_{2d\times 2d}.
\end{align}
It is easy to see that the unique solution of \eqref{Res11} is given by
\begin{align}\label{Res0}
\gR_{t,s}=\begin{pmatrix}\mI_{d\times d} & 0_{d\times d}\\
\int^t_s A_r\dif r & \mI_{d\times d} \end{pmatrix}.
\end{align}
From this expression, one sees that for any $s,t,r\geq 0$,
\begin{align}\label{RR1}
\gR_{t,s}^{-1}=\gR_{s,t},\ \ \gR_{t,r}\gR_{r,s}=\gR_{t,s}.
\end{align}
The Gram matrix associated with ${\bf A}_t$ and $\bTheta_t$ is defined by
$$
\K_{t,s}:=\int^t_s \gR_{t,r}\bTheta_r\bTheta_r^*\gR^*_{t,r}\dif r.
$$
The following lemma is well-known. For reader's convenience, we provide detailed proofs here.
\bl\label{Le25}
Under \eqref{CC1}, there is a constant $\kappa\geq1$ depending only on $\kappa_0$ and $\cE$ such that for all $0\leq s<t<\infty$ and $\x\in\mR^{2d}$,
\begin{align}\label{TT1}
|\K_{t,s}^{-1/2}\x|^2=\<\K_{t,s}^{-1}\x,\x\>\asymp_{\kappa}|\T^{-1}_{t-s}\x|^2,
\end{align}
and
\begin{align}\label{TT0}
|\T^{-1}_{t-s}\gR_{t,s}\x|\asymp_{\kappa} |\T^{-1}_{t-s}\x|.
\end{align}
\el
\begin{proof}
By the definition and the change of variables, it is easy to see that
\begin{align}\label{TT2}
\K_{t,s}=\T_{t-s}\hat\K_{1,0}\T_{t-s},
\end{align}
where
$$
\hat\K_{1,0}:=\int^1_0\hat{\gR}_{1,r}\hat\bTheta_r\hat\bTheta_r^*\hat\gR^*_{1,r}\dif r,
$$
and
$$
\hat\bTheta_r=\bTheta_{s+(t-s)r},\ \ \hat{\gR}_{1,r}
=\begin{pmatrix}\mI_{d\times d} & 0_{d\times d}\\
\int^1_r A_{s+(t-s)u}\dif u & \mI_{d\times d} \end{pmatrix}.
$$
Thus, without loss of generality, we may assume $s=0$, $t=1$ and $|\x|^2=|x_1|^2+|x_2|^2=1$. Clearly,
\begin{align}\label{EG7}
\<\K_{1,0}^{-1}\x,\x\>\asymp 1\Leftrightarrow \<\K_{1,0}\x,\x\>\asymp1.
\end{align}
Note that by \eqref{CC1},
$$
\<\K_{1,0}\x,\x\>\textcolor{black}{=\int_0^1 |\bTheta_r^*\hat\gR^*_{1,r} \x|^2 \dif r}
=\int_0^1 \left|\sigma_r^{\textcolor{black}{*}}\left[x_1+\int^1_rA_u^{\textcolor{black}{*}}x_2\dif u\right] \right|^2 \dif r
\asymp\int_0^1 \left|x_1+\int^1_rA_u^{\textcolor{black}{*}}x_2\dif u \right|^2 \dif r,
$$
and by the change of variable,
$$
\int_0^1 \left|x_1+\int^1_rA_u^{\textcolor{black}{*}}x_2\dif u \right|^2 \dif r
=\int_0^1 \left|x_1+\int^r_0A_{1-u}^{\textcolor{black}{*}}x_2\dif u \right|^2 \dif r.
$$
Since $A_u\in\cE$ and $\cE$ is a closed convex subset of $GL(\mR^d)$, we have for some $c_0\in(0,1)$,
$$
\inf_{A\in\cE}|A^{\textcolor{black}{*}}x_2|\geq c_0|x_2|\Rightarrow c_0r|x_2|\leq\left|\int^r_0A_u^{\textcolor{black}{*}}x_2\dif u\right|\leq c_0^{-1}r|x_2|.
$$
Recall $|x_1|^2+|x_2|^2=1$. If $|x_1|\leq \frac{c_0^{\textcolor{black}{3}}}{4}|x_2|$, then $|x_2|^2\geq (\frac{c_0^{\textcolor{black}{6}}}{16}+1)^{-1}$ and
\begin{align*}
\int_0^1 \left|x_1+\int^r_0A_u^{\textcolor{black}{*}}x_2\dif u \right|^2 \dif r
&\geq \int_0^1 \(|x_1|^2-2c_0^{-1}r|x_1||x_2|+c_0^{\textcolor{black}{2}}r^{\textcolor{black}{2}}|x_2|^2\)\dif r\\
&\geq c_0^{\textcolor{black}{2}} |x_2|^2\left(\int_{0}^1 \textcolor{black}{\big[r^2-\frac r2\big]}\dif r\right)\geq \frac{c_0^2}{12} \Big(\frac{c_0^{\textcolor{black}{6}}}{16}+1\Big)^{-1};
\end{align*}
if $|x_1|\geq \frac{c_0^{\textcolor{black}{3}}}{4}|x_2|$, then $|x_1|^2\geq(1+\frac{16}{c_0^{\textcolor{black}{6}}})^{-1}$ and
\begin{align*}
\int_0^1 \left|x_1+\int^r_0A_u^{\textcolor{black}{*}}x_2\dif u \right|^2 \dif r&\geq \int_0^{c^{\textcolor{black}{4}}_0/8} \left||x_1|-rc_0^{-1}|x_2|\right|^2 \dif r\ge |x_1|^2\int_0^{\textcolor{black}{c_0^4/8}} (1-r 4c_0^{-4})^2   \dif r\geq \frac{c_0^{\textcolor{black}{4}}}{32}\Big(1+\frac{16}{c_0^{\textcolor{black}{6}}}\Big)^{-1}.
\end{align*}
\textcolor{black}{We thus obtain \eqref{EG7}. As for \eqref{TT0}, it then readily follows from the scaling relation \eqref{TT2}}.
\end{proof}
\subsection{Control problem}

\textcolor{black}{In this subsection we show how the quantity $|\T^{-1}_{t-s}(\btheta_{t,s}(\x)-\y)|^2$ appearing in our main estimates can be, under \eqref{AX1} and \eqref{AX2}, related to a control problem (see \cite{dela:meno:10}) associated with $\F$ and $B$}.
More precisely, we consider the following deterministic \textcolor{black}{control} problem:
\begin{align}\label{ODE9}
\dot\bphi_{r,s}=\F(r,\bphi_{r,s})+B\varphi_r,\ \ r\in[s,t],\ \ \bphi_{s,s}=\x,\ \ \bphi_{t,s}=\y,
\end{align}
where $\varphi:[s,t]\to\mR^d$ is a \textcolor{black}{square integrable control} function.
Let $I(s,\x;t,\y)$ be the associated energy functional
$$
I(s,\x;t,\y)=\inf\left\{\left(\int^t_s|\varphi_r|^2\dif r\right)^{1/2},\ \ \bphi_{s,s}=\x,\ \ \bphi_{t,s}=\y\right\},
$$
where the infimum is taken \textcolor{black}{over all admissible controls $\varphi$}.
The following proposition plays a crucial role for proving the lower bound estimate of the heat kernel.
\bp\label{Le33_ACTION}
Under {\bf(H$^0_{\F}$)},
for any $T>0$, there exist constants $\kappa_5,\kappa_6\geq 1$ depending only on $T,\kappa_0,\kappa_1,d,\cE$ such that for all $0\leq s<t\leq T$ and $\x,\y\in\R^{2d}$,
\begin{align}\label{ES3}
\kappa_5^{-1} \(|\T^{-1}_{t-s}(\btheta_{t,s}(\x)-\y)|-1\)\leq I(s,\x;t,\y)\leq {\kappa_5} \(|\T^{-1}_{t-s}(\btheta_{t,s}(\x)-\y)|+1\).
\end{align}
Moreover, one can find a control $\varphi:[s,t]\to\mR^d$ and a solution $\phi_{r,s}$ to ODE \eqref{ODE9} such that 
\begin{align}\label{ES1}
\sup_{r\in[s,t]}|\varphi_r|\leq\kappa_6\(|\T^{-1}_{t-s}(\btheta_{t,s}(\x)-\y)|+1\)/\sqrt{t-s}.
\end{align}
\ep
\begin{proof}
Without loss of generality, we may assume $s=0$ and $t=t_0\leq\delta$, where $\delta$ is a small number only depending on
$T,\kappa_0,\kappa_1,d,\cE$. Let $(\bphi_t)_{t\in[0,t_0]}$ be any solution of control problem \eqref{ODE9}. Let
$$
\bpsi_t:=\bphi_t-\btheta_t,\ \ \btheta_t:=\btheta_{t,0}(\x).
$$
Then $\bpsi_t=(\psi_t^1,\psi_t^2)$ solves the following control problem:
\begin{align}\label{CN1}
\dot\bpsi_t=\bA(t,\bpsi_t)\psi_t+\wt\F(t,\bpsi_t)+B\varphi_t,\ \ \bpsi_0=0,\ \ \bpsi_{t_0}=\y-\btheta_{t_0},
\end{align}
where
$$
\wt\F(s,\x):=\begin{pmatrix}
F_1(s,\x+\btheta_t)-F_1(t,\btheta_t)\\
F_2(t,\btheta_t^1, x_2+\btheta_t^2)-F_2(t,\btheta^1_t,\btheta^2_t)
\end{pmatrix},
$$
and
$$
\bA(s,\x):=\begin{pmatrix}
0_{d\times d}&0_{d\times d}\\
\int^1_0\nabla_{x_1}F_2(t, u x_1+\btheta_t^1, x_2+\btheta^2_t)\dif u & 0_{d\times d} \end{pmatrix}.
$$
By \eqref{CN1} and \eqref{AX1}, we have
\begin{align*}
|\psi^1_t|\lesssim\int^t_0(1+|\psi_s|)\dif s+\int^t_0|\varphi_s|\dif s,
\end{align*}
and due to $\nabla_{x_1}F_2\in\cE$ and by \eqref{AX2},
\begin{align*}
|\psi^2_t|\lesssim\int^t_0|\psi^1_s|\dif s+\int^t_0(|\psi^2_s|^{\frac13}+|\psi^2_s|)\dif s.
\end{align*}
Thus by Lemma \ref{Le11} we have
\begin{align}\label{ES8}
|\psi^1_t|\lesssim t+\int^t_0|\varphi_s|\dif s,\ \ |\psi^2_t|\lesssim t^{3/2}+t\int^t_0|\varphi_s|\dif s,\ t\in[0,t_0].
\end{align}
Hence,
\begin{align*}
|\T^{-1}_{t_0}(\y-\btheta_{t_0})|&=t_0^{-\frac12}|\psi^1_{t_0}|+t_0^{-\frac32}|\psi^2_{t_0}|
\lesssim 1+t_0^{-\frac12}\int^{t_0}_0|\varphi_s|\dif s
\leq 1+\left(\int^{t_0}_0|\varphi_s|^2\dif s\right)^{1/2},
\end{align*}
which gives the left hand side estimate in \eqref{ES3}.

On the other hand, by Coron [Theorem 3.40], system \eqref{CN1} is controllable, and the exhibited control $\varphi$ is given by
\begin{align}\label{CN2}
\varphi_s=(\gR_{t_0,s}B)^*\K^{-1}_{t_0,0}\left(\y-\btheta_{t_0}-\int^{t_0}_0\gR_{t_0,s}\wt\F(s,\psi_s)\dif s\right),
\end{align}
where $\gR_{t,s}$ is the resolvent of $\bA(t,\psi_t)$ (see \eqref{Res11}), and
$$
\K_{t_0,0}:=\int^{t_0}_0 \gR_{t_0,s}BB^*\gR^*_{t_0,s}\dif s.
$$
Indeed, by Duhamel's formula, the above control $\varphi$ satisfies \eqref{CN1}.
Note that by \eqref{Res0},
$$
\gR_{t_0,s}B=\begin{pmatrix}\mI_{d\times d}\\
\int^{t_0}_s \int^1_0\nabla_{x_1}F_2(r, u \psi^1_r+\btheta_r^1, \psi^2_r+\btheta^2_r)\dif u\dif r \end{pmatrix}.
$$
By \eqref{TT2} and \eqref{TT1}, it is easy to see that
\begin{align*}
\sup_{s\in[0,t_0]}|\varphi_s|
&\lesssim t_0^{-\frac12}\left(|\T^{-1}_{t_0}(\y-\btheta_{t_0})|+\left|\T^{-1}_{t_0}\int^{t_0}_0\gR_{t_0,s}\wt\F(s,\psi_s)\dif s\right|\right)\\
&\stackrel{\eqref{AX2}}{\lesssim} 
t_0^{-\frac12}|\T^{-1}_{t_0}(\y-\btheta_{t_0})|+t_0^{-1}\int^{t_0}_0(1+|\psi_s|)\dif s+t_0^{-2}\int^{t_0}_0(|\psi^2_s|^{\frac13}+|\psi^2_s|)\dif s\\
&\stackrel{\eqref{ES8}}{\lesssim} t_0^{-\frac12}|\T^{-1}_{t_0}(\y-\btheta_{t_0})|+t_0^{-\frac12}+t_0\sup_{s\in[0,t_0]}|\varphi_s|
+t^{-\frac13}_0\sup_{s\in[0,t_0]}|\varphi_s|^{\frac13}\\
&\lesssim t_0^{-\frac12}|\T^{-1}_{t_0}(\y-\btheta_{t_0})|+\eps^{-1}t_0^{-\frac12}+(t_0+\eps)\sup_{s\in[0,t_0]}|\varphi_s|,
\end{align*}
where the last step is due to Young's inequality and the implicit constant only depends on $T,\kappa_0,\kappa_1,d,\cE$.
In particular, we can choose $\eps$ and $\delta$ small enough so that for all $t_0\in(0,\delta]$,
$$
\sup_{s\in[0,t_0]}|\varphi_s|\lesssim t_0^{-\frac12}\( |\T^{-1}_{t_0}(\y-\btheta_{t_0})|+1\).
$$
This in turn yields \eqref{ES1} as well as the right hand side estimate in \eqref{ES3}.
\end{proof}
\br\rm
When $\F$ is Lipschitz continuous, \cite{dela:meno:10} has shown that
$$
I(s,\x;t,\y)\asymp |\T^{-1}_{t-s}(\btheta_{t,s}(\x)-\y)|
$$
and
$$
\sup_{r\in[s,t]}|\varphi_r|\lesssim_C |\T^{-1}_{t-s}(\btheta_{t,s}(\x)-\y)|/\sqrt{t-s}.
$$
\textcolor{black}{The additional constant in the estimates of Proposition \ref{Le33_ACTION}, due to the rougher framework, will anyhow not perturb too much the analysis}.
\er
\section{Density bounds for SDEs with smooth coefficients} \label{SEC_DENS_SMOOTH}

In this section we always assume {\bf (H$^\gamma_\sigma$)} and {\bf (H$^{\gamma}_{\F}$)} for some $\gamma\in(0,1)$,
and temporarily assume that
\begin{align}\label{Smo1}
\|\nabla^j_\x \F\|_\infty<\infty,\ \ \|\nabla^j_\x \sigma\|_\infty<\infty,\ \ j\in\mN.
\end{align}
It is well known that
in the current smooth coefficients framework there exists a transition density $p(s,\x;t,\y)$ which is $C^\infty_b$- smooth in variables $\x,\y$ for all $s<t$, by H\"ormander's theorem.
Moreover, $p(s,\x;t,\y)$ satisfies the following backward Kolmogorov equation
 \begin{equation}\label{KOLM_DIFF}
\partial_s p(s,\x;t,\y)+\cL_{s,\x}p(s,\x;t,\y)=0, \quad
 {p}(s,\cdot; t,\y)\longrightarrow \delta_\y(\cdot)\; \text{weakly as}\;  s\uparrow t,
\end{equation}
and the forward Kolmogorov equation (Fokker-Planck equation):
\begin{equation}\label{FP_EQ_DIFF_SMOOTH_COEFF}
\partial_t p(s,\x;t,\y)\textcolor{black}{-}\cL^*_{(t,\y)}p(s,\x;t,\y)=0,\quad 
{p}(s,\x;t,\cdot)\longrightarrow \delta_\x(\cdot)\; \text{weakly as}\;  t\downarrow s,
\end{equation}
where, setting $a=\sigma\sigma^*/2$,
$$
\mathcal{L}_{s,\x}f(\x)=\tr \big(a(s,\x)\nabla_{x_1}^2 f(\x)\big)+\langle \gF(s,\x), \nabla_{\x} f(\x)\rangle,
$$
and 
$$
\mathcal{L}^*_{(t,\y)}f(\y)=\tr\(\nabla^2_{y_1}(a(t,\y)f(\y))\)-\div_{\y}(\gF(t,\y)f(\y)).
$$

The scope of the section is to obtain two-sided Aronson like bounds,
\green{where all the constants appearing below only depend on $\Theta_T$.}

\subsection{The Duhamel representation for $p(s,\x;t,\y)$}\label{SEC_DUHAMEL}

Fix now $(\t,\bxi)\in \R_+\times \R^{2d}$ as \textit{freezing parameters} to be chosen later on and let 
\begin{equation}
\dot \btheta_{t,\t}(\bxi)=\gF(t,\btheta_{t,\t}(\bxi)), \quad t\geq0, \quad \btheta_{\t,\t}(\bxi)=\bxi.
\end{equation}
We consider the stochastic linearized dynamics $(\wt \X_{t,s}^{(\tau,\bxi)})_{t\geq s} $:
\begin{equation}\label{FROZ}
\wt \X_{t,s}^{(\tau,\bxi)}=\x+\int^t_s [\gF(r,\btheta_{r,\t}(\bxi))+ \bA(r,\btheta_{r,\t}(\bxi))(\wt \X_{r,s}^{(\tau,\bxi)}-\btheta_{r,\t}(\bxi))]\dif r +
\int^t_sB\sigma(r,\btheta_{r,\t}(\bxi)) \dif W_r, 
 \end{equation} 
where, for all $\x=(x_1,x_2)\in \R^{2d}$, 
$$
\bA(r,\x)=\begin{pmatrix}0_{d\times d} & 0_{d\times d}\\
\nabla_{x_1}F_2(r,\x) & 0_{d\times d} \end{pmatrix}.
$$ 
Let $({\gR}^{(\tau,\bxi)}_{t,s})_{t\geq s}$ be the resolvent associated with $\bA(r,\btheta_{r,\t}(\bxi))$ (see \eqref{Res11}), which is explicitly given by
\begin{align}\label{Res00}
\gR^{(\tau,\bxi)}_{t,s}=\begin{pmatrix}\mI_{d\times d} & 0_{d\times d}\\
\int^t_s \nabla_{x_1}F_2(r,\btheta_{r,\t}(\bxi))\dif r & \mI_{d\times d} \end{pmatrix}.
\end{align}
If we define
\begin{align}\label{RR2}
\bvtheta_{t,s}^{(\tau,\bxi)}(\x):={\gR}^{(\tau,\bxi)}_{t,s}\x+
\int^t_s {\gR}^{(\tau,\bxi)}_{t,r}\Big(\gF(r,\btheta_{r,\t}(\bxi))-\bA(r,\btheta_{r,\t}(\bxi))\btheta_{r,\t}(\bxi)\Big)\dif r
\end{align}
and
\begin{align}\label{SB5}
\bTheta_r^{(\tau,\bxi)}:=B\sigma(r,\btheta_{r,\t}(\bxi)),
\end{align}
then by the variation formula of constants, $\wt  \X_{t,s}^{(\tau,\bxi)}(\x)$ is explicitly given by
 \begin{align}
\wt  \X_{t,s}^{(\tau,\bxi)}(\x) = \bvtheta_{t,s}^{(\tau,\bxi)}(\x)+ \int^t_s {\gR}^{(\tau,\bxi)}_{t,r}\bTheta_r^{(\tau,\bxi)} \dif W_r.\label{INTEGRATED}
 \end{align}
Clearly, the random variable $\wt  \X_{t,s}^{(\tau,\bxi)}(\x)$
 admits a Gaussian density $\wt p^{(\tau,\bxi)}(s,\x;t,\cdot) $ given by
\begin{equation}\label{CORRESP}
\wt p^{(\tau,\bxi)}(s,\x;t,\y)=\frac{1}{(2\pi)^{d}\det(\K_{t,s}^{(\tau,\bxi)})^{\frac 12}}
\exp\left( -\frac12\big|(\K_{t,s}^{(\tau,\bxi)})^{-\frac12} (\bvtheta_{t,s}^{(\tau,\bxi)}(\x)-\y)\big|^2\right),
\end{equation}
where
\begin{equation}\label{KK1}
\K_{t,s}^{(\tau,\bxi)}:=\int^t_s {\gR}^{(\tau,\bxi)}_{t,r}\bTheta_r^{(\tau,\bxi)}({\gR}^{(\tau,\bxi)}_{t,r}\bTheta_r^{(\tau,\bxi)})^* \dif r.
\end{equation}
In particular, $\wt p^{(\tau,\bxi)}(s,\x;t,\y)$ satisfies
\begin{equation}\label{Kolmogorov_frozen}
\partial_s \wt p^{(\tau,\bxi)}(s,\x;t,\y)+\wt \cL^{(\tau,\bxi)}_{s,\x}\wt p^{(\tau,\bxi)}(s,\x;t,\y)=0, \quad
\wt p^{(\tau,\bxi)}(s,\x;t,\y)\longrightarrow \delta_\y(\cdot)\; \text{weakly as}\;  s\uparrow t,
\end{equation}
where,  for $a=\sigma\sigma^*/2$,
\begin{equation}\label{frozen_gen}
\wt\cL^{(\tau,\bxi)}_{s,\x}=\tr\big(a(s,\btheta_{s,\tau}(\bxi))\cdot \nabla^2_{x_1}\big)+
\langle \left(\gF(s,\btheta_{s,\t}(\bxi))+ \bA(s,\btheta_{s,\t}(\bxi))(\x-\btheta_{s,\t}(\bxi))\right),\nabla_{\x}\rangle
\end{equation}
denotes the generator of the diffusion with frozen coefficients in \eqref{FROZ}.\\

The following proposition is a direct consequence of expression \eqref{CORRESP} and Lemma \ref{Le25}.

\bp[A priori controls for the frozen Gaussian densities]\label{PROP_Proxy}
Under {\bf(H$^0_\sigma$)}, {\bf(H$^{\0}_{\F}$)} and \eqref{Smo1}, for any $T>0$
and $j=(j_1,j_2)\in\mN_0^2$, there are constants
$\lambda_j, C_j\geq 1$ depending only on $\Theta_T$ such that
for all $0\leq s<t\leq T$, $\tau\in[0,T]$ and $\x,\y,\bxi\in\mR^{2d}$,
 \begin{align}
C_0^{-1}g_{\lambda_0^{-1}}\(t-s,\bvtheta^{(\tau,\bxi)}_{t,s}(\x)-\y\)\leq
\wt p^{(\tau,\bxi)}(s,\x;t,\y)\leq C_0g_{\lambda_0}\(t-s,\bvtheta^{(\tau,\bxi)}_{t,s}(\x)-\y\),\label{lower_proxy}
\end{align}
where $g_\lambda(t,\x)$ is defined by \eqref{GG2}, and
\begin{align}
|\nabla^{j_1}_{x_1}\nabla^{j_2}_{x_2}\wt p^{(\tau,\bxi)}(s,\x;t,\y)|&\leq C_j(t-s)^{-\frac{j_1+3j_2}{2}}g_{\lambda_j}\(t-s,\bvtheta^{(\tau,\bxi)}_{t,s}(\x)-\y\). \label{upper_proxy}
\end{align}
\ep

The starting point of our analysis is the following Duhamel type representation formula which readily follows in the current \textit{smooth coefficients} setting from \eqref{KOLM_DIFF}-\eqref{FP_EQ_DIFF_SMOOTH_COEFF} and \eqref{frozen_gen}:
\begin{align}
p(s,\x;t,\y)&=\wt{p}^{(\tau,\bxi)}(s,\x;t,\y)+\int^t_s\int_{\R^{2d}}\wt{p}^{(\tau,\bxi)}(s,\x;r,\z)(\cL_{r,\z}-\wt \cL^{(\tau,\bxi)}_{r,\z})p(r,\z; t,\y)\dif\z \dif r \label{D0}
\\&=\wt{p}^{(\tau,\bxi)}(s,\x;t,\y)+\int^t_s\int_{\R^{2d}}p(s,\x;r,\z)(\cL_{r,\z}-\wt\cL^{(\tau,\bxi)}_{r,\z}){\wt p}^{(\tau,\bxi)}(r,\z; t,\y)\dif\z \dif r,\label{D1}
\end{align}
If we take $(\tau,\bxi)=(s,\x)$ in \eqref{D0} and set $\wt p_0(s,\x;t,\y):=\wt{p}^{(s,\x)}(s,\x;t,\y)$, then we obtain the {\it backward} representation
$$
p(s,\x;t,\y)=\wt{p}_0(s,\x;t,\y)+\int^t_s\int_{\R^{2d}}\wt{p}_0(s,\x;r,\z)(\cL_{r,\z}-\wt\cL^{(s,\x)}_{r,\z})p(r,\z; t,\y)\dif\z \dif r. 
$$
If we take $(\tau,\bxi)=(t,\y)$ in \eqref{D1} and set $\wt p_1(s,\x;t,\y):=\wt{p}^{(t,\y)}(s,\x;t,\y)$, we then obtain the {\it forward} representation
$$
p(s,\x;t,\y)=\wt{p}_1(s,\x;t,\y)+\int^t_s\int_{\R^{2d}}p(s,\x;r,\z)(\cL_{r,\z}-\wt\cL^{(t,\y)}_{r,\z}){\wt p}_1(r,\z; t,\y)\dif\z \dif r.
$$

To give the estimates of $\wt{p}_i(s,\x;t,\y)$, $i=0,1$, we need the following lemma.
\bl\label{Le32}
For any $t\geq s$ and $\x,\y\in\mR^{2d}$, it holds that
$$
\bvtheta_{t,s}^{(s,\x)}(\x)=\btheta_{t,s}(\x),\ \ \bvtheta_{t,s}^{(t,\y)}(\x)-\y={\gR}^{(t,\y)}_{t,s}(\x-\btheta_{s,t}(\y)).
$$
Moreover, under {\bf (H$^{\gamma}_{\F}$)}, for any $T>0$, there is a constant $C:=C(\Theta_T)>0$ such that for all $0\leq s\leq t\leq T$ and $\x,\y\in\mR^{2d}$,
\begin{align}\label{TT8}
|\T^{-1}_{t-s}(\btheta_{t,s}(\x)-\bvtheta_{t,s}^{(t,\y)}(\x))|\leq C (t-s)^{\frac\gamma2}\(|\T^{-1}_{t-s}(\btheta_{t,s}(\x)-\y)|^{1+\gamma}+1\).
\end{align}
\el
\begin{proof}
(i) Since $\bvtheta_{t,s}^{(\tau,\bxi)}(\x)=\mE(\wt \X_{t,s}^{(\tau,\bxi)})$, by \eqref{FROZ} one sees that
\begin{align}\label{EQ1}
\bvtheta_{t,s}^{(\tau,\bxi)}(\x)=\x+\int^t_s [\gF(r,\btheta_{r,\tau}(\bxi))+ \bA(r,\btheta_{r,\tau}(\bxi))(\bvtheta_{r,s}^{(\tau,\bxi)}(\x)-\btheta_{r,\tau}(\bxi))]\dif r.
\end{align}
Hence,
\begin{align*}
|\bvtheta_{t,s}^{(s,\x)}(\x)-\btheta_{t,s}(\x)|
&\leq \int^t_s|\bA(r,\btheta_{r,s}(\x))(\bvtheta_{r,s}^{(s,\x)}(\x)-\btheta_{r,s}(\x))|\dif r\\
&\leq\|\nabla_{x_1}F_2\|_\infty\int^t_s|\bvtheta_{r,s}^{(s,\x)}(\x)-\btheta_{r,s}(\x)|\dif r,
\end{align*}
which implies by Gronwall's inequality that
$$
\bvtheta_{t,s}^{(s,\x)}(\x)=\btheta_{t,s}(\x),
$$
\textcolor{black}{which gives the first equality of the lemma}.

(ii) For the second one, note that  by \eqref{RR2} and \eqref{RR1},
\begin{align*}
{\gR}^{(t,\y)}_{s,t}\bvtheta_{t,s}^{(t,\y)}(\x)=\x+
\int^t_s {\gR}^{(t,\y)}_{s,r}\Big(\gF(r,\btheta_{r,t}(\y))-\bA(r,\btheta_{r,t}(\y))\btheta_{r,t}(\y)\Big)\dif r.
\end{align*}
Since $\partial_{s}\gR^{(t,\y)}_{s,t}=-\bA(s,\btheta_{s,t}(\y))\gR^{(t,\y)}_{s,t}$, by Duhamel's formula, we also have
\begin{align*}
\Gamma_s(\y)&:={\gR}^{(t,\y)}_{s,t}\y-\int^t_s {\gR}^{(t,\y)}_{s,r}\Big(\gF(r,\btheta_{r,t}(\y))-\bA(r,\btheta_{r,t}(\y))\btheta_{r,t}(\y)\Big)\dif r\\
&=\y+\int^t_s [\gF(r,\btheta_{r,t}(\y))- \bA(r,\btheta_{r,t}(\y))(\Gamma_r(\y)-\btheta_{r,t}(\y))]\dif r.
\end{align*}
As above, one has $\Gamma_s(\y)=\btheta_{s,t}(\y)$. Hence,
$$
{\gR}^{(t,\y)}_{s,t}\bvtheta_{t,s}^{(t,\y)}(\x)-\x={\gR}^{(t,\y)}_{s,t}\y-\btheta_{s,t}(\y).
$$

(iii) \textcolor{black}{Let us now turn to the proof of \eqref{TT8}}. Fix $0\leq s<u\leq t\leq T$. Note that by \eqref{EQ1},
\begin{align*}
\btheta_{t,s}(\x)-\bvtheta_{t,s}^{(t,\y)}(\x)=\int^t_s [\gF(r,\btheta_{r,s}(\x))-\gF(r,\btheta_{r,t}(\y))-\bA(r,\btheta_{r,t}(\y))(\bvtheta_{r,s}^{(t,\y)}(\x)-\btheta_{r,t}(\y))]\dif r,
\end{align*}
and by \eqref{EZ1},
\begin{align}\label{EZ11}
\begin{split}
&(t-s)^{-\frac12}|(\btheta_{r,s}(\x)-\btheta_{r,t}(\y))_1|+(t-s)^{-\frac32}|(\btheta_{r,s}(\x)-\btheta_{r,t}(\y))_2|\\
&\qquad\qquad\lesssim |\T^{-1}_{t-s}(\btheta_{t,s}(\x)-\y)|+1=:\cA.
\end{split}
\end{align}
Then we have
\begin{align*}
|(\btheta_{t,s}(\x)-\bvtheta_{t,s}^{(t,\y)}(\x))_1|
&\leq\int^t_s|F_1(r,\btheta_{r,s}(\x))-F_1(r,\btheta_{r,t}(\y))|\dif r\\
&\lesssim\int^t_s(1+|\btheta_{r,s}(\x)-\btheta_{r,t}(\y)|)\dif r\lesssim\green{(t-s)}\cA,
\end{align*}
and by definition,
\begin{align*}
|(\btheta_{t,s}(\x)-\bvtheta_{t,s}^{(t,\y)}(\x))_2|&\leq\int^t_s
\Big[\|\nabla_{x_1}F_2\|_\infty|(\btheta_{r,s}(\x)-\bvtheta_{r,s}^{(t,\y)}(\x))_1|+|\cT_{F_2(r)}(\btheta_{r,s}(\x),\btheta_{r,t}(\y))|\Big]\dif r\\
&\stackrel{\eqref{Taylor}}{\lesssim}\int^t_s\Big[|(\btheta_{r,s}(\x)-\bvtheta_{r,s}^{(t,\y)}(\x))_1|+|\btheta_{r,s}(\x)-\btheta_{r,t}(\y)|_{\bf d}^{1+\gamma}\Big]\dif r\\
&\stackrel{\eqref{EZ11}}{\lesssim}(t-s)^2\cA+ (t-s)^{\frac{3+\gamma}{2}}\cA^{1+\gamma}.
\end{align*}
Thus we obtain \eqref{TT8}. The proof is complete.
\end{proof}

For any $\lambda>0$, $0\leq s<t<\infty$ and $\x,\y\in\mR^{2d}$, let
\begin{align}\label{GG1}
\hat p_\lambda(s,\x;t,\y)&:=g_{\lambda}\(t-s,\btheta_{t,s}(\x)-\y\)=(t-s)^{-2d}\e^{-|\T^{-1}_{t-s}(\btheta_{t,s}(\x)-\y)|^2/(2\lambda)},
\end{align}
\green{with $\btheta_{t,s}(\x)$ the flow associated with $\gF$}.
\bl\label{Le33}
For any $T,\lambda>0$ and $\alpha>0$, there are $\lambda'>\lambda$ and $C>0$ depending only on $\Theta_T$ and $\lambda,\alpha$ such that 
for all $0\leq s<t\leq T$ and $\x,\y\in\mR^{2d}$,
\begin{align}
\(|\btheta_{t,s}(\x)-\y|_{\bf d}^{\alpha}+|\x-\btheta_{s,t}(\y)|_{\bf d}^{\alpha}\)\hat{p}_\l (s,\x;t,\y)&\le C (t-s)^{\frac\alpha 2}\hat{p}_{\l'}(s,\x;t,\y).\label{regularizing_proxy}
\end{align}
\el
\begin{proof}
Note that by \eqref{EQ_EQUIV_FLOW},
$$
g_{\lambda/\kappa_3}\(t-s,\x-\btheta_{s,t}(\y)\)\lesssim\hat p_\lambda(s,\x;t,\y)\lesssim g_{\kappa_3\lambda}\(t-s,\x-\btheta_{s,t}(\y)\).
$$
The desired estimate now follows by definition.
\end{proof}
We also need the following convolution type inequality for the Gaussian functions $\hat p_\l$. 

\bl[Reproduction property]
Under \eqref{AX1} and \eqref{AX2}, for any $T>0$, there are constants $\kappa_7, C_3\geq 1$ depending only on $\Theta_T$ such that for all $\x,\y\in\R^{2d}$, $0\leq s<r\leq t\leq T$
and $\l_1, \l_2>0$,
\begin{equation}\label{Convolution}
C_3^{-1}\hat{p}_{(\l_1\wedge \l_2)/\kappa_7}(s,\x;t,\y)\leq \int_{\R^{2d}}\hat{p}_{\l_1}(s,\x;r,\z)\hat{p}_{\l_2}(r,\z;t,\y)\dif\z\leq C_3\hat{p}_{\kappa_7(\l_1\vee\l_2)}(s,\x;t,\y).
\end{equation}
\el
\proof 
Recalling \eqref{GG1} and \eqref{GG2}, we first consider the convolution 
$$I:=\int_{\R^{2d}}g_{\l_1}(\epsilon_1,\x'-\z)g_{\l_2}(\epsilon_2,\z-\x'')\dif\z.$$
Clearly, up to a constant dependent on $d$ and $\l_1$, $\l_2$ , $I$ is the density of the sum of two independent Gaussian vectors 
with means $\x'$ and $\x''$ respectively. 
Therefore $I$ is a Gaussian function with mean $\x'-\x''$ and covariance matrix 
$\l_1\T^{-2}_{\epsilon_1}+\l_2\T^{-2}_{\epsilon_2}$.
We have 
$$
\tfrac{\l_1\wedge \l_2}8\T^{-2}_{\epsilon_1+\epsilon_2}
\leq \lambda_1\T^{-2}_{\epsilon_1}+\lambda_2\T^{-2}_{\epsilon_2}
\leq (\l_1\vee \l_2)\T^{-2}_{\epsilon_1+\epsilon_2},
$$
which yields 
\begin{equation}\label{Convolution2}
g_{(\l_1\wedge \l_2)/8}(\epsilon_1+\epsilon_2,\x'-\x'')\lesssim I\lesssim g_{\l_1\vee \l_2}(\epsilon_1+\epsilon_2,\x'-\x'').
\end{equation}
To prove \eqref{Convolution}, it suffice\green{s} to notice that by Lemma \ref{lemme:bilipflow},
$$g_{\l_2/\kappa_3}(t-r,\z-\btheta_{r,t}(\y))\lesssim \hat{p}_{\l_2}(r,\z; t,\y)\lesssim g_{\kappa_3\l_2}(t-r,\z-\btheta_{r,t}(\y)).$$
\textcolor{black}{The claim then follows from} \eqref{Convolution2} with $\x'=\btheta_{r,s}(\x)$ and $\x''=\btheta_{r,t}(\y)$ and using \eqref{EZ1}.
\endproof

The following lemma is a direct consequence of Proposition \ref{PROP_Proxy} and Lemma \ref{Le32}.
\bl\label{Le34}
Under {\bf(H$^0_\sigma$)}, {\bf(H$^0_{\F}$)} and \eqref{Smo1}, for any $T>0$
and $j=(j_1,j_2)\in\mN_0^2$, there are constants
$\lambda_j,C_j\geq 1$ depending only on $\Theta_T$ such that
for all $0\leq s<t\leq T$  and $\x,\y\in\mR^{2d}$,
\begin{align}
C_0^{-1}\hat p_{\lambda^{-1}_0}\(s,\x;t,\y)\leq
\wt p_i(s,\x;t,\y)\leq C_0\hat p_{\lambda_0}\(s,\x;t,\y),\ i=0,1,\label{11}
\end{align}
and
\begin{align}
|\nabla^{j_1}_{x_1}\nabla^{j_2}_{x_2}\wt p^{(\tau,\bxi)}(s,\x;t,\y)|_{(\tau,\bxi)=(s,\x)\, \mathrm{ or }\, (t,\y)}&\leq C_j(t-s)^{-\frac{j_1+3j_2}{2}}\hat p_{\lambda_j}\(s,\x;t,\y). \label{22}
\end{align}
Moreover, for any $\alpha\in [0,1]$, \textcolor{black}{recalling $\wt p_1(s,\x;t,\y):=\wt{p}^{(t,\y)}(s,\x;t,\y)$},
 \begin{align}
&|\nabla_{x_1}^{j_1} \nabla_{x_2}^{j_2}\left(\wt p_1(s,\x;t,\y)-\wt p_1(s,\x'; t,\y)\right)|\nonumber \\
& \qquad \quad \leq C_j |\x-\x'|_{\bf d}^{\alpha}(t-s)^{-\frac{\alpha + j_1+3j_2}{2}}
\left(\hat{p}_{\lambda_j}(s,\x;t,\y)+\hat{p}_{\lambda_j}(s,\x'; t,\y)\right). \label{proxy_regularity}
\end{align}
\el
\begin{proof}
Estimates \eqref{11} and \eqref{22} are direct by \eqref{lower_proxy}, \eqref{upper_proxy}, Lemma \ref{Le32} and \eqref{TT0}.
Here we prove \eqref{proxy_regularity} with $j=(1,0)$ for simplicity since the general case is analogous. 
The off-diagonal regime $|\x-\x'|_{\bf d}> (t-s)^\frac{1}{2}$ is straightforward from \eqref{22}. 
Assume now $|\x-\x'|_{\bf d}\leq (t-s)^\frac{1}{2}$. Then we have
\begin{align*}
|\nabla_{x_1}\left(\wt p_1(s,\x;t,\y)-\wt p_1(s,\x'; t,\y)\right)|
& \le \sum_{i=1}^2|x_i-x'_i|\sup_{\eta\in [0,1]}|\nabla_{x_i}\nabla_{x_1}\wt p_1(s,\x+\eta(\x'-\x); t,\y)|\\
& \lesssim \sum_{i=1}^2 |x_i-x'_i|(t-s)^{-i}\sup_{\eta\in [0,1]}\hat{p}_{\l_2}(s,\x+\eta(\x'-\x); t,\y)\\
&\lesssim|\x-\x'|_{\bf d}^{\alpha}\,(t-s)^{-\frac{\alpha + 1}{2}}\hat{p}_{\lambda_j}(s,\x;t,\y),
\end{align*}
where in the last inequality we have used $|\x-\x'|_{\bf d}\leq (t-s)^\frac{1}{2}$ and \eqref{VX2}.
The proof is complete.
\end{proof}


\subsection{Derivation of the upper bound}\label{UPPER}
For notational convenience, we write from now on 
\begin{align}\label{HH0}
\begin{split}
&H(s,\x;t,\y):=(\cL_{s,\x}-\wt\cL^{(t,\y)}_{s,\x}){\wt p}_1(s,\x;t,\y)=(a(s,\x)-a(s,\btheta_{s,t}(\y)))\cdot\nabla^2_{x_1}{\wt p}_1(s,\x;t,\y)\\
&\quad+(F_1(s,\x)- F_1(s,\btheta_{s,t}(\y)))\cdot \nabla_{x_1}{\wt p}_1(s,\x;t,\y)+\cT_{F_2(s)}(\x,\btheta_{s,t}(\y))\cdot \nabla_{x_2}{\wt p}_1(s,\x;t,\y),
\end{split}
\end{align}
where $\cT_{F_2(s)}$ is defined by \eqref{Taylor}, and
\begin{equation}\label{DE1}
(p\otimes H)(s,\x;t,\y)=\int^t_s\int_{\R^d}p(s,\x;r,\z)H(r,\z; t,\y)\dif\z \dif r .
\end{equation}
Thus, from the Duhamel representation \eqref{D1}, we have 
\begin{align}
p(s,\x;t,\y)=\wt{p}_1(s,\x;t,\y)+(p\otimes H)(s,\x;t,\y).\label{FR}
\end{align}
For $N\geq 2$, iterating $N-1$-times the identity \eqref{FR}, we obtain
\begin{align}
p(s,\x;t,\y)=\wt{p}_1(s,\x;t,\y)+\sum_{j=1}^{N-1}(\wt p_1\otimes H^{\otimes j})(s,\x;t,\y)+(p\otimes H^{\otimes N})(s,\x;t,\y).\label{FR1}
\end{align}

\bp[Control of the parametrix expansion]\label{Conv_Kernels}
For any $T>0$ and $N\in\N$, there are constants $C_N,\lambda_N>0$ depending only on $\Theta_T$ such that for all $\x,\y\in\R^{2d}$ and $0\leq s<t\leq T$,
\begin{equation}\label{ES6}
|H^{\otimes N}(s,\x;t,\y)|\leq C_N(t-s)^{-1+\frac{N {\gamma}}{2}} \hat p_{\lambda_N}(s,\x;t,\y),
\end{equation}
where $\lambda_N\to \infty$ as $N\to\infty$. In particular,
\begin{equation*}
p(s,\x;t,\y)\leq C_{N-1} \hat{p}_{\l_{N-1}}(s,\x;t,\y)+|(p\otimes H^{\otimes N})(s,\x;t,\y)|. 
\end{equation*}
\ep
\proof
By \eqref{HH0}, \eqref{upper_proxy}, \eqref{Taylor} and \eqref{regularizing_proxy}, we have
\begin{align*}
|H(s,\x;t,\y)|&\leq|a(s,\x)-a(s,\btheta_{s,t}(\y))|\cdot|\nabla^2_{x_1}{\wt p}_1(s,\x;t,\y)|\\
&\quad +|F_1(s,\x)- F_1(s,\btheta_{s,t}(\y))|\cdot |\nabla_{x_1}{\wt p}_1(s,\x;t,\y)|\\
&\quad +|\cT_{F_2(s)}(\x,\btheta_{s,t}(\y))|\cdot |\nabla_{x_2}{\wt p}_1(s,\x;t,\y)|\\
&\lesssim (t-s)^{-1}|\x-\btheta_{s,t}(\y)|_{\bf d}^{\gamma} \hat p_\lambda(s,\x;t,\y) \\
&\quad+ (t-s)^{-\frac 12}(1+|\x-\btheta_{s,t}(\y)|) \hat p_\lambda(s,\x;t,\y) \\
&\quad + (t-s)^{-\frac 32}|\x-\btheta_{s,t}(\y)|_{\bf d}^{1+\gamma} \hat p_\lambda(s,\x;t,\y)\\
&\lesssim (t-s)^{-1+\frac{ \gamma}{2}}\hat p_{\lambda'}(s,\x;t,\y).
\label{FIRST_CTR_H}
\end{align*}
This gives the stated estimate for $N=1$. 
For general $N\geq 2$, by induction, it is readily seen that:
\begin{align*}
|H^{\otimes N}(s,\x;t,\y)|&\lesssim\int^t_s (r-s)^{-1+\frac{(N-1)\gamma}{2}}(t-r)^{-1+\frac{ \gamma}{2}}
\left(\int_{\mR^{2d}}\hat{p}_{\l_{N-1}}(s,\x;r,\z)\hat{p}_{\l_1}(r,\z;t,\y)\dif z \right) \dif r\\
&\stackrel{\eqref{Convolution}}{\lesssim}\left(\int^t_s (r-s)^{-1+\frac{(N-1)\gamma}{2}}(t-r)^{-1+\frac{ \gamma}{2}}\dif r\right)\hat{p}_{\l_N}(s,\x;t,\y)\\
&\lesssim (t-s)^{-1+\frac{N\gamma}{2}} \hat{p}_{\l_N}(s,\x;t,\y).
\end{align*}
This completes the proof.
\endproof

We carefully remark that we have to stop the iteration at some fixed $N$ to avoid the explosion of $\l_N$ as $N$ goes to infinity.  
The following proposition provides a control for the remainder and concludes the proof of the upper bound.

\bp[Control of the remainder]\label{PROP_CONTR}
Let $N$ be large enough such that $$-1+\frac{N \gamma}{2}>2d.$$
Then there exist constants $C_0, \l_0>0$ such that for all $\x,\y\in\R^{2d}$ and $0\le s<t\le T$,
\begin{equation}
|(p\otimes H^{\otimes N})(s,\x;t,\y)|\le C_0\hat{p}_{\l_0}(s,\x;t,\y).
\end{equation}
\ep

\begin{proof}
We first recall that, from
 the intrinsic scaling properties of SDE \eqref{SDE0} we can restrict w.l.o.g to the case $s=0,t=1$ for the proof (we refer to \cite{dela:meno:10} for additional details).
 
Indeed, recall first from \eqref{DEF_T} that for any $\lambda>0$,  the \textit{intrinsic scale matrix} writes:
\begin{align*}
\T^{-1}_\lambda:=\begin{pmatrix} 
\lambda^{-\frac12}\mI_{d\times d} & 0_{d\times d}\\
0_{d\times d} & \lambda^{-\frac32}\mI_{d\times d}
\end{pmatrix}.
\end{align*}
Fix then $s\geq 0$ and $\lambda>0$. We define 
\begin{align}\label{AX8}
\X^{\l,s}_t:=\T^{-1}_\lambda\X_{\l {t+s}}, \ t\geq 0.
\end{align}
By the change of variables, it is easy to see that $(\X^{\l,s}_t)_{t\geq 0}$ satisfies 
\begin{align*}
\dif\X^{\l,s}_t=\gF^{\lambda,s}(t, \X^{\l,s}_t)\dif t+B\sigma^{\lambda,s}(t,\X^{\l,s}_t)\dif W^{\l}_t,
\end{align*}
where $W^{\l}_t:=\l^{-\frac 12}W_{\l t}$ is a still Brownian motion, and 
\begin{align}\label{FF1}
\gF^{\lambda,s}(t,\x):=\l\T^{-1}_\l \gF\(s+\l t, \T_\l \x\),\ \ \sigma^{\lambda,s}(t,\x):=\sigma\(s+\l t, \T_\l \x\).
\end{align}
In particular, let $p(s,\x;t,\y)$ (resp. $p^{\lambda,s}(\x;t,\y)$) be the density of $\X_t$ (resp. $\X_t^{\lambda,s}$) starting from $\x$ at time $s$ (resp. $0$). Then
\begin{align}\label{SCL2}
p(s,\x;t,\y)=\lambda^{2d}p^{\lambda,s}(\T^{-1}_\lambda\x;\lambda^{-1}t,\T^{-1}_\lambda\y).
\end{align}

From the scaling property \eqref{SCL2}, it thus suffices to consider the case $(s,t)=(0,1)$. 
First of all, by \eqref{ES6}, we have
\begin{align*}
\sI:=|(p\otimes H^{\otimes N})(0,\x;1,\y)|
&\lesssim\int^1_0(1-r)^{-1+\frac{N {\gamma}}{2}}\int_{\mR^{2d}}p(0,\x;r,\z) \hat p_{\lambda_N}(r,\z;1,\y)\dif \z\dif r\\
&=\int^1_0(1-r)^{-1+\frac{N {\gamma}}{2}}\mE\hat p_{\lambda_N}(r,\X_{r,0}(\x);1,\y)\dif r.
\end{align*}
Noting that by Lemma 2.8 in \cite{MPZ20},
\begin{align*}
\mE\hat p_{\lambda_N}(r,\X_{r,0}(\x);1,\y)
&\leq C_1\sup_{\z\in\mR^{2d}}\exp\Big\{\ln \hat p_{\lambda_N}(r,\z;1,\y)-C_2|\z-\btheta_{r,0}(\x)|^2\Big\},
\end{align*}
and
\begin{align*}
\ln \hat p_{\lambda_N}(r,\z;1,\y)
&\stackrel{\eqref{GG1}}{=}\ln (1-r)^{-2d}-\lambda_N|\T^{-1}_{1-r}(\btheta_{1,r}(\z)-\y)|^2\\
&\stackrel{\eqref{EQ_EQUIV_FLOW}}{\leq} \ln (1-r)^{-2d}-\lambda'_N|\T^{-1}_{1-r}(\z-\btheta_{r,1}(\y))|^2+C_3\\
&\,\,\leq \ln (1-r)^{-2d}-\lambda'_N|\z-\btheta_{r,1}(\y)|^2+C_3,
\end{align*}
we further have
\begin{align*}
\mE\hat p_{\lambda_N}(r,\X_{r,0}(\x);1,\y)
&\lesssim(1-r)^{-2d}\exp\left\{-\lambda''_N\inf_{\z\in\mR^{2d}}\Big\{|\z-\btheta_{r,1}(\y)|^2+|\z-\btheta_{r,0}(\x)|^2\Big\}\right\}\\
&=(1-r)^{-2d}\exp\left\{-\textcolor{black}{\frac{\lambda''_N}{2}}|\btheta_{r,1}(\y)-\btheta_{r,0}(\x)|^2\right\}\\
&\stackrel{\eqref{EQ_EQUIV_FLOW}}{\leq} (1-r)^{-2d}\exp\left\{-\lambda'''_N|\btheta_{1,0}(\x)-\y|^2\right\}.
\end{align*}
Hence,
\begin{align*}
\sI\lesssim \left(\int^1_0(1-r)^{-1+\frac{N {\gamma}}{2}-2d} \dif r \right)\e^{-\lambda'''_N|\btheta_{1,0}(\x)-\y|^2}\lesssim \e^{-\lambda'''_N|\btheta_{1,0}(\x)-\y|^2}.
\end{align*}
The proof is complete.
\end{proof}

\subsection{Derivation of the lower bound}\label{LOWER}

We first derive a local bound, starting from the one step parametrix expansion \eqref{FR}: 
by \eqref{upper_proxy} and the upper bound for $p$, we have 
\begin{align*}
p(s,\x;t,\y)&=\wt{p}_1(s,\x;t,\y)+(p\otimes H)(s,\x;t,\y)\\
&\geq c_0\hat{p}_{\l^{-1}_0}(s,\x;t,\y)-\int^t_s\int_{\R^{2d}}p(s,\x;r,\z)|H(r,\z; t,\y)|\dif\z \dif r\\
&\geq c_0\hat{p}_{\l^{-1}_0}(s,\x;t,\y)-C\int^t_s\int_{\R^{2d}}\hat{p}_{\l}(s,\x;r,\z)(t-r)^{\frac{\gamma}{2}-1}\hat{p}_{\l}(r,\z;t,\y)\dif\z \dif r\\\
&\stackrel{\eqref{Convolution}}{\geq} C^{-1}\hat{p}_{\l^{-1}}(s,\x;t,\y)-C(t-s)^{\frac{\gamma}{2}}\hat{p}_{\l}(s,\x;t,\y),
\end{align*}
for sufficiently large $\l$ and $C$, depending on $\Theta_T$.
Suppose that 
$$
|\T^{-1}_{t-s}(\btheta_{t,s}(\x)-\y)|\le \kappa_5(2\kappa_6+1)+2\kappa_3+1=:\kappa_9,
$$
where $\kappa_3,\kappa_5$ and $\kappa_6$ are from \eqref{EQ_EQUIV_FLOW}, \eqref{ES3} and \eqref{ES1}.
Then
\begin{equation}
p(s,\x;t,\y)\geq (t-s)^{-2d}(C^{-1}\e^{-\l \kappa_9^2 }-C(t-s)^{\frac{\gamma}{2}})\geq \frac{\e^{-\l\kappa_9^2 }}{2C}(t-s)^{-2d}=C_0(t-s)^{-2d}, \label{loc_lower}
\end{equation}
provided that 
$$
t-s \leq (2C^2)^{-\frac{2}{\gamma}}\e^{-\frac{2\l\kappa_9^2 }{\gamma}}=:T_0.
$$
Now for fixed $0\leq s<t\leq T$ with $t-s\leq T_0$ and $\x,\y\in\mR^{2d}$, we assume
$$
|\T^{-1}_{t-s}(\btheta_{t,s}(\x)-\y)|\geq \kappa_9.
$$
Let $M$ be the smallest integer such that
\begin{align}\label{DP1}
M-1\leq |\T^{-1}_{t-s}(\btheta_{t,s}(\x)-\y)|^2<M.
\end{align}
Define
$$
\delta:=\tfrac{t-s}{M},\ \  t_j:=s+j\delta,\ \ j=0,1,\cdots,M.
$$
Let $\bphi_{t,s}(\x)$ be the optimal curve in \textcolor{black}{Proposition} \ref{Le33_ACTION} so that the corresponding control $\varphi$ has the estimate
$$
\sup_{r\in[s,t]}|\varphi_r|\leq \kappa_6\( |\T^{-1}_{t-s}(\btheta_{t,s}(\x)-\y)|+1\)/\sqrt{t-s}\leq \kappa_6(\sqrt{M}+1)/\sqrt{t-s}.
$$ 
Define
$$
\bxi_j:=\bphi_{t_j,s}(\x),\ \ j=0,1,\cdots,M.
$$
For any $j=0,\cdots,M-1$, by Proposition \ref{Le33_ACTION}, we have
\begin{align}\label{EG5}
\begin{split}
|\T^{-1}_\delta\(\btheta_{t_{j+1},t_j}(\bxi_j)-\bxi_{j+1}\)|&\lesssim_{\kappa_5}
I(t_j,\bxi_j;t_{j+1}, \bxi_{j+1})+1\leq \left(\int^{t_{j+1}}_{t_j}|\varphi_r|^2\dif r\right)^{1/2}+1\\
&\leq(t_{j+1}-t_j)^{1/2}\sup_{r\in[s,t]}|\varphi_r|+1\leq 2\kappa_6+1,
\end{split}
\end{align}
\textcolor{black}{from the very definition of $\delta $ and the previous control on $\sup_{r\in[s,t]}|\varphi_r| $}.
Now set
$$
\Sigma_0:=\{\bxi_0\}=\{\x\},\ \ \Sigma_M:=\{\bxi_M\}=\{\y\},
$$
and for $j=1,\cdots, M-1$, 
$$
\Sigma_j:=\Big\{\z\in\mR^{2d}:  |\T^{-1}_\delta(\z-\bxi_{j})|\leq 1\Big\}.
$$ 
By \eqref{EG5} and \eqref{EQ_EQUIV_FLOW}, for any $j=0,\cdots,M-1$, 
we have that for $\z_j\in \Sigma_j$ and $\z_{j+1}\in \Sigma_{j+1}$,
\begin{align*}
|\T^{-1}_\delta(\btheta_{t_{j+1},t_j}(\z_j)-\z_{j+1})|
&\leq  |\T^{-1}_\delta(\btheta_{t_{j+1},t_j}(\bxi_j)-\bxi_{j+1})|
+|\T^{-1}_\delta(\z_{j+1}-\bxi_{j+1})|\\
&\quad+|\T^{-1}_\delta(\btheta_{t_{j+1},t_j}(\z_j)-\btheta_{t_{j+1},t_j}(\bxi_j))|
\\&\leq  |\T^{-1}_\delta(\btheta_{t_{j+1},t_j}(\bxi_j)-\bxi_{j+1})|+|\T^{-1}_\delta(\z_{j+1}-\bxi_{j+1})|\\
&\quad+\kappa_3(|\T^{-1}_\delta(\z_j-\bxi_j)|+1)\\
&\leq \kappa_5(2\kappa_6+1)+2\kappa_3+1=\kappa_9.
\end{align*}
This precisely means that the previous diagonal lower bound holds for $p(t_j, \z_j; t_{j+1},\z_{j+1} )$. 
Thus, by the Chapman-Kolmogorov equation and \eqref{loc_lower}, we have
\begin{align*}
p(s,\x;t,\y)&=\int_{\R^{2d}}\cdots\int_{\R^{2d}}p(t_0,\x; u,\z_1)\cdots p(t_{M-1}, \z_{M-1}; t_{M},\y)\dif \z_1\cdots \dif \z_{M-1}
\\&\geq \int_{\Sigma_1}\cdots\int_{\Sigma_{M-1}}p(t_0,\z_0; u,\z_1)\cdots p(t_{M-1}, \z_{M-1}; t_{M},\z_{M})\dif \z_1\cdots \dif \z_{M-1}
\\&\geq (C_0\delta^{-2d})^{M}\int_{\Sigma_1}\cdots\int_{\Sigma_{M-1}}\dif \z_1\cdots \dif \z_{M-1}
=(C_0\delta^{-2d})^{M} (C_4\delta^{2d})^{M-1},
\end{align*}
where $C_0$ is given in \eqref{loc_lower}, and the last equality is due to $|\Sigma_j|=C_4\delta^{2d}$ for some $C_4$ only depending on $d$.
Recalling $\delta=(t-s)/M$ and $M$ given in \eqref{DP1}, we finally have
\begin{align*}
p(s,\x;t,\y)&\geq C^{M}_0 (C_4)^{M-1}\delta^{-2d}=(t-s)^{-2d} M^{2d}\exp\{M\log(C_0C_4)\}/C_4
\\&\geq   C_5(t-s)^{-2d}\exp\{-\lambda_0|\T^{-1}_{t-s}(\btheta_{t,s}(\x)-\y)|^2\}.
\end{align*}
The lower bound is thus obtained for $t-s\leq T_0$ and all $\x,\y\in\mR^{2d}$. For general $0\leq s<t\leq T$, it again follows from the Chapman-Kolmogorov equation.
\subsection{The parametrix series expansion of $p(s,\x;t,\y)$}
\label{SEC_FULL_PARAM}
We introduce, for $\delta>0$ the SDE \eqref{SDE} with diffusion coefficient $\sigma(s,\x)=\delta \mI_{d\times d}$ and denote with 
$\bar p_{\delta}$ the corresponding density. 
By scaling and 
the two-sided density bounds of Sections \ref{UPPER} and \ref{LOWER} it holds that, for any $\l>0$, there exists $\delta=\delta(\l)$ large enough, and $C_{\delta}\geq 1$, $\l'>\l$ such that, 
for all $\x,\y\in \R^{2d}$ and $0\le s<t\le T$,
\begin{equation}\label{HK2}
C_{\delta}^{-1}\hat{p}_{\l}(s,\x;t,\y)\le \bar p_{\delta}(s,\x;t,\y)\le C_{\delta}\hat{p}_{\l'}(s,\x;t,\y).
\end{equation}
By Lemmas \ref{Le33} and \ref{Le34}, we may choose $\l$, and then $\delta=\delta(\l)$ such that, for any $j=(j_1,j_2)\in \N^2_0$ 
with $j_1+3j_2\le 3$, $\textcolor{black}{k\in \{1,2\}}$, 
$\beta\in [0,2]$ and for all $0\le s<t\le T$, $\x,\x',\y\in \R^{2d}$ with $|\x-\x'|_{\bf d}\le (t-s)^{\frac 12}$, 
\begin{equation}
|(\x-\btheta_{s,t}(\y))_k|^{\beta}
|\nabla_{x_1}^{j_1} \nabla_{x_2}^{j_2}\wt p_1(s,\x;t,\y)|\lesssim
(t-s)^{\beta (k-{\frac 12})-\frac{k_1+3j_2}{2}}\bar p_{\delta}(s,\x;t,\y), \label{Proxy_B1}
 \end{equation}
and also, for any $\alpha \in [0,1]$,
\begin{align} \label{Proxy_B1bis}
\begin{split}
&|(\x-\btheta_{s,t}(\y))_k|^{\beta} |\nabla_{x_1}^{j_1} \nabla_{x_2}^{j_2}({\wt p}_1(s,\x;t,\y)-{\wt p}_1(s,\x'; t,\y))|\\
&\qquad\lesssim|\x-\x'|_{\bf d}^{\alpha}(t-s)^{\beta (k-{\frac 12})-\frac{\alpha +j_1+3j_2}{2}}\bar p_{\delta}(s,\x;t,\y).
 \end{split}
\end{align}
We shall fix from now on $\delta$ such that \eqref{Proxy_B1} and \eqref{Proxy_B1bis} hold and for simplicity we write $\bar p=\bar p_{\delta}$.
Importantly, $\bar{p}$ enjoys the Chapman-Kolmogorov property, namely, for all $\x,\y\in \R^{2d}$ and $0\le t< r<s\le T$,
\begin{equation}\label{CK}
\int_{\R^{2d}}\bar p(s,\x;r,\z)\bar p(r,\z; t,\y)\dif\z=\bar p(s,\x;t,\y).
\end{equation}

\bp\label{Prop_expansion}
Under {\bf (H$^\gamma_\sigma$)}, {\bf (H$^{\gamma}_{\F}$)} and \eqref{Smo1}, the density $p(s,\x;t,\y)$ admits the following parametrix expansion 
\begin{equation}\label{parametrix_expansion}
p(s,\x;t,\y)=\wt{p}_1(s,\x;t,\y)+\wt{p}_1\otimes \cH(s,\x;t,\y),
\end{equation}
where $\cH:=\sum_{k\geq 1}H^{\otimes k}$ enjoys the following estimates: for all $0\le s<t\leq T$ and $\x,\y\in\R^{2d}$,
\begin{align}
|\cH(s,\x;t,\y)|&\lesssim (t-s)^{\frac{\gamma}{2}-1}\bar p(s,\x;t,\y).\label{Series_Ctr}
\end{align}
If in addition $F_1$ also satisfies \eqref{Re2}, then for all $\eps\in(0,\gamma)$ and $\x,\x',\y\in\R^{2d}$,
\begin{align}
|\cH(s,\x;t,\y)-\cH(s,\x'; t,\y)|&
\lesssim |\x-\x'|_{\bf d}^{\gamma-\eps}(t-s)^{\textcolor{black}{\frac \eps 2}-1}\left(\bar p(s,\x;t,\y)+\bar p(s,\x'; t,\y)\right).\label{Series_Ctr_H}
\end{align}
\ep
\proof 
Let us first prove \eqref{Series_Ctr} and \eqref{Series_Ctr_H}.
From the definition \eqref{HH0} of $H$, the proof of Proposition \ref{Conv_Kernels} and \eqref{Proxy_B1}, we derive that 
for all $\x,\y\in\R^{2d}$ and $0\le s<t\le T$, 
\begin{equation}\label{kernel_first}
|H(s,\x;t,\y)|\le C (t-s)^{\frac{\gamma}{2}-1}\bar p (s,\x;t,\y).
\end{equation}
The point here is that $\bar p$ is a true density and enjoys the reproduction property \eqref{CK}. 
This allows to manage the iteration procedure without deteriorating
the constants at every step. Indeed by direct induction, for $k\geq 1$, 
\begin{equation}\label{Kernel_iterated}
|H^{\otimes k}(s,\x;t,\y)|\le C^k \frac{\Gamma^k\left(\gamma/2\right)}{\Gamma\left(k\gamma/2\right)} (t-s)^{\frac{k\gamma}{2}-1}\bar p (s,\x;t,\y),
\end{equation}
where $\Gamma$ is the Euler-Gamma function. Estimate \eqref{Series_Ctr} easily follows. 

Next let us consider \eqref{Series_Ctr_H}. When $|\x-\x'|_{\bf d}> (t-s)^{\frac 12}$ (off-diagonal regime) the estimate directly follows from 
\eqref{Series_Ctr}. When $|\x-\x'|_{\mathbf d}\leq (t-s)^{\frac 12}$ we restart from the first term of the expansion. By \eqref{HH0}, we have
\begin{align*}
&|H(s,\x;t,\y)-H(s,\x'; t,\y)|\\
&\quad =\left|(\cL_{s,\x}-\wt\cL^{(t,\y)}_{s,\x})\wt{p}_1(s,\x;t,\y)-(\cL_{s,\x}-\wt\cL^{(t,\y)}_{s,\x})\wt{p}_1(s,\x'; t,\y) \right|\\
&\quad \leq \Big[ |a(s,\x)-a(s,\btheta_{s,t}(\y))| |\nabla^2_{x_1}\left(\wt{p}_1(s,\x;t,\y)-\wt{p}_1(s,\x'; t,\y)\right)|\\
&\quad \quad +|a(s,\x)-a(s,\x')| |\nabla^2_{x_1}\wt{p}_1(s,\x'; t,\y)|  \Big]\\
&\quad \quad+ \Big[|F_1(s,\x)-F_1(s,\btheta_{s,t}(\y))| |\nabla_{x_1}\left(\wt{p}_1(s,\x;t,\y)-\wt{p}_1(s,\x'; t,\y)\right)| \\
&\quad \quad + |F_1(s,\x)-F_1(s,\x')| |\nabla_{x_1}\wt{p}_1(s,\x'; t,\y)| \Big]\\
&\quad \quad+ \Big[|\cT_{F_2(s)}(\x,\btheta_{s,t}(\y))|\,|\nabla_{x_2}\left(\wt{p}_1(s,\x;t,\y)-\wt{p}_1(s,\x'; t,\y)\right)| \\
&\quad\quad  + |\cT_{F_2(s)}(\x,\btheta_{s,t}(\y))-\cT_{F_2(s)}(\x',\btheta_{s,t}(\y))|\, |\nabla_{x_2}\wt{p}_1(s,\x'; t,\y)|\Big]\\
&\quad =: I_1+I_2+I_3.
\end{align*}
By \eqref{Proxy_B1}, \eqref{Proxy_B1bis} and using that $|\x-\x'|_{\bf d}\leq (t-s)^{\frac 12}$ we have
\begin{align*}
I_1&\lesssim |\x-\x'|_{\bf d}^{\gamma}\left(\frac{|\x-\btheta_{s,t}(\y)|_{\bf d}^{\gamma}}{(t-s)^{1+\frac{\gamma}{2}}}\hat p_{\l_2}(s,\x;t,\y)+\frac{1}{t-s}\hat p_{\l_2}(s,\x'; t,\y)\right)\\
&\lesssim |\x-\x'|_{\bf d}^{\gamma-\eps}(t-s)^{\frac\eps2-1}\left(\bar p(s,\x;t,\y)+\bar p(s,\x'; t,\y)\right).
\end{align*}
Similarly, by \eqref{F11} we also have 
\begin{align*}
I_2&\lesssim |\x-\x'|_{\bf d}^{\gamma}(t-s)^{-\frac 12}\left(\bar p(s,\x;t,\y)+\bar p(s,\x'; t,\y)\right),
\end{align*}
and by \eqref{Taylor}, \eqref{Proxy_B1} and \eqref{Proxy_B1bis},
\begin{align*}
I_3&\lesssim 
|\x-\btheta_{s,t}(\y)|_{\bf d}^{1+\gamma} (t-s)^{-\frac{3+\gamma}{2}}|\x-\x'|_{\bf d}^{\gamma}\hat p_{\l_3}(s,\x;t,\y)\\
&\quad+|\x-\x'|_{\bf d}^{1+\gamma} (t-s)^{-\frac 32}\hat p_{\l_3}(s,\x'; t,\y)\\
&\lesssim |\x-\x'|_{\bf d}^{\gamma-\eps}(t-s)^{\textcolor{black}{\frac \eps 2}-1}\left(\bar p(s,\x;t,\y)+\bar p(s,\x'; t,\y)\right).
\end{align*}
Finally, we have
\begin{align*}
&|\cH(s,\x;t,\y)- \cH(s,\x'; t,\y)|\\
&\qquad \leq |H(s,\x;t,\y)- H(s,\x'; t,\y)|+\int^t_s\int_{\R^{2d}}|H(s,\x;r,\z)- H(s,\x';r,\z)||\cH(r,\z; t,\y)|\dif\z \dif r\\
&\qquad \lesssim  |\x-\x'|_{\bf d}^{\gamma-\eps}(t-s)^{\textcolor{black}{\frac \eps2}-1}\left(\bar p(s,\x;t,\y)+\bar p(s,\x'; t,\y)\right)\\
&\qquad \quad + \int^t_s\frac{|\x-\x'|_{\bf d}^{\gamma-\frac \eps 2}}{(r-s)^{1-\textcolor{black}{\frac \eps 2}}(t-r)^{1-\textcolor{black}{\frac \eps 2}}}
\int_{\R^{2d}}\left(\bar p(s,\x;r,\z)+\bar p(s,\x';r,\z)\right)\bar p(r,\z; t,\y)\dif\z \dif r \\
&\qquad \lesssim  |\x-\x'|_{\bf d}^{\gamma-\eps}(t-s)^{\textcolor{black}{\frac \eps 2}-1}\left(\bar p(s,\x;t,\y)+\bar p(s,\x'; t,\y)\right).
\end{align*}
The parametrix expansion \eqref{parametrix_expansion} follows by letting $N\longrightarrow \infty$ in \eqref{FR1} and proving that the remainder $R^N:=p\otimes H^{\otimes N}$ converges to zero uniformly in every variable. 
Let $M_N:=C^{N}{\Gamma^N\left(\frac{\gamma}{2}\right)}\Gamma\left(\frac{N\gamma}{2}\right)^{-1}$ with $C$ as in \eqref{Kernel_iterated}. 
By \eqref{Kernel_iterated} we have 
\begin{align*}
|(p\otimes H^{\otimes N})(s,\x;t,\y)|&\le C M_N \int^t_s (t-r)^{\frac{N\gamma}{2}-1}\int_{\R^{2d}}\bar p(s,\x;r,\z)\bar p(r,\z; t,\y)\dif\z \dif r\\
&\le C\frac{M_N}{N}(t-s)^{\frac{N\gamma}{2}}\bar p(s,\x;t,\y)\longrightarrow 0.
\end{align*}
The proof is completed.
\endproof

\subsection{Sensitivities of the frozen densities with respect to spatial parameters}
We provide here a lemma which quantifies the sensitivity of the frozen densities w.r.t. different (spatial) freezing parameters. 
\bl[Sensitivity of the semigroup w.r.t. the freezing parameters]\label{Lem_sensitivity}
Under {\bf (H$^\gamma_\sigma$)}, {\bf (H$^{\gamma}_{\F}$)} and \eqref{Smo1}, for any $T>0$
and $j=(j_1,j_2)\in\mN_0^2$, there exist constants
$\lambda_j,C_j\geq 1$ depending only on $\Theta_T$ such that
for all $0\leq s<t\leq T$  and $\x,\y\in\mR^{2d}$,
\begin{align}\label{THE_CTR_SENSI_DIFF_PARAM}
\begin{split}
&\big|\nabla^{j_1}_{x_1}\nabla_{x_2}^{j_2}\(\wt{p}_1-\wt{p}^{(\tau,\bxi)}\)(s,\x;t,\y)\big|_{(\t,\bxi)=(s,\x)}
\leq C_j(t-s)^{\frac{\gamma-(j_1+3j_2)}{2}}\hat p_{\lambda_j}(s,\x;t,\y);
\end{split}
\end{align}
Moreover, for any $\x'\in\R^{2d}$ 
\begin{equation}\label{THE_CTR_SENSI_DIFF_PARAMbis}
\big|\nabla^{j_1}_{x_1}\nabla_{x_2}^{j_2}\(\wt{p}^{(\tau,\bxi)}-\wt{p}^{(\tau,\bxi')}\)(s,\x;t,\y)\big|_{(\t,\bxi,\bxi')=(s,\x,\x')}
\leq C_j|\x-\x'|^{\gamma}(t-s)^{\frac{-(j_1+3j_2)}{2}}\hat p_{\lambda_j}(s,\x;t,\y).
\end{equation}
\el


\begin{proof}
We only prove \eqref{THE_CTR_SENSI_DIFF_PARAM} for $j_1=j_2=0$. For general $j_1, j_2\in\mN_0$, the statement follows by the chain rule and tedious but similar calculations. The control \eqref{THE_CTR_SENSI_DIFF_PARAMbis} is derived analogously. For  notational simplicity, we introduce
$$
\cC_1=\K^{(s,\x)}_{t,s},\ \cC_2=\K^{(t,\y)}_{t,s},\ \w_1=\btheta_{t,s}(\x)-\y,\ \w_2=\bvtheta^{(t,\y)}_{t,s}(\x)-\y,
$$ 
and
$$
\cA:=|\T^{-1}_{t-s}(\btheta_{t,s}(\x)-\y)|+1.
$$
Using the above notations and by definition \eqref{CORRESP}, we have
\begin{align*}
&\big|\(\wt{p}_1-\wt{p}^{(\tau,\bxi)}\)(s,\x;t,\y)\big|_{(\t,\bxi)=(s,\x)}\\
&\qquad=(2\pi)^{-d}\Big|\det(\cC_2)^{-\frac 12}
\exp\Big( -\tfrac12\big|\cC_2^{-\frac12}\w_2\big|^2\Big)
-\det(\cC_1)^{-\frac 12}
\exp\Big( -\tfrac12\big|\cC_1^{-\frac12}\w_1\big|^2\Big)\Big|.
\end{align*}
To show \eqref{THE_CTR_SENSI_DIFF_PARAM}, it suffices to prove that for some $\lambda>0$,
\begin{align}
|(\det \cC_1)^{-\frac 12}-(\det \cC_2)^{-\frac 12}|&\lesssim (t-s)^{-2d+\frac{\gamma}{2}}\cA^{\textcolor{black}{\gamma}},\label{e12}\\
\Big|\exp\(-\tfrac{1}{2}|\cC_1^{-\frac12}\w_1|^2\)-\exp\(-\tfrac{1}{2}|\cC_2^{-\frac12}\w_2|\)\Big|
&\lesssim (t-s)^{\frac{\gamma}{2}}\exp\left(-\lambda\cA^2\right).\label{e13}
\end{align}
Note that by \eqref{TT2},
\begin{align}\label{TT6}
\cC_i=\textcolor{black}{\T_{t-s}\hat\cC_i\T_{t-s}},
\end{align}
where
$$
\hat\cC_1:=\int^1_0\hat\gR^{(s,\x)}_{1,r}\hat\bTheta^{(s,\x)}_r(\hat\gR^{(s,\x)}_{1,r}\hat\bTheta^{(s,\x)}_r)^*\dif r,\
\hat\cC_2:=\int^1_0\hat\gR^{(t,\y)}_{1,r}\hat\bTheta^{(t,\y)}_r(\hat\gR^{(t,\y)}_{1,r}\hat\bTheta^{(t,\y)}_r)^*\dif r,
$$
and for $\eta(r):=s+(t-s)r$,
$$
\hat\bTheta^{(\tau,\bxi)}_r:=B\sigma(\eta(r),\btheta_{\eta(r),\tau}(\bxi)),
$$
and
$$
\hat\gR^{(\tau,\bxi)}_{1,r}:=\begin{pmatrix}\mI_{d\times d} & 0_{d\times d}\\
\int^1_r \nabla_{x_1}F_2(\eta(u),\btheta_{\eta(u),\t}(\bxi))\dif u & \mI_{d\times d} \end{pmatrix}.
$$
Since $\sigma,\nabla_{x_1}F_2\in\sC_{\bf d}^\gamma$, by \eqref{EZ1} we have for any $r\in[0,1]$,
$$
|\hat\bTheta^{(s,\x)}_r-\hat \bTheta^{(t,\y)}_{r}|\lesssim|\btheta_{\eta(r),s}(\x)-\btheta_{\eta(r),t}(\y)|_{\bf d}^\gamma\lesssim(t-s)^{\frac\gamma2}\cA^\gamma,
$$
and
$$
|\hat\gR^{(s,\x)}_{1,r}-\hat \gR^{(t,\y)}_{1,r}|\lesssim(t-s)^{\frac\gamma2}\cA^\gamma.
$$
Hence, by \eqref{TT1},
\begin{align}\label{TT9}
|\hat\cC_1-\hat\cC_2|\lesssim (t-s)^{\frac\gamma2}\cA^\gamma,\ \ |\hat\cC_i|\lesssim 1,\ \ \det\hat \cC_i\asymp 1,
\end{align}
and
\begin{align*}
|(\det \cC_1)^{-\frac 12}-(\det \cC_1)^{-\frac 12}|&\lesssim (t-s)^{-2d}|(\det \hat\cC_1)^{-\frac 12}-(\det \hat\cC_1)^{-\frac 12}|\\
&\lesssim (t-s)^{-2d}|\det \hat\cC_1-\det \hat\cC_1|\\
&\lesssim (t-s)^{-2d}|\hat\cC_1-\hat\cC_2|\lesssim (t-s)^{-2d+\frac{\gamma}{2}}\cA^{\gamma}.
\end{align*}
Thus we obtain \eqref{e12}. For proving \eqref{e13}, without loss of generality, we may assume
$$
|\cC_1^{-\frac12}\w_1|\leq |\cC_2^{-\frac12}\w_2|.
$$
Then by $1-\e^{-x}\leq x$, we have
\begin{align*}
&\Big|\exp\(-\tfrac{1}{2}|\cC_1^{-\frac12}\w_1|^2\)-\exp\(-\tfrac{1}{2}|\cC_2^{-\frac12}\w_2|^{\textcolor{black}{2}}\)\Big|\leq 
\tfrac12(|\cC_2^{-\frac12}\w_2|^2-|\cC_1^{-\frac12}\w_1|^2)\exp\(-\tfrac{1}{2}|\cC_1^{-\frac12}\w_1|^2\).
\end{align*}
Since
\begin{align*}
&|\cC_2^{-\frac12}\w_2|^2-|\cC_1^{-\frac12}\w_1|^2\stackrel{\eqref{TT6}}{=}\<\hat\cC_2^{-1}\T^{-1}_{t-s}\w_2,\T^{-1}_{t-s}\w_2\>-\<\hat\cC_1^{-1}\T^{-1}_{t-s}\w_1,\T^{-1}_{t-s}\w_1\>\\
&\qquad\leq|\hat\cC_1^{-1}||\T^{-1}_{t-s}(\w_2-\w_1)|\(|\T^{-1}_{t-s}\w_2|+|\T^{-1}_{t-s}\w_1|\)+|\hat\cC_2^{-1}-\hat\cC_1^{-1}||\T^{-1}_{t-s}\w_1|^2,
\end{align*}
by \eqref{TT8} and \eqref{TT9}, we get
$$
|\cC_2^{-\frac12}\w_2|^2-|\cC_1^{-\frac12}\w_1|^2\lesssim(t-s)^{\frac\gamma2}\cA^{2+\gamma}.
$$
Therefore, we have \eqref{e13}. The proof is complete.
\end{proof}

\subsection{Gradient bounds in the non-degenerate direction $\x_1$}\label{SEC_GRAD_X_1}
In this Section we provide pointwise controls for the derivatives of the densities in $\x_1$ in the current smooth setting. 
Importantly the controls do not depend on the smoothing procedure, and the constants therein only depend on $\Theta_T$.

We are ready to prove the following

\bp\label{gradient_x1}
Under {\bf (H$^\gamma_\sigma$)}, {\bf (H$^{\gamma}_{\F}$)} and \eqref{Smo1}, for any $T>0$ and $j=0,1$,
there exists $\lambda_j, C_j\geq 1$ depending on $\Theta_T,\eta$ such that for any $0\le s<t\le T$ and $\x,\y\in\R^{2d}$, 
\begin{equation}\label{EST_GRAD_X1}
|\nabla^{j}_{x_1} p(s,\x;t,\y)|\leq C_j(t-s)^{-\frac{j}{2}}\hat p_{\lambda_j} (s,\x;t,\y),
\end{equation}
and for $\eta_0,\eta_1\in(0,1)$, and any $\x,\x',\y\in \R^{2d}$,
\begin{align}
|\nabla^j_{x_1}(p(s,\x;t,\y)-p(s,\x'; t,\y))|&\le C_j|\x-\x'|_{\bf d}^{\eta_j}(t-s)^{-\frac{j+\eta_j}{2}}\(\hat p_{\lambda_j}(s,\x;t,\y)+ \hat p_{\lambda_j}(s,\x';t,\y)\).
\label{HOLD_MOD_D_X1_2}
\end{align}
If in addition $F_1$ also satisfies \eqref{Re2}, then \eqref{EST_GRAD_X1} and \eqref{HOLD_MOD_D_X1_2} also \textcolor{black}{hold} for $j=2$ and $\eta_2\in(0,\gamma)$.
\ep
\begin{proof}
We only prove the case $j=2$. In this case some time singularities \textcolor{black}{appear} in the integrals. 
A way to overcome such a problem is to exploit cancellation properties of the derivatives of the Gaussian kernels. 
In the following estimates, for simplicity, we use the same $\lambda$ to denote possible different constants in different places.

(i) We first look at \eqref{EST_GRAD_X1}.  Starting from \eqref{parametrix_expansion}, for $u=\frac{t+s}{2}$ and $(\tau,\bxi)\in[s,t]\times\mR^{2d}$, we write
\begin{align*}
(\nabla^2_{x_1}\wt p_1\otimes \cH)(s,\x;t,\y)&=\int^t_u\int_{\R^{2d}}\nabla^2_{x_1}\wt{p}_1(s,\x;r,\z)\cH(r,\z; t,\y)\dif\z \dif r\\
&+\int^u_s\int_{\R^{2d}}\nabla^2_{x_1}\(\wt{p}_1-\wt{p}^{(\tau,\bxi)}\)(s,\x;r,\z)\cH(r,\z; t,\y)\dif\z \dif r\notag\\
&  + \int^u_s\int_{\R^{2d}}\nabla^2_{x_1}\wt{p}^{(\tau,\bxi)}(s,\x;r,\z)\left(\cH(r,\z; t,\y)-\cH(r,\btheta_{r,\tau}(\bxi); t,\y)\right)\dif\z \dif r\\
&=:I_1+I^{(\tau,\bxi)}_2+I^{(\tau,\bxi)}_3.
\end{align*}
For $I_1$, by \eqref{Proxy_B1} and \eqref{Series_Ctr} we have
\begin{align*}
I_1&\lesssim \int^t_u(r-s)^{-1}\int_{\R^{2d}}\bar{p}(s,\x;r,\z)(t-r)^{\frac{\gamma}{2}-1}\bar p(r,\z; t,\y)\dif\z \dif r\\
&\leq2 (t-s)^{-1}\left(\int^t_u(t-r)^{\frac{\gamma}{2}-1}\dif r\right)\bar p(s,\x; t,\y)\lesssim(t-s)^{\frac{\gamma}{2}-1}\bar p(s,\x; t,\y).
\end{align*}
For $I^{(\tau,\bxi)}_2$, taking ${(\tau,\bxi)}=(s,\x)$, by \eqref{THE_CTR_SENSI_DIFF_PARAM}, \eqref{Series_Ctr} and \eqref{Convolution}, we have
\begin{align*}
I^{(s,\x)}_2&\lesssim \int^u_s\int_{\R^{2d}}(r-s)^{\frac{\gamma}{2}-1}\hat p_{\lambda}(s,\x;r,\z)(t-r)^{\frac{\gamma}{2}-1}\bar p(r,\z; t,\y)\dif\z \dif r\\
&\lesssim \left(\int^u_s(r-s)^{\frac{\gamma}{2}-1}(t-r)^{\frac{\gamma}{2}-1} \dif r\right)\hat p_{\lambda}(s,\x; t,\y)\\
&\lesssim (t-s)^{\gamma-1}\hat p_{\lambda}(s,\x; t,\y).
\end{align*}
For $I^{(\tau,\bxi)}_3$, \textcolor{black}{with the former choice} ${(\tau,\bxi)}=(s,\x)$, by \eqref{THE_CTR_SENSI_DIFF_PARAM}, \eqref{Series_Ctr_H}, \eqref{Proxy_B1} and \eqref{Convolution}, we have 
\begin{align*}
I^{(s,\x)}_3&\lesssim\int^u_s\int_{\R^{2d}}(r-s)^{-1}\hat p_\lambda (s,\x;r,\z)|\z-\btheta_{r,s}(\x)|_{\bf d}^{\frac{\gamma}{2}}(t-r)^{\frac{\gamma}{2}-1}\\
&\qquad\times\Big[\bar p(r,\z; t,\y)+\bar p(r,\btheta_{r,s}(\x); t,\y)\Big]\dif\z \dif r\\
&\lesssim\int^u_s\int_{\R^{2d}}(r-s)^{\frac{\gamma}{4}-1}\hat p_{\lambda} (s,\x;r,\z)(t-r)^{\frac{\gamma}{2}-1}\\
&\qquad\times\Big[\hat p_\lambda(r,\z; t,\y)+\hat p_\lambda(r,\btheta_{r,s}(\x); t,\y)\Big]\dif\z \dif r\\
&\lesssim (t-s)^{\frac\gamma2-1}\left[\hat p_{\l}(s,\x;t,\y)+\int^u_s(r-s)^{\frac{\gamma}{4}-1}\hat p_\lambda(r,\btheta_{r,s}(\x); t,\y) \dif r\right]\\
&\lesssim(t-s)^{\frac\gamma2-1}\hat p_{\l}(s,\x;t,\y),
\end{align*}
where in the last step we have used that
\begin{align}\label{VX1}
\hat p_{\lambda} (r,\btheta_{r,s}(\x); t,\y)=g_{\lambda}\(t-r,\btheta_{t,r}\circ\btheta_{r,s}(\x)-\y\)\stackrel{\eqref{GG1}}{=}g_{\lambda}\(t-r,\btheta_{t,s}(\x)-\y\).
\end{align}
Finally, note that
\begin{align*}
|\nabla^2_{x_1}p(s,\x;t,\y)|&\lesssim (t-s)^{-1}\bar p(s,\x;t,\y)+|(\nabla^2_{x_1}\wt p_1\otimes \cH)(s,\x;t,\y)|,
\end{align*}
which combining the above calculations, yields \eqref{EST_GRAD_X1} for $j=2$.

(ii) Next we look at \eqref{HOLD_MOD_D_X1_2}. \textcolor{black}{We can assume w.l.o.g. that $|\x-\x'|\le (t-s)^{\frac 12} $ since otherwise the control readily follows from the previous one}.
Starting from \eqref{parametrix_expansion}, for $u=\frac{t+s}{2}$  we write
\begin{align*}
&\nabla^2_{x_1}p(s,\x;t,\y)-\nabla^2_{x_1}p(s,\x';t,\y)=\nabla^2_{x_1}\wt p_1(s,\x;t,\y)-\nabla^2_{x_1}\wt p_1(s,\x';t,\y)\\
&\qquad+\int^t_s\int_{\mR^{2d}}(\nabla^2_{x_1}(\wt p_1(s,\x;r,\z)-\nabla^2_{x_1}\wt p_1(s,\x';r,\z))\cH(r,z;t,\y)\dif z\dif r\\
&\qquad=:I_1+\int^t_u\cJ(r,\x,\x'; t,\y)\dif r+\int^u_s\cJ(r,\x,\x'; t,\y)\dif r=:I_1+I_2+I_3.
\end{align*}
For $I_1$, by \eqref{Proxy_B1bis} we have
$$
I_1\lesssim |\x-\x'|_{\bf d}^{\gamma}(t-s)^{-1-\frac{\gamma}{2}}\(\hat p_{\lambda}(s,\x;t,\y)+ \hat p_{\lambda}(s,\x';t,\y)\).
$$
For $I_2$, by \eqref{Proxy_B1bis} and \eqref{Series_Ctr}, we have
\begin{align*}
I_2&\lesssim |\x-\x'|_{\bf d}^{\gamma}\int^t_u(r-s)^{-1-\frac{\gamma}{2}}(t-r)^{\frac\gamma2-1}\int_{\mR^{2d}}\(\bar p(s,\x;r,\z)+ \bar p(s,\x';r,\z)\)\bar p(r,\z;t,\y)\dif z\dif r\\
&\lesssim |\x-\x'|_{\bf d}^{\gamma}(t-s)^{-1}\(\bar p(s,\x;t,\y)+ \bar p(s,\x';t,\y)\).
\end{align*}
Next we treat the hard term $I_3$.
Fix $r\in[s,u]$. We \textcolor{black}{ handle $\cJ$ according to the current \textit{diagonal/off-diagonal} regime of $|\x-\x'|_{\mathbf d} $ w.r.t. to the current integration time} and consider two cases: 
$$
\mbox{Case \textcolor{black}{(I)}: } |\x-\x'|_{\mathbf d}>(r-s)^{\frac 12};\ \mbox{Case \textcolor{black}{(II)}: } |\x-\x'|_{\mathbf d}\leq (r-s)^{\frac 12}.
$$
In Case (I), for any $(\tau,\bxi), (\tau,\bxi')\in[s,t]\times\mR^{2d}$,  we write
\begin{align*}
\cJ&=\int_{\R^{2d}}\nabla^2_{x_1}\(\wt{p}_1-\wt{p}^{(\tau,\bxi)}\)(s,\x;r,\z)\cH(r,\z; t,\y)\dif\z\\
& \quad +\int_{\R^{2d}}\nabla^2_{x_1}\wt{p}^{(\tau,\bxi)}(s,\x;r,\z)\(\cH(r,\z; t,\y)-\cH(r,\btheta_{r,\tau}(\bxi); t,\y)\)\dif\z\\
&\hspace*{5pt}-\int_{\R^{2d}}\nabla^2_{x_1}\(\wt{p}_1-\wt{p}^{(\t,\bxi')}\)(s,\x';r,\z)\cH(r,\z; t,\y)\dif\z\\
& \quad + \int_{\R^{2d}}\nabla^2_{x_1}\wt{p}^{(\t,\bxi')}(s,\x';r,\z)\left(\cH(r,\z; t,\y)-\cH(r,\btheta_{t,\tau}(\bxi'); t,\y)\right)\dif\z\\
&=:\sum_{i=1,2}\Big(\sI^{(\tau,\bxi)}_i(r,\x; t,\y)+\sN_i^{(\tau,\bxi')}(r,\x'; t,\y)\Big).
\end{align*}
By \eqref{THE_CTR_SENSI_DIFF_PARAM} and \eqref{Series_Ctr}, we have
\begin{align*}
|\sI^{(s,\x)}_1(r,\x; t,\y)|&\lesssim (r-s)^{\frac{\gamma}{2}-1}(t-r)^{\frac{\gamma}{2}-1}\int_{\mR^{2d}}\bar p(s,\x;r,\z)\bar p(r,\z;t,\y)\dif z\\
&= (r-s)^{\frac{\gamma}{2}-1}(t-r)^{\frac{\gamma}{2}-1}\bar p(s,\x;t,\y).
\end{align*}
By \eqref{THE_CTR_SENSI_DIFF_PARAM} and \eqref{Series_Ctr_H}, we have
\begin{align*}
|\sI^{(s,\x)}_2(r,\x; t,\y)|&\lesssim (r-s)^{-1}(t-r)^{\eps-1}
\int_{\mR^{2d}}\hat p_\lambda (s,\x;r,\z)|\z-\btheta_{r,s}(\x)|_{\bf d}^{\gamma-\eps}\\
&\quad\times\Big[\bar p(r,\z; t,\y)+\bar p(r,\btheta_{r,s}(\x); t,\y)\Big]\dif z\\
&\stackrel{\eqref{regularizing_proxy}}{\lesssim} (r-s)^{\frac{\gamma-\eps}{2}-1}(t-r)^{\eps-1}
\int_{\mR^{2d}}\hat p_\lambda (s,\x;r,\z)\\
&\quad\times\Big[\hat p_{\lambda} (r,\z; t,\y)+\hat p_{\lambda} (r,\btheta_{r,s}(\x); t,\y)\Big]\dif z\\
&\stackrel{\eqref{regularizing_proxy},\eqref{VX1}}{\lesssim} (r-s)^{\frac{\gamma-\eps}{2}-1}(t-r)^{\eps-1}\hat p_{\lambda} (s,\x; t,\y).
\end{align*}
Similarly, we have
\begin{align*}
|\sN^{(s,\x')}_1(r,\x'; t,\y)|&\lesssim(r-s)^{\frac{\gamma}{2}-1}(t-r)^{\frac{\gamma}{2}-1}\bar p(s,\x';t,\y),\\
|\sN^{(s,\x')}_2(r,\x'; t,\y)|&\lesssim(r-s)^{\frac{\gamma-\eps}{2}-1}(t-r)^{\eps-1}\hat p_{\lambda} (s,\x'; t,\y).
\end{align*}
Combining the above estimates and taking $\eps$ small enough, we get 
\begin{align*}
\int^u_s\1_{|\x-\x'|_{\mathbf d}>(r-s)^{1/2}}|\cJ(r,\x,\x'; t,\y)|\dif r\lesssim
|\x-\x'|_{\mathbf d}^\eta(t-s)^{-1}\(\hat p_\lambda(s,\x;t,\y)+\hat p_\lambda(s,\x';t,\y)\),
\end{align*}
\textcolor{black}{for all $\eta\in (0,\gamma) $}.

In Case (II), for any $(\tau,\bar\bxi)\in[s,t]\times\mR^{2d}$,  we write
\begin{align*}
\cJ&=\int_0^1\Bigg[ \int_{\R^{2d}}(\x-\x')\cdot\nabla_\x \nabla^2_{x_1}
\(\wt{p}_1-\wt{p}^{(\t,\bar \bxi)}\)(s,\x'+\rho (\x-\x');r,\z)\cH(r,\z; t,\y)\dif\z \\
& + \int_{\R^{2d}}(\x-\x')\cdot\nabla_\x\nabla^2_{x_1}\wt{p}^{(\t,\bar\bxi)}(s,\x'+\rho (\x-\x');r,\z)\left(\cH(r,\z; t,\y)-\cH(r,\btheta_{r,s}(\bar \bxi); t,\y)\right)\dif\z\Bigg]\dif \rho.
\end{align*}
In the above bracket, taking $\bar\bxi=\bar\bxi_\rho:=\x'+\rho (\x-\x')$ and by \eqref{THE_CTR_SENSI_DIFF_PARAM} and \eqref{Series_Ctr_H}, we have
\begin{align*}
|\cJ|&\lesssim\int_0^1\Bigg[\sum_{i=1,2} \int_{\R^{2d}}|(\x-\x')_i|(r-s)^{\frac{\gamma}{2}-1-\frac{2i-1}{2}}\hat p_{\lambda}(s,\bar\bxi_\rho;r,\z)
(t-r)^{\frac{\gamma}{2}-1}\bar p(r,\z; t,\y)\dif\z \\
& \qquad + \sum_{i=1,2} \int_{\R^{2d}}|(\x-\x')_i|(r-s)^{-1-\frac{2i-1}{2}}p_{\lambda}(s,\bar\bxi_\rho;r,\z)|z-\btheta_{r,s}(\bar\bxi_\rho)|_{\bf d}^{\gamma-\eps}\\
&\qquad\qquad\times(t-r)^{\textcolor{black}{\frac \eps 2}-1}\Big[\hat p_{\lambda} (r,\z; t,\y)+\hat p_{\lambda_1} (r,\btheta_{r,s}(\bar\bxi_\rho); t,\y)\Big]\dif\z\Bigg]\dif \rho\\
&\lesssim|\x-\x'|_{\bf d}^\eta(r-s)^{\frac{\gamma-\eta}{2}-1}(t-s)^{\textcolor{black}{\frac \eps 2}-1}\int^1_0\hat p_{\lambda}(s,\bar\bxi_\rho;t,\y)\dif \rho\\
&\lesssim|\x-\x'|_{\bf d}^\eta(r-s)^{\frac{\gamma-\eta}{2}-1}(t-s)^{\textcolor{black}{\frac \eps 2}-1}\hat p_{\lambda}(s,\x';t,\y),
\end{align*}
where in the second step we have used $r\in[s,u]$, $|\x-\x'|_{\mathbf d}\leq(r-s)^{1/2}$, \eqref{Convolution} and \eqref{VX1}. \textcolor{black}{In} the last step we have used \eqref{VX2}.
Therefore,
\begin{align*}
\int^u_s\1_{|\x-\x'|_{\mathbf d}\leq(r-s)^{1/2}}|\cJ(r,\x,\x'; t,\y)|\dif r\lesssim
|\x-\x'|_{\mathbf d}^\eta(t-s)^{-1}\(\hat p_\lambda(s,\x;t,\y)+\hat p_\lambda(s,\x';t,\y)\).
\end{align*}
The whole proof is thus complete.
\end{proof} 
 
 \subsection{Gradient bound in degenerate direction $\x_2$}\label{GRAD_X_2}

The aim of this section is to show the following a priori gradient estimate.


\bp\label{gradient_x2}
Under {\bf (H$^\gamma_\sigma$)}, {\bf (H$^{\gamma}_{\F}$)}, \eqref{Re1}, \eqref{Re2}  and \eqref{Smo1}, for any $T>0$,
there exist constants $\lambda, C\geq 1$ depending on $\Theta_T$ and $\bar\kappa_0, \bar\kappa_1$ such that for any $0\le s<t\le T$ and $\x,\y\in\R^{2d}$, 
\begin{equation}\label{EST_GRAD_X2}
|\nabla_{x_2} p(s,\x;t,\y)|\leq C(t-s)^{-\frac{3}{2}}\hat p_{\lambda} (s,\x;t,\y).
\end{equation}
\ep
\proof 
For $0\leq s\leq t\leq T$, $\lambda>0$ and $f\in C^\infty_b(\mR^{2d})$, we define
$$
P_{s,t}f(\x):=\int_{\R^{2d}}p(s,\x;t,\z)f(\z)\dif\z,\ \ \hat{P}^\lambda_{s,t}f(\x):=\int_{\R^{2d}}\hat{p}_\lambda(s,\x;t,\z)f(\z)\dif\z,
$$
and for $(\tau,\bxi)\in[s,t]\times\mR^{2d}$,
\begin{equation}\label{Semigr_tilde}
\wt{P}^{(\tau, \bxi)}_{s,t}f(\x):=\int_{\R^{2d}}\wt{p}^{(\tau, \bxi)}(s,\x;t,\z)f(\z)\dif\z.
\end{equation}
Under \eqref{Smo1}, it is standard to derive that
$$
P_{s,t}f,\ \hat{P}^\lambda_{s,t}f,\ \wt{P}^{(\tau, \bxi)}_{s,t}f\in C^\infty_b(\mR^{2d}).
$$
Thus all the calculations below are rigorous. Now we set
$$
u:=(s+t)/2.
$$
Let $f\in C^\infty_b(\mR^{2d})$ be nonnegative. Noting that
$$
\wt{P}^{(\tau, \bxi)}_{s,t}f=\wt{P}^{(\tau, \bxi)}_{s,u}\wt{P}^{(\tau, \bxi)}_{u,t}f,
$$
by \eqref{D0}, it is easy to see that for any $(\tau,\bxi)\in[s,t]\times\mR^{2d}$,
\begin{align}\label{SG_LOC_IN_U}
P_{s,t}f(\x)=\wt{P}^{(\tau, \bxi)}_{s,u}P_{u,t}f(\x)+\int_s^u\wt{P}^{(\tau, \bxi)}_{s,r}
(\cL_{r}-\cL^{(\tau,\bxi)}_{r})P_{r,t}f(\x) \dif r.
\end{align}
Hence,
\begin{align}
\nabla_{x_2}P_{s,t}f(\x)=
\nabla_{x_2}\wt{P}^{(\tau,\bxi)}_{s,u}P_{u,t}f(\x)+\int_s^{u}\nabla_{x_2}\wt{P}^{(\tau,\bxi)}_{s,r}(\cL_{r}-\cL^{(\tau,\bxi)}_{r})P_{r,t}f(\x) \dif r.
\label{Duhamel_change_freez}
\end{align}
Below we shall take $(\tau,\bxi)=(s,\x)$. Observe that by \eqref{CORRESP} and \eqref{Res00},
\begin{align*}
\nabla_{x_2}\wt{p}^{(\tau, \bxi)}(s,\x;r,\z)\big|_{(\tau,\bxi)=(s,\x)}&=-\nabla_{z_2}\wt{p}^{(\tau, \bxi)}(s,\x;r,\z)\big|_{(\tau,\bxi)=(s,\x)}=-\nabla_{z_2}\wt{p}_0(s,\x;r,\z).
\end{align*}
In particular, 
$$
\nabla_{x_2}\wt{P}^{(\tau, \bxi)}_{s,r}f(\x)|_{(\tau,\bxi)=(s,\x)}=-\int_{\mR^{2d}}\nabla_{z_2}\wt{p}_0(s,\x;r,\z)f(\z)\dif\z=:Q_{s,r}f(\x),
$$
and for any bounded function $f(z_1)$ of the first variable,
\begin{align}\label{CAN0}
Q_{s,r}f(\x)=-\int_{\mR^{2d}}\nabla_{z_2}\wt{p}_0(s,\x;r,\z)f(z_1)\dif \z\equiv0.
\end{align}
\textcolor{black}{Equation \eqref{CAN0} is precisely what we call a partial cancellation property. 
}

Moreover, by \eqref{upper_proxy},
\begin{align}\label{CAN}
|\nabla_{z_2}\wt{p}_0(s,\x;r,\z)|\lesssim (r-s)^{-\frac32}\hat p_\lambda(s,\x;r,\z)\Rightarrow |Q_{s,r}f(\x)|\lesssim(r-s)^{-\frac32}\hat P^\lambda_{s,r}f(\x).
\end{align}
Here and below, $\lambda$ may \textcolor{black}{vary from line to line} \textcolor{black}{but} all the constants only depend on $\Theta_T$ and $\bar\kappa_0,\bar\kappa_1$.
For \textcolor{black}{notational convenience}, we write
\begin{align}\label{PRO2}
\cL_r-\cL^{(s,\x)}_r=\sD_r+\sB_r+\sK_r,
\end{align}
where
$$
\sD_r:=\tr[(a(s,\z)-a(r,\btheta_{r,s}(\x)))\cdot\nabla^2_{z_1}]
$$
and
$$
\sB_r:=(F_1(r,\z)-F_1(r,\btheta_{r,s}(\x)))\cdot \nabla_{z_1},\ \ \sK_r:=\cT_{F_2(r)}(\z,\btheta_{r,s}(\x))\cdot \nabla_{z_2}.
$$
Using the above notations, we can rewrite \eqref{Duhamel_change_freez} as
\begin{align}
\nabla_{x_2}P_{s,t}f(\x)=
Q_{s,u}P_{u,t}f(\x)+\int_s^{u}Q_{s,r}(\sD_{r}+\sB_{r}+\sK_{r})P_{r,t}f(\x) \dif r.
\label{Duh}
\end{align}
From \eqref{CAN}, the upper bound for $p$ and \eqref{Convolution}, we have
\begin{align*}
|Q_{s,u}P_{u,t}f(\x)|
&\lesssim (u-s)^{-\frac 32}\hat P^\lambda_{s,u}P_{u,t}f(\x)
\lesssim (t-s)^{-\frac 32}\hat P^\lambda_{s,t}f(\x).
\end{align*}
To fully make up the time singularity in \eqref{CAN}, we exploit the cancellation property \eqref{CAN0}. For $r\in[s,u]$, let $\textcolor{black}{\Phi_r}:\mR^d\to[0,1]$ be a \textcolor{black}{\textit{smooth}} cutoff function with
\begin{align}\label{CUT_OFF_CAN}
\textcolor{black}{\Phi}_r(z_2)=\left\{
\begin{aligned}
&1,\ \ |z_2-\btheta^2_{r,s}(\x)|\le (t-r)^{\frac32},\\
&0,\ \ |z_2-\btheta^2_{r,s}(\x)|\ge 2(t-r)^{\frac32},\ 
\end{aligned}
\right.\ \ 
|\nabla_{z_2}\textcolor{black}{\Phi_r}(z_2)|\leq 4(t-r)^{-\frac32}.
\end{align}
First of all we look at the term containing $\sD_{r}$ in \eqref{Duh} and write
\begin{align*}
Q_{s,r}\sD_{r}P_{r,t}f(\x)&=-\int_{\mR^{2d}}\nabla_{z_2}p_0(s,\x;r,\z)(1-\textcolor{black}{\Phi}_r(z_2))\sD_{r}P_{r,t}f(\z)\dif\z\\
&\quad+\int_{\mR^{2d}}p_0(s,\x;r,\z)\nabla_{z_2} \textcolor{black}{\Phi}_r(z_2)\sD_{r}P_{r,t}f(\z)\dif\z\\
&\quad-\int_{\mR^{2d}}\nabla_{z_2}(p_0(s,\x;r,\z)\textcolor{black}{\Phi}_r(z_2))\sD_{r}P_{r,t}f(\z)\dif\z\\
&=:I_1+I_2-I_3.
\end{align*}
For $I_1$, by \eqref{CAN} and the upper bound estimate \eqref{EST_GRAD_X1} of $p$, \textcolor{black}{since $f$ is nonnegative}, we have
\begin{align*}
|I_1|&\lesssim(r-s)^{-\frac32}\int_{\mR^{2d}}\hat p_\lambda(s,\x;r,\z) (1-\textcolor{black}{\Phi}_r(z_2))|\z-\btheta_{r,s}(\x)|_{\bf d}^\gamma |\nabla^2_{z_1}P_{r,t}f(\z)|\dif\z\\
&\lesssim(r-s)^{\frac{\gamma-3}2}(t-r)^{-1}\int_{\mR^{2d}}\hat p_\lambda(s,\x;r,\z)\1_{\{|z_2-\btheta^2_{r,s}(\x)|\ge (t-r)^{3/2}\}} \hat P^\lambda_{r,t}f(\z)\dif\z\\
&\leq(r-s)^{\frac{\gamma-3}2}(t-r)^{-\frac{3}{2}}\int_{\mR^{2d}}|z_2-\btheta^2_{r,s}(\x)|^{\frac13}\hat p_\lambda(s,\x;r,\z) \hat P^\lambda_{r,t}f(\z)\dif\z\\
&\lesssim(r-s)^{\frac{\gamma}2-1}(t-r)^{-\frac{3}{2}}\int_{\mR^{2d}}\hat p_\lambda(s,\x;r,\z) \hat P^\lambda_{r,t}f(\z)\dif\z\\
&\lesssim(r-s)^{\frac{\gamma}2-1}(t-s)^{-\frac{3}{2}}\hat P^\lambda_{s,t}f(\x).
\end{align*}
Similarly, for $I_2$, we have
\begin{align*}
|I_2|&\lesssim \int_{\mR^{2d}}\hat p_\lambda(s,\x;r,\z) |\nabla_{z_2}\textcolor{black}{\Phi}_r(z_2)|\,|\z-\btheta_{r,s}(\x)|_{\bf d}^\gamma |\nabla^2_{z_1}P_{r,t}f(\z)|\dif\z\\
&\lesssim (t-r)^{-\frac52}(r-s)^{\frac\gamma2}\int_{\mR^{2d}}\hat p_\lambda(s,\x;r,\z)\hat P^\lambda_{r,t}f(\z)\dif\z\\
&\lesssim(t-s)^{\frac{\gamma-5}2}\hat P^\lambda_{s,t}f(\x).
\end{align*}
For $I_3$, \textcolor{black}{similarly to} the \textcolor{black}{partial} cancellation \eqref{CAN0} we can write
\begin{align*}
I_3&=\int_{\mR^{2d}}\nabla_{z_2}(p_0(s,\x;r,\z)\textcolor{black}{\Phi}_r(z_2))
\Big[\tr\left( \left(a(r,\z)-a(r,\btheta_{r,s}(\x))\right) \cdot\nabla^2_{z_1} P_{r,t}f(\z)\right)\\
&- \tr\left( \left(a(r,z_1,\btheta^2_{r,s}(\x))-a(r,\btheta_{r,s}(\x))\right) 
\cdot\nabla^2_{z_1} P_{r,t}f(z_1,\btheta^2_{r,s}(\x))\right)\Big]\dif\z\\
&=\int_{\mR^{2d}}\nabla_{z_2}(p_0(s,\x;r,\z)\textcolor{black}{\Phi}_r(z_2))
\Big[\tr\left((a(r,\z)-a(r,z_1,\btheta^2_{r,s}(\x)))\cdot \nabla^2_{z_1} P_{r,t}f(\z)\right)\\
&+ \tr\big( (a(r,z_1,\btheta^2_{r,s}(\x))-a(r,\btheta_{r,s}(\x)))\cdot
\nabla^2_{z_1}(P_{r,t}f(\z)-P_{r,t}\textcolor{black}{f}(z_1,\btheta^2_{r,s}(\x)))\big)\Big]\dif\z.
\end{align*}
Thus, by \eqref{CAN}, \textcolor{black}{\eqref{CUT_OFF_CAN}}, \eqref{HOLD_MOD_D_X1_2} and \eqref{Re1}, for any $\eta\in(0,\gamma)$, we have
\begin{align*}
|I_{\textcolor{black}{3}}|&\lesssim\Big[(r-s)^{-\frac32}+(t-s)^{-\frac32}\Big]\int_{\mR^{2d}}\hat p_\lambda(s,\x;r,\z)\1_{\{|z_2-\btheta^2_{r,s}(\x)|\le 2(t-r)^{3/2}\}}\\
&\times\Big[|z_2-\btheta^2_{r,s}(\x)|^{\frac{1+\gamma}{3}}|\nabla^2_{z_1}P_{r,t}f(\z)|+|z_1-\btheta^1_{r,s}(\x)|^\alpha
|\nabla^2_{z_1}(P_{r,t}f(\z)-P_{r,t}\textcolor{black}{f}(z_1,\btheta^2_{r,s}(\x)))|\Big]\dif\z\\
&\lesssim\Big[(r-s)^{-\frac32}+(t-s)^{-\frac32}\Big]\int_{\mR^{2d}}\hat p_\lambda(s,\x;r,\z)\1_{\{|z_2-\btheta^2_{r,s}(\x)|\le 2(t-r)^{3/2}\}}\\
&\quad\times\Big[|z_2-\btheta^2_{r,s}(\x)|^{\frac{1+\gamma}{3}}(t-r)^{-1}\hat P^\lambda_{r,t}f(\z)+|z_1-\btheta^1_{r,s}(\x)|^\alpha|z_2-\btheta^2_{r,s}(\x)|^{\frac\eta3}\\
&\qquad\times(t-r)^{-1-\frac\eta2}(\hat P^\lambda_{r,t}f(\z)+\hat P^\lambda_{r,t}f(z_1,\btheta^2_{r,s}(\x)))\Big]\dif\z\\
&\lesssim\Big[(r-s)^{-\frac32}+(t-s)^{-\frac32}\Big]\Big[(r-s)^{\frac{\gamma+1}2}(t-r)^{-1}+(r-s)^{\frac{\alpha+\eta}2}(t-r)^{-1-\frac\eta2}\Big]\\
&\qquad\times\int_{\mR^{2d}}\hat p_\lambda(s,\x;r,\z)\hat P^\lambda_{r,t}f(\z)\dif z,
\end{align*}
where the last step is due to \eqref{regularizing_proxy} and \eqref{VX2}. Since $\alpha\in((1-\gamma)\vee\gamma,1]$, one can choose $\eta$ close to $\gamma$ so that for some $\eps>0$,
$$
|I_{\textcolor{black}{3}}|\lesssim[(r-s)^{\eps-1}+(t-s)^{\eps-1}](t-s)^{-\frac32}\hat P^\lambda_{s,t}f(\x).
$$
Next we treat the term containing $\sB_{r}$, and similarly write
\begin{align*}
Q_{s,r}\sB_{r}P_{r,t}f(\x)&=-\int_{\mR^{2d}}\nabla_{z_2}p_0(s,\x;r,\z)(1-\textcolor{black}{\Phi}_r(z_2))\sB_{r}P_{r,t}f(\z)\dif\z\\
&\quad+\int_{\mR^{2d}}p_0(s,\x;r,\z)\nabla_{z_2} \textcolor{black}{\Phi}_r(z_2)\sB_{r}P_{r,t}f(\z)\dif\z\\
&\quad-\int_{\mR^{2d}}\nabla_{z_2}(p_0(s,\x;r,\z)\textcolor{black}{\Phi}_r(z_2))\sB_{r}P_{r,t}f(\z)\dif\z\\
&=:J_1+J_2-J_3.
\end{align*}
For $J_1$ and $J_2$, the \textcolor{black}{analysis is similar to the one performed for $I_1$ and $I_2$}. By \eqref{CAN} and the upper bound estimate \eqref{EST_GRAD_X1} of $p$, we have
\begin{align*}
|J_1|+|J_2|&\lesssim\Big[(r-s)^{\frac{\gamma}2-1}(t-r)^{-1}+(t-s)^{\frac\gamma2-2}\Big]\hat P^\lambda_{s,t}f(\x).
\end{align*}
For $J_3$,  as for $I_3$, \textcolor{black}{similarly to} \eqref{CAN0} we can write
\begin{align*}
J_{\textcolor{black}{3}}&=\int_{\mR^{2d}}\nabla_{z_2}(p_0(s,\x;r,\z)f_r(z_2))
\Big[(F_1(r,\z)-F_1(r,z_1,\btheta^2_{r,s}(\x)))\cdot \nabla_{z_1}P_{r,t}f(\z)\\
&\quad+  (F_1(r,z_1,\btheta^2_{r,s}(\x))-F_1(r,\btheta_{r,s}(\x)))\cdot
\nabla_{z_1}(P_{r,t}f(\z)-P_{r,t}\textcolor{black}{f}(z_1,\btheta^2_{r,s}(\x)))\Big]\dif\z.
\end{align*}
Thus, by \eqref{CAN}, \textcolor{black}{\eqref{CUT_OFF_CAN}}, \eqref{HOLD_MOD_D_X1_2} and \textcolor{black}{\eqref{Re2}}, for any $\eta\in(0,1)$, we have
\begin{align*}
|J_3|&\lesssim\Big[(r-s)^{-\frac32}+(t-s)^{-\frac32}\Big]\int_{\mR^{2d}}\hat p_\lambda(s,\x;r,\z)\1_{\{|z_2-\btheta^2_{r,s}(\x)|\le 2(t-r)^{3/2}\}}\\
&\times\Big[|z_2-\btheta^2_{r,s}(\x)|^{\frac{1+\gamma}{3}}|\nabla_{z_1}P_{r,t}f(\z)|+|z_1-\btheta^1_{r,s}(\x)|^\gamma
|\nabla_{z_1}(P_{r,t}f(\z)-P_{r,t}\textcolor{black}{f}(z_1,\btheta^2_{r,s}(\x)))|\Big]\dif\z\\
&\lesssim\Big[(r-s)^{-\frac32}+(t-s)^{-\frac32}\Big]\int_{\mR^{2d}}\hat p_\lambda(s,\x;r,\z)\1_{\{|z_2-\btheta^2_{r,s}(\x)|\le 2(t-r)^{3/2}\}}\\
&\quad\times\Big[|z_2-\btheta^2_{r,s}(\x)|^{\frac{1+\gamma}{3}}(t-r)^{-\textcolor{black}{\frac 12}}\hat P^\lambda_{r,t}f(\z)+|z_1-\btheta^1_{r,s}(\x)|^\gamma|z_2-\btheta^2_{r,s}(\x)|^{\frac\eta3}\\
&\qquad\times(t-r)^{-\textcolor{black}{\frac 12}-\frac\eta2}(\hat P^\lambda_{r,t}f(\z)+\hat P^\lambda_{r,t}f(z_1,\btheta^2_{r,s}(\x)))\Big]\dif\z\\
&\lesssim\Big[(r-s)^{-\frac32}+(t-s)^{-\frac32}\Big]\Big[(r-s)^{\frac{\gamma+1}2}(t-r)^{-\textcolor{black}{\frac 12}}+(r-s)^{\frac{\gamma+\eta}2}(t-r)^{-\textcolor{black}{\frac 12}-\frac\eta2}\Big]\\
&\qquad\times\int_{\mR^{2d}}\hat p_\lambda(s,\x;r,\z)\hat P^\lambda_{r,t}f(\z)\dif z,
\end{align*}
In particular, one can choose $\eta$ close to $1$ so that for some $\eps>0$,
$$
|J_{\textcolor{black}{3}}|\lesssim  [(r-s)^{\eps-\textcolor{black}{ 1}}+(t-s)^{\eps-\textcolor{black}{ 1}}](t-s)^{-\textcolor{black}{1}}\hat P^\lambda_{s,t}f(\x).
$$
\textcolor{black}{The term $Q_{s,r}\sB_{r}P_{r,t}f(\x)$ is therefore not critical in terms of the time singularities, since it appears to be less singular than $Q_{s,r}\sD_{r}P_{r,t}f(\x)$. However, the additional regularity assumption \eqref{Re2} is really needed to derive an integrable singularity in the variable $r$.}

On the other hand, we have
\begin{align*}
|Q_{s,r}\sK_{r}P_{r,t}f(\x)|
&\lesssim(r-s)^{-\frac32}\int_{\mR^{2d}}\hat p_\lambda(s,\x;r,\z)|\z-\btheta_{r,s}(\x)|_{\bf d}^{1+\gamma}|\nabla_{z_2}P_{r,t}f(\z)|\dif\z\\
&\lesssim(r-s)^{\frac{\gamma}{2}-1}\int_{\mR^{2d}}\hat p_\lambda(s,\x;r,\z)|\nabla_{z_2}P_{r,t}f(\z)|\dif\z\\
&=(r-s)^{\frac{\gamma}{2}-1}\hat P^\lambda_{s,r}|\nabla_{z_2}P_{r,t}f|(\x).
\end{align*}
Gathering together all the \textcolor{black}{previous} controls, we eventually get
\begin{align}
|\nabla_{x_2}P_{s,t}f(\x)|&\lesssim(t-s)^{-\frac 32}\hat P^\lambda_{s,t}f(\x)+\int^u_s(r-s)^{\frac{\gamma}{2}-1}\hat P^\lambda_{s,r}|\nabla_{z_2}P_{r,t}f|(\x)\dif r\no\\
&\lesssim(t-s)^{-\frac 32}\bar P^\delta_{s,t}f(\x)+\int^u_s(r-s)^{\frac{\gamma}{2}-1}\bar P^\delta_{s,r}|\nabla_{z_2}P_{r,t}f|(\x)\dif r,\label{HK33}
\end{align}
where $\delta$ is chosen as in \eqref{HK2}.
Fix $s_0<t$. For any $s\in(s_0,t]$, define
$$
\Phi^f_{s_0,t}(s):=(t-s)^{\frac 32}\bar P^\delta_{s_0,s}|\nabla_{x_2}P_{s,t}f|(\x).
$$
Using $\bar P^\delta_{s_0,s}$ act on both sides of \eqref{HK33} and by the \textcolor{black}{Chapman-Kolmogorov equation}, we obtain
\begin{align*}
\Phi^f_{s_0,\textcolor{black}{t}}(s)&\lesssim\bar P^\delta_{s_0,t}f(\x)+(t-s)^{\frac 32}\int^u_s(r-s)^{\frac{\gamma}{2}-1}\bar P^\delta_{s_0,r}|\nabla_{z_2}P_{r,t}f|(\x)\dif r\\
&\leq\bar P^\delta_{s_0,t}f(\x)+(t-s)^{\frac 32}\int^u_s(r-s)^{\frac{\gamma}{2}-1}(t-r)^{-\frac 32}\Phi^f_{s_0,\textcolor{black}{t}}(\textcolor{black}{r})\dif r\\
&\leq\bar P^\delta_{s_0,t}f(\x)+2^{\frac32}\int^u_s(r-s)^{\frac{\gamma}{2}-1}\Phi^f_{s_0,\textcolor{black}{t}}(r)\dif r\\
&\leq\bar P^\delta_{s_0,t}f(\x)+2^{\frac32}\int^t_s(r-s)^{\frac{\gamma}{2}-1}\Phi^f_{s_0,\textcolor{black}{t}}(r)\dif r.
\end{align*}
By the Volterra type Gronwall inequality we obtain 
\begin{equation}
\Phi^f_{s_0,\textcolor{black}{t}}(s)\lesssim\bar P^\delta_{s_0,t}f(\x) \; \Longrightarrow \; \bar P^\delta_{s_0,s}|\nabla_{x_2}P_{s,t}f|(\x)\lesssim(t-s)^{-\frac32}\bar P^\delta_{s_0,t}f(\x).
\end{equation}
Letting $s_0\uparrow s$, we get
$$
|\nabla_{x_2}P_{s,t}f|(\x)\lesssim(t-s)^{-\frac32}\bar P^\delta_{s,t}f(\x).
$$
Finally, for fixed $t'>t$ and $\y\in\mR^{2d}$, we let $f(\x):=p(t,\x;t',\y)\in C^\infty_b(\mR^{2d})$, then by the Chapman-Kolmogorov equation and \eqref{Convolution}, we obtain
$$
|\nabla_{x_2} p(s,\x;t',\y)|\lesssim_{C_3} (t-s)^{-\frac32}\hat p_\lambda(s,\x;t',\y).
$$
This then readily gives estimate \eqref{EST_GRAD_X2}. The proof is complete.
\endproof

\bp\label{gradient_x2_HOLD}
Under {\bf (H$^\gamma_\sigma$)}, {\bf (H$^{\gamma}_{\F}$)}, \eqref{Re1}, \eqref{Re2}  and \eqref{Smo1}, for any $T>0$,
there exist constants $\lambda, C\geq 1$ and $\eta\in(0,(\alpha-\gamma)\wedge (\alpha+\gamma-1))$,depending on $\Theta_T$ and $\bar\kappa_0, \bar\kappa_1$ such that for any $0\le s<t\le T$ and $\x,\y\in\R^{2d}$, 
\begin{align}
|\nabla_{x_2}(p(s,\x;t,\y)-p(s,\x'; t,\y))|&\le C|\x-\x'|_{\bf d}^{\eta}(t-s)^{-\frac{3+\eta}{2}}\(\hat p_{\lambda}(s,\x;t,\y)+ \hat p_{\lambda}(s,\x';t,\y)\).
\label{HOLD_MOD_D_X2}
\end{align}
\ep
\proof 
If $|\x-\x'|^2_{\bf d}>t-s$, \textcolor{black}{then} by \eqref{EST_GRAD_X2} we clearly have
\begin{equation*}
|\nabla_{x_2}(p(s,\x;t,\y)-p(s,\x'; t,\y))|\lesssim (t-s)^{-\frac{3}{2}}(\hat p_{\lambda}(s,\x;t,\y)+ \hat p_{\lambda}(s,\x';t,\y))\lesssim \text{r.h.s. of \eqref{HOLD_MOD_D_X2}}.
\end{equation*}
Next we restrict to the \textit{global diagonal case} $|\x-\x'|^2_{\bf d}\le t-s$.
We first write a localized version of the Duhamel formula.
For any freezing couple $(\bar \tau, \bar \bxi)$ we have taking formally $f=\delta_\y $ in \eqref{SG_LOC_IN_U}
\begin{equation}\label{Duhamel}
p(s,\x;t,\y)=\tilde{P}^{(\bar\tau,\bar\bxi)}_{\textcolor{black}{s},u}p(\textcolor{black}{u},\cdot;t,\y)(\x)+\int_s^u\int_{\R^{2d}}\tilde{p}^{(\bar\tau,\bar\bxi)}(s,\x;r,\z)
\bar{\mathcal{L}}^{(\bar\tau,\bar\bxi)}p(r,\z;t,\y)\dif\z\dif r, 
\end{equation}
where $\wt{P}^{(\bar\tau,\bar\bxi)}_{t,u}f $ is as in \eqref{Semigr_tilde} and we denoted 
$\bar{\mathcal{L}}^{(\bar\tau,\bar\bxi)}=\mathcal{L}-\mathcal{L}^{(\bar\tau,\bar\bxi)}$.
Let us now differentiate w.r.t $u$: we obtain 
\begin{equation}\label{Duhamel_loc}
0=\p_u[\wt{P}^{(\bar\tau,\bar\bxi)}_{\textcolor{black}{s},u}p(\textcolor{black}{u},\cdot;t,\y)(\x)]+\int_{\R^{2d}}\tilde{p}^{(\bar\tau,\bar\bxi)}(s,\x;u,\z)
\bar{\mathcal{L}}^{(\bar\tau,\bar\bxi)}p(u,\z;t,\y)\dif\z.
\end{equation}
We integrate the previous equation taking: 
\begin{enumerate}
\item $(\bar \tau, \bar \bxi)=(\tau_0,\bxi_0)$ between $s$ and $s_1=\frac{t+s}{2}$; 
\item $(\bar \tau, \bar \bxi)=(\tau_1,\bxi_1)$ between $s_1$ and $t$;
\end{enumerate}
Then we get 
\begin{align*}
0&=\tilde{P}^{(\tau_0,\bxi_0)}_{s,s_1}p(s_1,\cdot;t,\y)(\x)-p(s,\x;t,\y)
+\int_s^{s_1}\int_{\R^{2d}}\tilde{p}^{(\tau_0,\bxi_0)}(s,\x;r,\z)\bar{\mathcal{L}}^{(\tau_0,\bxi_0)}p(r,\z;t,\y)\dif\z\dif r,\\
0&=\tilde{p}^{(\tau_1,\bxi_1)}(s,\x;t,\y)-\tilde{P}^{(\tau_1,\bxi_1)}_{s,s_1}p(s_1,\cdot;t,\y)(\x)
+\int_{s_1}^t\int_{\R^{2d}}\tilde{p}^{(\tau_1,\bxi_1)}(s,\x;r,\z)\bar{\mathcal{L}}^{(\tau_1,\bxi_1)}p(r,\z;t,\y)\dif\z\dif r.
\end{align*}
Next we use expansion \eqref{Duhamel_loc} for $p(t,\x';s,\y)$ and integrate the equation taking: 
\begin{enumerate}
\item $(\bar \tau, \bar \bxi)=(\textcolor{black}{\tau_0,\bxi_0'})$ between $s$ and $s_0=t+c\textcolor{black}{|\x-\x'|_{\mathbf{d}}^2}$; 
\item  $(\bar \tau, \bar \bxi)=(\textcolor{black}{\tau_0,\bxi_0})$ between $s_0$ and $s_1$;
\item  $(\bar \tau, \bar \bxi)=(\tau_1,\bxi_1)$ between $s_1$ and $t$;
\end{enumerate}
Then we get 
\begin{align*}
0&=\tilde{P}^{(\textcolor{black}{\tau_0,\bxi_0'})}_{s,s_0}p(s_0,\cdot;t,\y)(\x')-p(s,\x';t,\y)
+\int_s^{s_0}\int_{\R^{2d}}\tilde{p}^{(\textcolor{black}{\tau_0,\bxi_0'})}(s,\x';r,\z)\bar{\mathcal{L}}^{(\textcolor{black}{\tau_0,\bxi_0'})}p(r,\z;s,\y)\dif\z\dif r,\\
0&=\tilde{P}^{(\textcolor{black}{\tau_0,\bxi_0})}_{s,s_1}p(s_1,\cdot;t,\y)(\x')-\tilde{P}^{(\textcolor{black}{\tau_0,\bxi_0})}_{s,s_0}p(s_0,\cdot;t,\y)(\x')
+\int_{s_0}^{s_1}\int_{\R^{2d}}\tilde{p}^{(\textcolor{black}{\tau_0,\bxi_0})}(s,\x';r,\z)\bar{\mathcal{L}}^{(\textcolor{black}{\tau_0,\bxi_0})}p(r,\z;t,\y)\dif\z\dif r.\\
0&=\tilde{p}^{(\tau_1,\bxi_1)}(s,\x';t,\y)-\tilde{P}^{(\tau_1,\bxi_1)}_{s,s_1}p(s_1,\cdot;t,\y)(\x')
+\int_{s_1}^t\int_{\R^{2d}}\tilde{p}^{(\tau_1,\bxi_1)}(s,\x';r,\z)\bar{\mathcal{L}}^{(\tau_1,\bxi_1)}p(r,\z;t,\y)\dif\z\dif r.
\end{align*}

Notice that it \textcolor{black}{suffices} to take $0<c<\frac{1}{2}$ to ensure $s_0<s_1$. \textcolor{black}{We then  have:} 
\begin{align*}
&p(s,\x;t,\y)-p(s,\x';t,\y)\\
&\qquad =\tilde{p}^{(\tau_1,\bxi_1)}(s,\x;t,\y)-\tilde{p}^{(\tau_1,\bxi_1)}(s,\x';t,\y)\\
&\qquad \quad +\Delta\tilde{P}^{(\tau_0,\bxi_0,\bxi'_0)}(s_0,t;\x',\x',\y)+
\Delta\tilde{P}^{(\tau_0,\bxi_0,\bxi_0)}(s_1,t;\x,\x',\y)+\Delta\tilde{P}^{(\tau_1,\bxi_1,\bxi_1)}(s_1,t;\x',\x,\y)\\
&\qquad \quad+\Delta_{\rm off-diag}^{(\tau_0,\bxi_0,\bxi'_0)}(s,t;\x,\x',\y)
+\Delta_{\rm diag}^{(\tau_0,\bxi_0,\bxi_0)}(s,t;\x,\x',\y)
\end{align*}
where
\begin{align*}
\Delta\tilde{P}^{(\tau_0,\bxi_0,\bxi'_0)}(s_0,t;\x',\x',\y)&=\textcolor{black}{-}
\tilde{P}^{(\tau_0,\bxi'_0)}_{s,s_0}p(s_0,\cdot;t,\y)(\x')\textcolor{black}{+}\tilde{P}^{(\tau_0,\bxi_0)}_{s,s_0}p(s_0,\cdot;t,\y)(\x'),\\
\Delta\tilde{P}^{(\tau_0,\bxi_0,\bxi_0)}(s_1,t;\x,\x',\y)&=
\tilde{P}^{(\tau_0,\bxi_0)}_{s,s_1}p(s_1,\cdot;t,\y)(\x)-\tilde{P}^{(\tau_0,\bxi_0)}_{s,s_1}p(s_1,\cdot;t,\y)(\x'),\\
\Delta\tilde{P}^{(\tau_1,\bxi_1,\bxi_1)}(s_1,t;\x',\x,\y)&=
\tilde{P}^{(\tau_1,\bxi_1)}_{s,s_1}p(s_1,\cdot;t,\y)(\x')-\tilde{P}^{(\tau_1,\bxi_1)}_{s,s_1}p(s_1,\cdot;t,\y)(\x),
\end{align*}
and
\begin{align*}
\Delta_{\rm off-diag}^{(\tau_0,\bxi_0,\bxi'_0)}(s,t;\x,\x',\y)&=\int_s^{s_0}\int_{\R^{2d}}\tilde{p}^{(\tau_0,\bxi_0)}(s,\x;r,\z)
\bar{\mathcal{L}}^{(\tau_0,\bxi_0)}p(r,\z;t,\y)\dif\z\dif r\\
&\qquad - \int_s^{s_0}\int_{\R^{2d}}\tilde{p}^{(\tau_0,\bxi'_0)}(s,\x';r,\z)\bar{\mathcal{L}}^{(\tau_0,\bxi'_0)}p(r,\z;t,\y)\dif\z\dif r,\\
\Delta_{\rm diag}^{(\tau_0,\bxi_0,\tau_1,\bxi_1)}(s,t;\x,\x',\y)&=\int_{s_0}^{s_1}
\int_{\R^{2d}}(\tilde{p}^{(\tau_0,\bxi_0)}(s,\x;r,\z)-\tilde{p}^{(\tau_0,\bxi_0)}(s,\x';r,\z))\bar{\mathcal{L}}^{(\tau_0,\bxi_0)}p(r,\z;t,\y)\dif\z\dif r,\\
&\qquad +\int_{s_1}^{t}
\int_{\R^{2d}}(\tilde{p}^{(\tau_1,\bxi_1)}(s,\x;r,\z)-\tilde{p}^{(\tau_1,\bxi_1)}(s,\x';r,\z))\bar{\mathcal{L}}^{(\tau_1,\bxi_1)}p(r,\z;t,\y)\dif\z\dif r\\
&=:\Delta_{\rm diag,1}^{(\tau_0,\bxi_0)}(s,t;\x,\x',\y)+\Delta_{\rm diag,2}^{(\tau_1,\bxi_1)}(s,t;\x,\x',\y)\textcolor{black}{.}
\end{align*}
Let us start with the term $\Delta_{\rm off-diag}^{(\tau_0,\bxi_0,\bxi'_0)}(s,t;\x,\x',\y)$. \textcolor{black}{Differentiating} with respect to $x_2$ and taking $(\tau,\bxi_0,\bxi'_0)=(s,\x,\x')$ \textcolor{black}{\textbf{after}} differentiation, we have 
\begin{align*}
&|\nabla_{x_2}\Delta_{\rm off-diag}^{(\tau_0,\bxi_0,\bxi'_0)}(s,t;\x,\x',\y)|_{(\tau_0,\bxi_0,\bxi'_0)=(s,\x,\x')}\\
&\quad \leq \left|\int_s^{s_0}\int_{\R^{2d}}\nabla_{\x_2} \tilde{p}^{(\tau_0,\bxi_0)}(s,\x;r,\z)\bar{\mathcal{L}}^{(\tau_0,\bxi_0)}p(r,\z;t,\y)\dif\z\dif r\right|_{(\tau_0,\bxi_0)=(s,\x)}\\
&\qquad +\left|\int_s^{s_0}\int_{\R^{2d}}\nabla_{\x_2} \tilde{p}^{(\tau_0,\bxi'_0)}(s,\x';r,\z)\bar{\mathcal{L}}^{(\tau_0,\bxi'_0)}p(r,\z;t,\y)\dif\z\dif r\right|_{(\tau_0,\bxi'_0)=(s,\x')},
\end{align*}
where both the terms in the r.h.s. can be controlled separately as in the proof of Proposition \ref{gradient_x2}. We derive: 
\begin{align*}
&|\nabla_{x_2}\Delta_{\rm off-diag}^{(\tau_0,\bxi_0,\bxi'_0)}(s,t;\x,\x',\y)|_{(\tau_0,\bxi_0,\bxi'_0)=(s,\x,\x')}\\
&\quad \lesssim (\bar{p}(s,\x;t,\y)+\bar{p}(s,\x';t,\y))\int_s^{s+c|\x-\x'|_{\bf{d}}^2}{\[(r-s)^{-1+\frac{\gamma}{2}}+(r-s)^{\frac{\alpha+\eta-3}{2}}\]}(t-r)^{-\frac{3}{2}}\dif r
\intertext{(since $\alpha>(\gamma\vee 1-\gamma)$, we can choose any $0<\epsilon< (\alpha-\gamma) \wedge (\alpha+\gamma-1)$ such that)}
&\quad \lesssim (\bar{p}(s,\x;t,\y)+\bar{p}(s,\x';t,\y))\int_s^{s+c|\x-\x'|_{\bf{d}}^2}(r-s)^{-1+\frac{\epsilon}{2}}(t-r)^{-\frac{3}{2}}\dif r\\
&\quad \lesssim (t-r)^{-\frac{3}{2}}|\x-\x'|_{\bf{d}}^{\epsilon}(\bar{p}(s,\x;t,\y)+\bar{p}(s,\x';t,\y)).
\end{align*}
Let us now turn to $\Delta_{\rm diag}^{(\tau_0,\bxi_0,\tau_1,\bxi_1)}(s,t;\x,\x',\y)$: similarly to \eqref{proxy_regularity}, in the \textit{local diagonal case} $|\x-\x'|_{\bf{d}}^2\le r-s$ we have, for any $\alpha\in [0,1]$,
\begin{align*}
&|\nabla_{x_2}\left(\wt p^{(\tau_0,\bxi_0)}(s,\x;r,\z)-\wt p^{(\tau_0,\bxi_0)}(s,\x'; r,\z)\right)|_{(\tau_0,\bxi_0)=(s,\x)}\\
& \qquad \le \sum_{i=1}^2|x_i-x'_i|\sup_{\eta\in [0,1]}|\nabla_{x_i}\nabla_{x_2}\wt p^{(\tau_0,\bxi_0)}(s,\x+\eta(\x'-\x); r,\z)|_{(\tau_0,\bxi_0)=(s,\x)}\\
&\qquad \lesssim|\x-\x'|_{\bf d}^{\alpha}\,(r-s)^{-\frac{3+\alpha}{2}}\hat{p}_{\lambda}(s,\x;r,\z).
\end{align*}
Therefore we have 
\begin{align*}
&|\nabla_{x_2}\Delta_{\rm diag,1}^{(\tau_0,\bxi_0)}(s,t;\x,\x',\y)|_{(\tau_0,\bxi_0)=(s,\x)}\\
&\qquad \lesssim |\x-\x'|_{\bf{d}}^{\eta}\int_{s_0}^{s_1}(r-s)^{-\frac{3+\eta}{2}}\int_{\R^{2d}}\hat{p}_{\lambda}(s,\x;r,\z)|\bar{\mathcal{L}}^{(s,\x)}p(r,\z;t,\y)|\dif\z\dif r, 
\end{align*}
which is controlled again, as in the proof of Proposition \ref{gradient_x2}. We eventually derive, 
for any $\epsilon$ as above, and $\eta<\epsilon$ 
\begin{align*}
|\nabla_{x_2}\Delta_{\rm diag,1}^{(\tau_0,\bxi_0)}(s,t;\x,\x',\y)|_{(\tau_0,\bxi_0)=(s,\x)}
&\lesssim |\x-\x'|_{\bf{d}}^{\eta}\hat{p}_{\lambda}(s,\x;t,\y)\int_{s_0}^{s_1}(r-s)^{-\frac{1+\eta-\epsilon}{2}}(t-r)^{-\frac{3}{2}}\dif r\\
&\lesssim |\x-\x'|_{\bf{d}}^{\eta}(t-s)^{-\frac{3}{2}}\hat{p}_{\lambda}(s,\x;t,\y).
\end{align*}
\textcolor{black}{The control for} $\Delta_{\rm diag,2}^{(\tau_1,\bxi_1)}(s,t;\x,\x',\y)$ is more direct. Indeed, taking $(\tau_1,\bxi_1)=(t,\y)$ we have, by \eqref{proxy_regularity}, 
\begin{align*}
&|\nabla_{x_2}\Delta_{\rm diag,2}^{(\tau_1,\bxi_1)}(s,t;\x,\x',\y)|_{(\tau_1,\bxi_1)=(t,\y)}\\
&\qquad \lesssim \int_{s_1}^{t}
\int_{\R^{2d}}\frac{|\x-\x'|^{\eta}_{\bf{d}}}{(r-s)^{\frac{3+\eta}{2}}}(\hat{p}_{\lambda}(s,\x;r,\z)+\hat{p}_{\lambda}(s,\x';r,\z))
\frac{\hat{p}_{\lambda}(r,\z;t,\z)}{(t-r)^{1-\frac{\gamma}{2}}}\dif\z\dif r\\
&\qquad \lesssim |\x-\x'|^{\eta}_{\bf{d}}(t-s)^{-\frac{3+\eta}{2}}(\hat{p}_{\lambda'}(s,\x;t,\y)+\hat{p}_{\lambda'}(s,\x';t,\y))
\end{align*}
where we used \eqref{Convolution} in the last inequality. 

It remains to check the terms $\Delta\tilde{P}$, which arise from the change of freezing couples. 
By \eqref{THE_CTR_SENSI_DIFF_PARAMbis}, we have 
\begin{align*}
&|\nabla_{x_2}\Delta\tilde{P}^{(\tau_0,\bxi_0,\bxi'_0)}(s_0,t;\x',\x',\y)|_{(\tau_0,\bxi_0,\bxi'_0)=(s,\x,\x')}\\
& \qquad =\left|\nabla_{x_2}\int_{\R^{2d}}\left(\tilde{p}^{(\tau_0,\bxi_0)}(s,\x';s_0,\z)-\tilde{p}^{(\tau_0,\bxi_0')}(s,\x';s_0,\z)\right)
p(s_0,\z;t,\y)\dif\z \right|_{(\tau_0,\bxi_0,\bxi'_0)=(s,\x,\x')}\\
&\qquad \lesssim |\x-\x'|^{\gamma-\epsilon}\int_{\R^{2d}}(s_0-s)^{-\frac{3-\epsilon}{2}}\hat{p}_{\lambda}(s,\x';s_0,\z)
(p(s_0,\z;t,\y)-p(s_0,\z_1,\btheta^2_{s_0,s}(\x');t,\y))\dif\z\\
&\qquad \lesssim\frac{|\x-\x'|^{\gamma-\epsilon}}{(s_0-s)^{\frac{3-\epsilon}{2}}}
\int_{\R^{2d}}\hat{p}_{\lambda}(s,\x';s_0,\z)\frac{|\z_2-\btheta^2_{s_0,s}(\x')|}{(t-s_0)^{\frac{3}{2}}}
\sup_{\rho\in [0,1]}\hat{p}_{\lambda}(s_0,\z_1,\z_2+\rho(\btheta^2_{s_0,s}(\x')-\z_2);t,\y)\dif\z
\intertext{(by \eqref{regularizing_proxy}, choosing $\epsilon<\gamma$, and since $2(t-s_0)>t-s$)}
&\qquad \lesssim\frac{|\x-\x'|^{\gamma-\epsilon}}{(t-s)^{\frac{3}{2}}}\hat{p}_{\lambda'}(s,\x';t,\y).
\end{align*}
The remaining terms  $\Delta\tilde{P}^{(\tau_0,\bxi_0,\bxi_0)}(s_1,t;\x,\x',\y)$ and $\Delta\tilde{P}^{(\tau_1,\bxi_1,\bxi_1)}(s_1,t;\x',\x,\y)$ are handled similarly, 
following \eqref{regularizing_proxy} and Lemma \ref{THE_CTR_SENSI_DIFF_PARAM}.
\endproof

\section{Proof of the main theorem}\label{SEC_PROOF}
\label{SEC_COMP}

{\color{black}
In this section we assume {\bf (H$^\gamma_\sigma$)} and {\bf (H$^{\gamma}_{\F}$)} for some $\gamma\in(0,1)$.
For $\eps\in(0,1)$, we define
$$
\gF^{(\eps)}(t,\x):=\gF(t,\cdot)*\rho_\eps(\x)=(F^{(\eps)}_1, F^{(\eps)}_2)(t,\x)
$$
and
$$
\sigma^{(\eps)}(s,\x):=\sigma(t,\cdot)*\rho_\eps(\x).
$$
Let $\X^{\eps}_{t,s}(\x)$ be the solution of
\eqref{SDE0} corresponding to $(\F^{(\eps)},\sigma^{(\eps)})$ and $p_{\eps}$ be the corresponding density.
By {Theorem 11.1.4 in \cite{stro:vara:79} and Theorem 1 in \cite{Raynal2017RegularizationEO}} , under \textbf{(H$_{\gF}^{\gamma}$)} and \textbf{(H$_\sigma^{\gamma}$)}, for any $f\in C^\infty_c(\mR^{2d})$,
\begin{align}\label{AAM1}
\lim_{\eps\to 0}\mE f(\X^{\eps}_{t,s}(\x))=\mE f(\X_{t,s}(\x)).
\end{align}
Moreover, from Propositions \ref{Conv_Kernels}, \ref{PROP_CONTR} and Section \ref{LOWER} we have the following uniform estimate: 
there exist constants $\lambda_0,C>0$ depending only on $\Theta_T$ such that for all $\eps\in(0,1)$, 
\begin{align}\label{AAM9}
C^{-1}g_{\lambda^{-1}_0}\(t-s, \btheta^{(\eps)}_{t,s}(\x)-\y\) \le p_{\eps}(s,\x;t,\y)\le C g_{\lambda_0}\(t-s, \btheta^{(\eps)}_{t,s}(\x)-\y\),
\end{align}
where $\btheta^{(\eps)}_{t,s}(\x)$ is the unique solution of the following ODE
$$
{\dot\btheta}^{(\eps)}_{t,s}(\x)=\gF^{(\eps)}(t,\btheta^{(\eps)}_{t,s}(\x)), \quad \btheta^{(\eps)}_{s,s}(\x)=\x.
$$
In order to take limits $\eps\to 0$, we need the following important estimate.
\bl\label{lemme:bilipflow}
For any $T>0$, there exists a constant $C\geq 1$ only depending on $\Theta_T$  such that for all $0\leq s<t\leq T $, $\x\in\mR^{2d}$
and $\eps\leq (t-s)^{3/2}$,
\begin{equation}
\label{EQ2}
|\T_{t-s}^{-1}\big(\btheta^{(\eps)}_{t,s}(\x)-\wt\btheta_{t,s}(\x)\big)| \leq C,
\end{equation}
where $\wt\btheta_{t,s}(\x)$ is defined by ODE \eqref{FLOW_MACRO}.
\el
\begin{proof}
For simplicity of notations, we assume $s=0$ and write for $\x\in\R^{2d}$ and $t\geq 0$,
$$
\ell_i(t):=\big|\(\btheta^{(\eps)}_{t,0}(\x)-\wt\btheta_{t,0}(\x)\)_i\big|,\ \ i=1,2.
$$ 
For $i=1$, noting that
$$
|F_1^{(\eps)}(t,\x)-F_1(t,\x)|\leq \kappa_1(1+\eps),
$$
by definition we have
\begin{align*}
\ell_1(t)&\leq\int^t_0|F^{(\eps)}_1(r,\btheta^{(\eps)}_r(\x))-F^{(1)}_1(r,\wt\btheta_r(\x))|\dif r\\
&\leq\int^t_0|F^{(\eps)}_1(r,\btheta^{(\eps)}_r(\x))-F^{(1)}_1(r,\btheta^{(\eps)}_r(\x))|\dif r\\
&\quad+\int^t_0|F^{(1)}_1(r,\btheta^{(\eps)}_r(\x))-F^{(1)}_1(r,\wt\btheta_r(\x))|\dif r\\
&\leq 2\kappa_1 t+\|\nabla F^{(1)}_1\|_\infty\int^t_0(\ell_1(r)+\ell_2(r))\dif r,
\end{align*}
which implies by Gronwall's inequality
\begin{align}\label{AX7}
\ell_1(t)\lesssim t+\int^t_0\ell_2(r)\dif r.
\end{align}
For $i=2$, note that
$$
|F_2^{(\eps)}(t,\x)-F^{(\eps)}_2(t,\y)|\leq \kappa_2(|x_1-y_1|+|x_2-y_2|^{(1+\gamma)/3}+|x_2-y_2|)
$$
and
$$
|F_2^{(\eps)}(t,\x)-F_2(t,\x)|\leq \kappa_2(\eps^{(1+\gamma)/3}+2\eps).
$$
Below we fix $t\in(0,T]$ and $\eps\leq t^{3/2}$.
By definition we have for all $s\in[0,t]$,
\begin{align*}
\ell_2(s)&\leq\int^s_0|F^{(\eps)}_2(r,\btheta^{(\eps)}_r(\x))-[F_2(r,\cdot)*\rho_{r^{3/2}}](\wt\btheta_r(\x))|\dif r\\
&\leq\int^s_0|F^{(\eps)}_2(r,\btheta^{(\eps)}_r(\x))-F^{(\eps)}_2(r,\wt\btheta_r(\x))|\dif r\\
&\quad+\int^s_0|F^{(\eps)}_2(r,\wt\btheta_r(\x))-[F_2(r,\cdot)*\rho_{r^{3/2}}](\wt\btheta_r(\x))|\dif r\\
&\lesssim\int^s_0\(\ell_1(r)+\ell_2(r)^{\frac{1+\gamma}3}+\ell_2(r)\)\dif r+s t^{(1+\gamma)/2}+\int^s_0r^{(1+\gamma)/2}\dif r\\
&\lesssim t^{(3+\gamma)/2}+\int^s_0\(\ell_2(r)^{\frac{1+\gamma}3}+\ell_2(r)\)\dif r.
\end{align*}
By Lemma \ref{Le11}, we obtain
$$
\sup_{s\in[0,t]}\ell_2(s)\lesssim t^{3/2},
$$
which together with \eqref{AX7} yields \eqref{EQ2}.
\end{proof}

Now by \eqref{AAM9} and \eqref{EQ2}, there is a constant $C_0>0$ such that for any nonnegative $f\in C_b(\mR^{2d})$ and $\eps\leq (t-s)^{3/2}$,
$$
C^{-1}_0\int_{\R^{2d}}g_{\lambda^{-1}_0}\(t-s, \wt\btheta_{t,s}(\x)-\y\)f(\y)\dif \y\leq \mE f(\X^\eps_{t,s}(\x))
\leq C_0\int_{\R^{2d}}g_{\lambda_0}\(t-s, \wt\btheta_{t,s}(\x)-\y\)f(\y)\dif\y,
$$
which together with \eqref{AAM1} yields that
$$
C^{-1}_0\int_{\R^d}g_{\lambda^{-1}_0}\(t-s, \wt\btheta_{t,s}(\x)-\y\)f(\y)\dif \y\leq \mE f(\X_{t,s}(\x))
\leq C_0\int_{\R^d}g_{\lambda_0}\(t-s, \wt\btheta_{t,s}(\x)-\y\)f(\y)\dif\y.
$$
In particular, this implies that $\X_{t,s}(\x)$ has a density $p(s,\x;t,\y)$ having lower and upper bound as in \eqref{Density_bounds_THM}. 
This proves point \textit{(i)} of Theorem \ref{MAIN_THM}.
 
Similarly, we derive from Propositions \ref{gradient_x1}, \ref{gradient_x2} and \ref{gradient_x2_HOLD} that 
under \textbf{(H$_{\gF}^{\gamma}$)} and \textbf{(H$_\sigma^{\gamma}$)}, 
\begin{align}\label{Gr0}
\sup_{\eps\leq(t-s)^{3/2}} \left|\nabla_{x_1}p_{\eps}(s,\x;t,\y) \right|\le {C_1}(t-s)^{-\frac{1}{2}}g_{\lambda_1}\(t-s, \wt\btheta_{t,s}(\x)-\y\), 
\end{align}
under \textbf{(H$_{\gF}^{\gamma}$)}, \textbf{(H$_\sigma^{\gamma}$)} with \eqref{F11},
\begin{align}\label{Gr}
\sup_{\eps\leq(t-s)^{3/2}} \left|\nabla^2_{x_1}p_{\eps}(s,\x;t,\y) \right|&\le {C_2}(t-s)^{-1}g_{\lambda_1}\(t-s, \wt\btheta_{t,s}(\x)-\y\),
\end{align}
and under \textbf{(H$_{\gF}^{\gamma}$)}, \textbf{(H$_\sigma^{\gamma}$)} with \eqref{Re1}, \eqref{Re2},  
\begin{align}\label{Gr2}
\sup_{\eps\leq(t-s)^{3/2}} \left|\nabla_{x_2}p_{\eps}(s,\x;t,\y) \right|&\le {C_3}(t-s)^{-\frac{3}{2}}g_{\lambda_3}\(t-s,\wt\btheta_{t,s}(\x)-\y\),
\end{align}
where in the above equations \eqref{Gr0}-\eqref{Gr2}, the constants $C_1$-$C_3$ only depend on $\Theta_T$.
Uniform H\"older controls in $\x$ also follow from the previously recalled Propositions. Equicontinuity w.r.t. the variable $\y$ could be established similarly.

Then, from the Ascoli-Arzel\`a theorem, one can find a subsequence $\eps_k$ such that for each $\x,\y\in\R^{2d}$,
$$\nabla^j_{x_1} p_{\eps_k}(s,\x;t,\y)\to \nabla^j_{x_1} p(s,\x;t,\y),\ \ j=0,1,2,\ \nabla_{x_2} p_{\eps_k}(s,\x;t,\y)\to \nabla_{x_2} p(s,\x;t,\y).$$
The gradient estimates follow, under the previously recalled additional assumptions when needed, from \eqref{Gr0}-\eqref{Gr2}.
}

\section*{Acknowledgements}

For the second author this research has been partially funded by RSF grant 17-11-01098 for this work which began in May 2021.

\bibliographystyle{acm}
\bibliography{bibli}

\end{document}